\documentclass[reqno,12pt,a4paper]{amsart}

\usepackage{mathrsfs}
\usepackage{amssymb}
\usepackage{amsfonts}
\usepackage{latexsym}
\usepackage{amsthm}
\usepackage{graphicx}

\def\lb{\label}

\newcommand{\er}[1]{\textrm{(\ref{#1})}}

\begin{document}


\renewcommand{\theequation}{\arabic{section}.\arabic{equation}}
\theoremstyle{plain}
\newtheorem{theorem}{\bf Theorem}[section]
\newtheorem{lemma}[theorem]{\bf Lemma}
\newtheorem{corollary}[theorem]{\bf Corollary}
\newtheorem{proposition}[theorem]{\bf Proposition}
\newtheorem{definition}[theorem]{\bf Definition}
\newtheorem{example}[theorem]{\bf Example}

\theoremstyle{remark}
\newtheorem{remark}[theorem]{\bf Remark}
\newtheorem{example1}[theorem]{\bf Example}

\def\a{\alpha}  \def\cA{{\mathcal A}}     \def\bA{{\bf A}}  \def\mA{{\mathscr A}}
\def\b{\beta}   \def\cB{{\mathcal B}}     \def\bB{{\bf B}}  \def\mB{{\mathscr B}}
\def\g{\gamma}  \def\cC{{\mathcal C}}     \def\bC{{\bf C}}  \def\mC{{\mathscr C}}
\def\G{\Gamma}  \def\cD{{\mathcal D}}     \def\bD{{\bf D}}  \def\mD{{\mathscr D}}
\def\d{\delta}  \def\cE{{\mathcal E}}     \def\bE{{\bf E}}  \def\mE{{\mathscr E}}
\def\D{\Delta}  \def\cF{{\mathcal F}}     \def\bF{{\bf F}}  \def\mF{{\mathscr F}}
\def\c{\chi}    \def\cG{{\mathcal G}}     \def\bG{{\bf G}}  \def\mG{{\mathscr G}}
\def\z{\zeta}   \def\cH{{\mathcal H}}     \def\bH{{\bf H}}  \def\mH{{\mathscr H}}
\def\e{\eta}    \def\cI{{\mathcal I}}     \def\bI{{\bf I}}  \def\mI{{\mathscr I}}
\def\p{\psi}    \def\cJ{{\mathcal J}}     \def\bJ{{\bf J}}  \def\mJ{{\mathscr J}}
\def\vT{\Theta} \def\cK{{\mathcal K}}     \def\bK{{\bf K}}  \def\mK{{\mathscr K}}
\def\k{\kappa}  \def\cL{{\mathcal L}}     \def\bL{{\bf L}}  \def\mL{{\mathscr L}}
\def\l{\lambda} \def\cM{{\mathcal M}}     \def\bM{{\bf M}}  \def\mM{{\mathscr M}}
\def\L{\Lambda} \def\cN{{\mathcal N}}     \def\bN{{\bf N}}  \def\mN{{\mathscr N}}
\def\m{\mu}     \def\cO{{\mathcal O}}     \def\bO{{\bf O}}  \def\mO{{\mathscr O}}
\def\n{\nu}     \def\cP{{\mathcal P}}     \def\bP{{\bf P}}  \def\mP{{\mathscr P}}
\def\r{\rho}    \def\cQ{{\mathcal Q}}     \def\bQ{{\bf Q}}  \def\mQ{{\mathscr Q}}
\def\s{\sigma}  \def\cR{{\mathcal R}}     \def\bR{{\bf R}}  \def\mR{{\mathscr R}}
\def\S{\Sigma}  \def\cS{{\mathcal S}}     \def\bS{{\bf S}}  \def\mS{{\mathscr S}}
\def\t{\tau}    \def\cT{{\mathcal T}}     \def\bT{{\bf T}}  \def\mT{{\mathscr T}}
\def\f{\phi}    \def\cU{{\mathcal U}}     \def\bU{{\bf U}}  \def\mU{{\mathscr U}}
\def\F{\Phi}    \def\cV{{\mathcal V}}     \def\bV{{\bf V}}  \def\mV{{\mathscr V}}
\def\P{\Psi}    \def\cW{{\mathcal W}}     \def\bW{{\bf W}}  \def\mW{{\mathscr W}}
\def\o{\omega}  \def\cX{{\mathcal X}}     \def\bX{{\bf X}}  \def\mX{{\mathscr X}}
\def\x{\xi}     \def\cY{{\mathcal Y}}     \def\bY{{\bf Y}}  \def\mY{{\mathscr Y}}
\def\X{\Xi}     \def\cZ{{\mathcal Z}}     \def\bZ{{\bf Z}}  \def\mZ{{\mathscr Z}}

\def\vr{\varrho}
\def\be{{\bf e}} \def\bc{{\bf c}}
\def\bx{{\bf x}} \def\by{{\bf y}}
\def\bv{{\bf v}} \def\bu{{\bf u}}
\def\Om{\Omega} \def\bp{{\bf p}}
\def\bbD{\pmb \Delta}
\def\mm{\mathrm m}
\def\mn{\mathrm n}

\newcommand{\mc}{\mathscr {c}}

\newcommand{\gA}{\mathfrak{A}}          \newcommand{\ga}{\mathfrak{a}}
\newcommand{\gB}{\mathfrak{B}}          \newcommand{\gb}{\mathfrak{b}}
\newcommand{\gC}{\mathfrak{C}}          \newcommand{\gc}{\mathfrak{c}}
\newcommand{\gD}{\mathfrak{D}}          \newcommand{\gd}{\mathfrak{d}}
\newcommand{\gE}{\mathfrak{E}}
\newcommand{\gF}{\mathfrak{F}}           \newcommand{\gf}{\mathfrak{f}}
\newcommand{\gG}{\mathfrak{G}}           
\newcommand{\gH}{\mathfrak{H}}           \newcommand{\gh}{\mathfrak{h}}
\newcommand{\gI}{\mathfrak{I}}           \newcommand{\gi}{\mathfrak{i}}
\newcommand{\gJ}{\mathfrak{J}}           \newcommand{\gj}{\mathfrak{j}}
\newcommand{\gK}{\mathfrak{K}}            \newcommand{\gk}{\mathfrak{k}}
\newcommand{\gL}{\mathfrak{L}}            \newcommand{\gl}{\mathfrak{l}}
\newcommand{\gM}{\mathfrak{M}}            \newcommand{\gm}{\mathfrak{m}}
\newcommand{\gN}{\mathfrak{N}}            \newcommand{\gn}{\mathfrak{n}}
\newcommand{\gO}{\mathfrak{O}}
\newcommand{\gP}{\mathfrak{P}}             \newcommand{\gp}{\mathfrak{p}}
\newcommand{\gQ}{\mathfrak{Q}}             \newcommand{\gq}{\mathfrak{q}}
\newcommand{\gR}{\mathfrak{R}}             \newcommand{\gr}{\mathfrak{r}}
\newcommand{\gS}{\mathfrak{S}}              \newcommand{\gs}{\mathfrak{s}}
\newcommand{\gT}{\mathfrak{T}}             \newcommand{\gt}{\mathfrak{t}}
\newcommand{\gU}{\mathfrak{U}}             \newcommand{\gu}{\mathfrak{u}}
\newcommand{\gV}{\mathfrak{V}}             \newcommand{\gv}{\mathfrak{v}}
\newcommand{\gW}{\mathfrak{W}}             \newcommand{\gw}{\mathfrak{w}}
\newcommand{\gX}{\mathfrak{X}}               \newcommand{\gx}{\mathfrak{x}}
\newcommand{\gY}{\mathfrak{Y}}              \newcommand{\gy}{\mathfrak{y}}
\newcommand{\gZ}{\mathfrak{Z}}             \newcommand{\gz}{\mathfrak{z}}

\def\ve{\varepsilon}   \def\vt{\vartheta}    \def\vp{\varphi}    \def\vk{\varkappa}

\def\A{{\mathbb A}} \def\B{{\mathbb B}} \def\C{{\mathbb C}}
\def\dD{{\mathbb D}} \def\E{{\mathbb E}} \def\dF{{\mathbb F}} \def\dG{{\mathbb G}} \def\H{{\mathbb H}}\def\I{{\mathbb I}} \def\J{{\mathbb J}} \def\K{{\mathbb K}} \def\dL{{\mathbb L}}\def\M{{\mathbb M}} \def\N{{\mathbb N}} \def\O{{\mathbb O}} \def\dP{{\mathbb P}} \def\R{{\mathbb R}}\def\S{{\mathbb S}} \def\T{{\mathbb T}} \def\U{{\mathbb U}} \def\V{{\mathbb V}}\def\W{{\mathbb W}} \def\X{{\mathbb X}} \def\Y{{\mathbb Y}} \def\Z{{\mathbb Z}}


\def\la{\leftarrow}              \def\ra{\rightarrow}            \def\Ra{\Rightarrow}
\def\ua{\uparrow}                \def\da{\downarrow}
\def\lra{\leftrightarrow}        \def\Lra{\Leftrightarrow}


\def\lt{\biggl}                  \def\rt{\biggr}
\def\ol{\overline}               \def\wt{\widetilde}
\def\ul{\underline}
\def\no{\noindent}


\let\ge\geqslant                 \let\le\leqslant
\def\lan{\langle}                \def\ran{\rangle}
\def\/{\over}                    \def\iy{\infty}
\def\sm{\setminus}               \def\es{\emptyset}
\def\ss{\subset}                 \def\ts{\times}
\def\pa{\partial}                \def\os{\oplus}
\def\om{\ominus}                 \def\ev{\equiv}
\def\iint{\int\!\!\!\int}        \def\iintt{\mathop{\int\!\!\int\!\!\dots\!\!\int}\limits}
\def\el2{\ell^{\,2}}             \def\1{1\!\!1}
\def\sh{\sharp}
\def\wh{\widehat}
\def\bs{\backslash}
\def\intl{\int\limits}

\def\na{\mathop{\mathrm{\nabla}}\nolimits}
\def\sh{\mathop{\mathrm{sh}}\nolimits}
\def\ch{\mathop{\mathrm{ch}}\nolimits}
\def\where{\mathop{\mathrm{where}}\nolimits}
\def\all{\mathop{\mathrm{all}}\nolimits}
\def\as{\mathop{\mathrm{as}}\nolimits}
\def\Area{\mathop{\mathrm{Area}}\nolimits}
\def\arg{\mathop{\mathrm{arg}}\nolimits}
\def\const{\mathop{\mathrm{const}}\nolimits}
\def\det{\mathop{\mathrm{det}}\nolimits}
\def\diag{\mathop{\mathrm{diag}}\nolimits}
\def\diam{\mathop{\mathrm{diam}}\nolimits}
\def\dim{\mathop{\mathrm{dim}}\nolimits}
\def\dist{\mathop{\mathrm{dist}}\nolimits}
\def\Im{\mathop{\mathrm{Im}}\nolimits}
\def\Iso{\mathop{\mathrm{Iso}}\nolimits}
\def\Ker{\mathop{\mathrm{Ker}}\nolimits}
\def\Lip{\mathop{\mathrm{Lip}}\nolimits}
\def\rank{\mathop{\mathrm{rank}}\limits}
\def\Ran{\mathop{\mathrm{Ran}}\nolimits}
\def\Re{\mathop{\mathrm{Re}}\nolimits}
\def\Res{\mathop{\mathrm{Res}}\nolimits}
\def\res{\mathop{\mathrm{res}}\limits}
\def\sign{\mathop{\mathrm{sign}}\nolimits}
\def\span{\mathop{\mathrm{span}}\nolimits}
\def\supp{\mathop{\mathrm{supp}}\nolimits}
\def\Tr{\mathop{\mathrm{Tr}}\nolimits}
\def\BBox{\hspace{1mm}\vrule height6pt width5.5pt depth0pt \hspace{6pt}}


\newcommand\nh[2]{\widehat{#1}\vphantom{#1}^{(#2)}}
\def\dia{\diamond}

\def\Oplus{\bigoplus\nolimits}



\def\qqq{\qquad}
\def\qq{\quad}
\let\ge\geqslant
\let\le\leqslant
\let\geq\geqslant
\let\leq\leqslant
\newcommand{\ca}{\begin{cases}}
\newcommand{\ac}{\end{cases}}
\newcommand{\ma}{\begin{pmatrix}}
\newcommand{\am}{\end{pmatrix}}
\renewcommand{\[}{\begin{equation}}
\renewcommand{\]}{\end{equation}}
\def\eq{\begin{equation}}
\def\qe{\end{equation}}
\def\[{\begin{equation}}
\def\bu{\bullet}

\makeatletter
\@namedef{subjclassname@2020}{\textup{2020} Mathematics Subject Classification}
\makeatother

\title[Trace formulas for magnetic Schr\"odinger operators on periodic graphs]
{Trace formulas for magnetic Schr\"odinger operators on periodic graphs and their applications}

\date{\today}
\author[Evgeny Korotyaev]{Evgeny Korotyaev}
\address{Academy for Advanced Interdisciplinary Studies Northeast Normal University, Changchun, 130024 Jilin, China, and HSE University, 3A Kantemirovskaya ulitsa, St. Petersburg,
194100, Russia, \ korotyaev@gmail.com}
\author[Natalia Saburova]{Natalia Saburova}
\address{Northern (Arctic) Federal University, Severnaya Dvina Emb. 17, Arkhangelsk, 163002, Russia,
 \ n.saburova@gmail.com, \ n.saburova@narfu.ru}


\begin{abstract}
We consider Schr\"odinger operators with periodic magnetic and electric potentials on periodic discrete graphs. The spectrum of such operators consists of a finite number of bands. Some of them and even all may be degenerate. We determine trace formulas for the magnetic Schr\"odinger operators. The traces of the fiber operators are expressed as finite Fourier series of the quasimomentum. The coefficients of the Fourier series are given in terms of the magnetic fluxes, electric potentials and cycles in the quotient graph from some specific cycle sets. Using the trace formulas we obtain new lower estimates of the total bandwidth for the magnetic Schr\"odinger operator in terms of geometric parameters of the graph, magnetic fluxes and electric potentials. We show that these estimates are sharp.
\end{abstract}

\maketitle

\section {\lb{Sec1}Introduction}
\setcounter{equation}{0}

Laplace and Schr\"odinger operators on graphs have a lot of
applications in physics,  chemistry and engineering. They are used
to study properties of different periodic media, e.g., nanomedia, see \cite{NG04} (two authors of this survey, Novoselov and Geim, won the Nobel Prize for discovering graphene). A discrete analogue of the magnetic Laplacian on $\R^2$ was introduced by Harper \cite{H55}.
It describes the behavior of an electron moving on the square lattice $\Z^2$ in an external uniform magnetic field perpendicular to the lattice in the so-called tight-binding approximation. Discrete magnetic Laplacians on arbitrary finite and infinite graphs have been studied by many authors (see, e.g., \cite{BGK20,DM06,G14,GKS16,LL93,Sh94,S94} and references therein).

Although the magnetic field is an external structure to graph, the discrete magnetic Laplacian can also be used to study some structural properties of the underlying graphs. For example, the discrete magnetic Laplacians with specific magnetic potentials were used in \cite{MSS15} for construction of Ramanujan expanders playing an important role in computer science. In \cite{FLP22} the authors related the spectrum of the discrete magnetic Laplacian on a finite graph with two structural properties of the graph: the existence of a perfect matching and the existence of a Hamiltonian cycle of the underlying graph. The Cheeger constant for discrete magnetic Laplacian was discussed in \cite{FLP20,LLPP15}. Moreover, the discrete Laplacian on periodic graphs without magnetic fields can be reduced to a family of discrete magnetic Laplacians on the finite quotient graph, where the magnetic potential is interpreted as a Floquet parameter of the periodic graph, see \cite{FLP18} and Remark \ref{RFD}.\emph{iii}.

In this paper we consider discrete Schr\"odinger operators with periodic magnetic and electric potentials on periodic graphs. Such operators were investigated in \cite{FL19,FLP20,HS99a,HS99b,KS17}. Their spectrum consists of an absolutely continuous part (a~union of a finite number of non-degenerate bands) and a finite number of flat bands, i.e., eigenvalues of infinite multiplicity. Higuchi-Shirai \cite{HS99a} proved the analyticity of the bottom of the spectrum near zero magnetic potential and represented the second derivative of the bottom of the spectrum in terms of geometry of a graph. The authors showed in \cite{HS99b} that if the Betti number of the quotient graph is equal to the dimension of the periodic one, then the magnetic Laplacian on the periodic graph is unitarily equivalent to the Laplacian without magnetic fields. Korotyaev-Saburova \cite{KS17} estimated the Lebesgue measure of the spectrum in terms of Betti numbers. They also obtained estimates of a variation of the spectrum of the Schr\"odinger operators under a perturbation by a magnetic field in terms of magnetic fluxes and estimates of effective masses associated with the ends of each spectral band for magnetic Laplacians in terms of geometric parameters of the graphs. Fabila-Carrasco-Lled\'o-Post \cite{FL19,FLP20} developed a band localization for the magnetic Laplacian on periodic graphs and gave sufficient conditions for the existence of gaps in its spectrum.

\medskip

This paper is a continuation of our previous works \cite{KS22,KS22a} devoted to the trace formulas for the Schr\"odinger operator \emph{without magnetic fields} and the corresponding estimates of its total bandwidth. Spectral analysis of the \emph{magnetic} Schr\"odinger operators is much more complicated problem. For example, in contrast to the non-magnetic case, the spectrum of the magnetic Schr\"odinger operator may consist of only flat bands, i.e., the absolutely continuous spectrum may be empty (see Example \ref{Exa}.\emph{ii}). Thus, the problem of estimating the total bandwidth of the magnetic Schr\"odinger operator \emph{from below} is more delicate than the similar problem for the non-magnetic case.

Our main goals are the following:

\begin{itemize}
  \item to determine \emph{trace formulas} connecting the spectra of the magnetic Schr\"odinger operators on periodic graphs with geometric parameters of the graphs, magnetic fluxes and electric potentials;
\end{itemize}
and, as applications of these trace formulas,
\begin{itemize}
  \item to obtain \emph{lower} bounds for the total bandwidth of the magnetic Schr\"odinger operators on periodic graphs;
  \item to formulate sufficient conditions on the magnetic fluxes and geometry of the periodic graph $\cG$ under which the spectrum of the magnetic Schr\"odinger operator on $\cG$ has at least one non-degenerate band.
\end{itemize}

The problem of estimating the total bandwidth for the Schr\"odinger operators on periodic graphs was studied by the authors during last 10 years \cite{KS14,KS17,KS18,KS20,KS22a}. First, in \cite{KS14} we obtained the \emph{upper} bound for the non-magnetic case in terms of the Betti number of the quotient graph. Then this upper estimate was generalized to the magnetic Schr\"odinger operators \cite{KS17}. Later these upper estimates were improved in terms of some geometric invariant of the periodic graph \cite{KS18,KS20}.
Finally, we obtained \emph{lower} estimates for the total bandwidth of the Schr\"odinger operators without magnetic fields \cite{KS22a}. It became possible due to the trace formulas for the non-magnetic case presented in our recent article \cite{KS22}. In the present paper we obtain \emph{lower} bounds for the total bandwidth of the magnetic Schr\"odinger operators. These lower bounds are new and we do not know any other lower estimates for the total bandwidth of the magnetic Schr\"odinger operators on arbitrary periodic graphs. If the magnetic field is absent, then the obtained estimates coincide with the estimates for the total bandwidth of the Schr\"odinger operator without magnetic fields \cite{KS22a}, and in the case of the one-dimensional lattice $\Z$ they coincide with the known Last's estimate (see \cite{L92}), which is sharp. For more details see Section \ref{S2.3}.

\subsection{Periodic graphs.} \lb{Sec.1}
Let $\cG=(\cV,\cE)$ be a connected infinite graph, possibly  having
loops and multiple edges and embedded into the space $\R^d$. Here
$\cV$  is the set of its vertices and $\cE$ is the set of its
unoriented edges. If a vertex is an endpoint of an edge, they are  \emph{incident}. Considering each edge in $\cE$ to have two orientations, we introduce the set $\cA$ of all oriented edges. An
edge starting at a vertex $x$ and ending at a vertex $y$ from $\cV$
will be denoted as the ordered pair $(x,y)\in\cA$. Let $\ul\be=(y,x)$ be the inverse edge of $\be=(x,y)\in\cA$. We define the \emph{degree} ${\vk}_x$ of the vertex $x\in\cV$ as the number of all unoriented edges from $\cE$ incident to the vertex $x$ (with loops counted twice),  or, equivalently, as the number of all oriented edges from $\cA$ starting at $x$.

Let $\G$ be a lattice of rank $d$ in $\R^d$ with a basis $\{\ga_1,\ldots,\ga_d\}$, i.e.,
$$
\G=\Big\{\ga : \ga=\sum_{s=1}^dn_s\ga_s, \; (n_s)_{s=1}^d\in\Z^d\Big\},
$$
and let
\[\lb{fuce}
\Omega=\Big\{\bx\in\R^d : \bx=\sum_{s=1}^d x_s\ga_s, \; (x_s)_{s=1}^d\in
[0,1)^d\Big\}
\]
be the \emph{fundamental cell} of the lattice $\G$. We define an equivalence relation on $\R^d$ by:
$$
\bx\equiv \by \; (\hspace{-4mm}\mod \G) \qq\Leftrightarrow\qq \bx-\by\in\G \qqq \forall\, \bx,\by\in\R^d.
$$
Identifying equivalent points of $\R^d$, we obtain the $d$-dimensional torus $\R^d/\G$.

We consider \emph{locally finite $\G$-periodic graphs} $\cG$, i.e.,
graphs satisfying the following conditions:
\begin{itemize}
  \item $\cG$ is invariant under translation by any vector $\ga\in\G$;
  \item the quotient graph  $\cG_*=\cG/\G$ is finite.
\end{itemize}
The basis $\ga_1,\ldots,\ga_d$ of the lattice $\G$ is called the {\it
periods}  of $\cG$. We also call the quotient graph $\cG_*=\cG/\G$
the \emph{fundamental graph} of the periodic graph $\cG$. The
fundamental graph $\cG_*$ is a graph on the $d$-dimensional torus
$\R^d/\G$. The graph $\cG_*=(\cV_*,\cE_*)$ has the vertex set
$\cV_*=\cV/\G$, the set $\cE_*=\cE/\G$ of unoriented edges and the
set $\cA_*=\cA/\G$ of oriented edges which are finite.

\setlength{\unitlength}{1.0mm}
\begin{figure}
\centering
\unitlength 1mm 
\linethickness{0.4pt}
\ifx\plotpoint\undefined\newsavebox{\plotpoint}\fi 
\begin{picture}(120,42)(0,0)

\put(25,25.0){$\Omega$}
\put(7,10){\line(1,0){36.00}}
\put(7,20){\line(1,0){36.00}}
\put(7,30){\line(1,0){36.00}}
\put(7,40){\line(1,0){36.00}}
\put(10,7){\line(0,1){36.00}}
\put(20,7){\line(0,1){36.00}}
\put(30,7){\line(0,1){36.00}}
\put(40,7){\line(0,1){36.00}}

\bezier{15}(20.5,20)(20.5,25)(20.5,30)
\bezier{15}(21,20)(21,25)(21,30)
\bezier{15}(21.5,20)(21.5,25)(21.5,30)
\bezier{15}(22,20)(22,25)(22,30)
\bezier{15}(22.5,20)(22.5,25)(22.5,30)
\bezier{15}(23,20)(23,25)(23,30)
\bezier{15}(23.5,20)(23.5,25)(23.5,30)
\bezier{15}(24,20)(24,25)(24,30)
\bezier{15}(24.5,20)(24.5,25)(24.5,30)
\bezier{15}(25,20)(25,25)(25,30)
\bezier{15}(25.5,20)(25.5,25)(25.5,30)
\bezier{15}(26,20)(26,25)(26,30)
\bezier{15}(26.5,20)(26.5,25)(26.5,30)
\bezier{15}(27,20)(27,25)(27,30)
\bezier{15}(27.5,20)(27.5,25)(27.5,30)
\bezier{15}(28,20)(28,25)(28,30)
\bezier{15}(28.5,20)(28.5,25)(28.5,30)
\bezier{15}(29,20)(29,25)(29,30)
\bezier{15}(29.5,20)(29.5,25)(29.5,30)

\put(10,10){\circle{1}}
\put(20,10){\circle{1}}
\put(30,10){\circle{1}}
\put(40,10){\circle{1}}

\put(10,20){\circle{1}}
\put(20,20){\circle*{1}}
\put(30,20){\circle{1}}
\put(40,20){\circle{1}}

\put(10,30){\circle{1}}
\put(20,30){\circle{1}}
\put(30,30){\circle{1}}
\put(40,30){\circle{1}}

\put(10,40){\circle{1}}
\put(20,40){\circle{1}}
\put(30,40){\circle{1}}
\put(40,40){\circle{1}}

\put(20,20){\vector(1,0){10.00}}
\put(20,20){\vector(0,1){10.00}}
\put(5,42){$\bS$}
\put(15.7,17.5){$\scriptstyle O=x$}
\put(26,31){$\scriptstyle x+\ga_1+\ga_2$}
\put(13,31){$\scriptstyle x+\ga_2$}
\put(31,18){$\scriptstyle x+\ga_1$}
\put(24,18){$\scriptstyle  \ga_1$}
\put(16.5,25){$\scriptstyle \ga_2$}
\put(24,20.5){$\scriptstyle  \be_1$}
\put(20.3,25){$\scriptstyle \be_2$}
\put(0,5){\emph{a})}
\put(46,5){\emph{b})}
\put(70,40){$\bS_*$}
\put(55,10){\vector(1,0){20.00}}
\put(55,10){\vector(0,1){20.00}}

\bezier{30}(55.5,10)(55.5,20)(55.5,30)
\bezier{30}(56,10)(56,20)(56,30)
\bezier{30}(56.5,10)(56.5,20)(56.5,30)
\bezier{30}(57,10)(57,20)(57,30)
\bezier{30}(57.5,10)(57.5,20)(57.5,30)
\bezier{30}(58,10)(58,20)(58,30)
\bezier{30}(58.5,10)(58.5,20)(58.5,30)
\bezier{30}(59,10)(59,20)(59,30)
\bezier{30}(59.5,10)(59.5,20)(59.5,30)
\bezier{30}(60,10)(60,20)(60,30)
\bezier{30}(60.5,10)(60.5,20)(60.5,30)
\bezier{30}(61,10)(61,20)(61,30)
\bezier{30}(61.5,10)(61.5,20)(61.5,30)
\bezier{30}(62,10)(62,20)(62,30)
\bezier{30}(62.5,10)(62.5,20)(62.5,30)
\bezier{30}(63,10)(63,20)(63,30)
\bezier{30}(63.5,10)(63.5,20)(63.5,30)
\bezier{30}(64,10)(64,20)(64,30)
\bezier{30}(64.5,10)(64.5,20)(64.5,30)
\bezier{30}(65,10)(65,20)(65,30)
\bezier{30}(65.5,10)(65.5,20)(65.5,30)
\bezier{30}(66,10)(66,20)(66,30)
\bezier{30}(66.5,10)(66.5,20)(66.5,30)
\bezier{30}(67,10)(67,20)(67,30)
\bezier{30}(67.5,10)(67.5,20)(67.5,30)
\bezier{30}(68,10)(68,20)(68,30)
\bezier{30}(68.5,10)(68.5,20)(68.5,30)
\bezier{30}(69,10)(69,20)(69,30)
\bezier{30}(69.5,10)(69.5,20)(69.5,30)
\bezier{30}(70,10)(70,20)(70,30)
\bezier{30}(70.5,10)(70.5,20)(70.5,30)
\bezier{30}(71,10)(71,20)(71,30)
\bezier{30}(71.5,10)(71.5,20)(71.5,30)
\bezier{30}(72,10)(72,20)(72,30)
\bezier{30}(72.5,10)(72.5,20)(72.5,30)
\bezier{30}(73,10)(73,20)(73,30)
\bezier{30}(73.5,10)(73.5,20)(73.5,30)
\bezier{30}(74,10)(74,20)(74,30)
\bezier{30}(74.5,10)(74.5,20)(74.5,30)
\bezier{30}(75,10)(75,20)(75,30)

\put(55,10){\circle*{1}}
\put(55,30){\circle{1}}
\put(75,30){\circle{1}}
\put(75,10){\circle{1}}

\multiput(55,30)(4,0){5}{\line(1,0){2}}
\multiput(75,10)(0,4){5}{\line(0,1){2}}

\put(52.0,9.0){$x$}
\put(76.5,9.0){$x$}
\put(52.5,31.0){$x$}
\put(75.5,31.0){$x$}
\put(63.5,7.5){$\ga_1$}
\put(51.0,20.0){$\ga_2$}
\put(63.5,11.5){$\be_1$}
\put(55.5,20.0){$\be_2$}

\put(64,22){$\Omega$}

\bezier{120}(85,10)(95,5)(105,10)
\bezier{90}(85,10)(76,15)(85,20)
\bezier{90}(105,10)(114,15)(105,20)
\bezier{120}(85,20)(95,25)(105,20)

\bezier{120}(88,15)(95,10)(102,15)
\bezier{120}(90,14)(95,18)(100,14)
\put(93.6,10){\circle*{1}}
\bezier{15}(95,7.5)(98,10)(95,12.5)
\bezier{70}(95,7.5)(92,10)(95,12.5)

\bezier{120}(86,12.5)(95,7.5)(104,12.5)
\bezier{90}(86,12.5)(82,15)(86,17.5)
\bezier{90}(104,12.5)(108,15)(104,17.5)
\bezier{120}(86,17.5)(95,22.5)(104,17.5)
\put(91.5,8.3){$\scriptstyle x$}
\put(103.0,14.7){$\scriptstyle \be_1$}
\put(96.5,10.9){$\scriptstyle \be_2$}

\put(93.0,30){\circle*{1}}
\put(90,28){$x$}
\put(102,29){$\be_1$}
\put(91,40){$\be_2$}
\bezier{200}(93.0,30)(106,37)(107.0,30)
\bezier{200}(93.0,30)(106,23)(107.0,30)
\bezier{200}(93.0,30)(100,43)(93.0,44)
\bezier{200}(93.0,30)(86,43)(93.0,44)
\put(108,29){$\scriptstyle (1,0)$}
\put(83,37){$\scriptstyle (0,1)$}
\end{picture}
\vspace{-5mm}
\caption{\scriptsize\emph{a}) The square lattice $\bS$; $\ga_1,\ga_2$ are the periods of $\bS$; the fundamental cell $\Omega=[0,1)^2$ is shaded; \; \emph{b})~the fundamental graph $\bS_\ast$ of $\bS$; the indices are shown near the edges.}
\label{FSqL}
\end{figure}
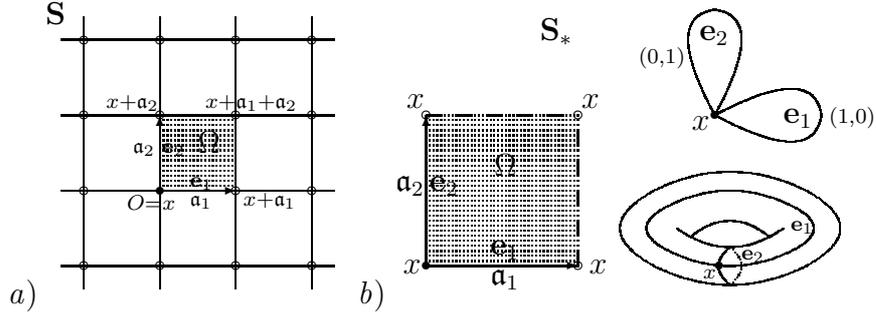

\begin{example1}\lb{ExPG} We consider some examples of periodic graphs.

\emph{i}) The square lattice $\bS=(\cV,\cE)$, where the vertex set $\cV$ and the edge set $\cE$ are given by
$$
\cV=\Z^2,\qq
\cE=\big\{(x,x+\ga_s), \;
\forall\,x\in\Z^2, \; s=1,2\big\}, \qq \ga_1=(1,0), \qq \ga_2=(0,1),
$$
see Fig.\ref{FSqL}\emph{a}. The square lattice $\bS$ is a $\Z^2$-periodic graph with periods $\ga_1,\ga_2$. The fundamental cell $\Omega=[0,1)^2$ is shaded in the figure. Identifying both pairs of opposite sides of $\Omega$ we get the two-dimensional torus $\R^2/\Z^2$, and the fragment of the square lattice $\bS$ in $\Omega$ becomes the fundamental graph $\bS_*$ on the torus $\R^2/\Z^2$, Fig.\ref{FSqL}\emph{b}. The fundamental graph $\bS_*$ of the square lattice $\bS$ consists of one vertex $x$ with degree $\vk_x=4$ and two loop edges $\be_1,\be_2$ at this vertex $x$.

\emph{ii}) The hexagonal lattice $\bG=(\cV,\cE)$ is shown in Fig.\ref{ff.0.3}\emph{a}, the vectors $\ga_1,\ga_2$ are the periods of $\bG$. The vertex set $\cV$ and the edge set $\cE$ over the basis $\ga_1,\ga_2$ are given by
$$
\textstyle \cV=\Z^2\cup\big(\Z^2+\big(\frac13\,,\frac13\big)\big),
$$
$$
\textstyle \cE=\big\{\big(x,x+\big(\frac13\,,\frac13\big)\big),
\big(x,x+\big(-\frac23\,,\frac13\big)\big),
\big(x,x+\big(\frac13\,,-\frac23\big)\big),\quad\forall\,x\in\Z^2\big\}.
$$
The fundamental cell $\Omega$ is the parallelogram spanned by the periods $\ga_1,\ga_2$. The fundamental graph $\bG_*$ of $\bG$ consists of two vertices $x_1,x_2$ and three multiple edges $\be_1,\be_2,\be_3$ connecting these vertices (Fig.\ref{ff.0.3}\emph{b}).

\emph{iii}) The periodic graph $\bK=(\cV,\cE)$ shown in Fig.\ref{FKaL}\emph{a} is called the \emph{Kagome lattice}. It is invariant under translation by vectors $\ga_1,\ga_2$. Its fundamental graph $\bK_*$ consists of three vertices $x_1,x_2,x_3$ and six edges $\be_1,\ldots,\be_6$ (Fig.~\ref{FKaL}\emph{b}).
\end{example1}

\begin{figure}
\centering
\unitlength 1mm 
\linethickness{0.4pt}
\ifx\plotpoint\undefined\newsavebox{\plotpoint}\fi 
\begin{picture}(125,43)(0,0)
\put(5,5){\emph{a})}
\bezier{20}(24.5,16.3)(24.5,22.3)(24.5,28.3)
\bezier{20}(25.0,16.6)(25.0,22.6)(25.0,28.6)
\bezier{20}(25.5,16.9)(25.5,22.9)(25.5,28.9)
\bezier{20}(26.0,17.2)(26.0,23.2)(26.0,29.2)
\bezier{20}(26.5,17.5)(26.5,23.5)(26.5,29.5)
\bezier{20}(27.0,17.8)(27.0,23.8)(27.0,29.8)
\bezier{20}(27.5,18.1)(27.5,24.1)(27.5,30.1)
\bezier{20}(28.0,18.4)(28.0,24.4)(28.0,30.4)
\bezier{20}(28.5,18.7)(28.5,24.7)(28.5,30.7)
\bezier{20}(29.0,19.0)(29.0,25.0)(29.0,31.0)

\bezier{20}(29.5,19.3)(29.5,25.3)(29.5,31.3)
\bezier{20}(30.0,19.6)(30.0,25.6)(30.0,31.6)
\bezier{20}(30.5,19.9)(30.5,25.9)(30.5,31.9)
\bezier{20}(31.0,20.2)(31.0,26.2)(31.0,32.2)
\bezier{20}(31.5,20.5)(31.5,26.5)(31.5,32.5)
\bezier{20}(32.0,20.8)(32.0,26.8)(32.0,32.8)
\bezier{20}(32.5,21.1)(32.5,27.1)(32.5,33.1)
\bezier{20}(33.0,21.4)(33.0,27.4)(33.0,33.4)
\bezier{20}(33.5,21.7)(33.5,27.7)(33.5,33.7)
\bezier{20}(34,22)(34,28)(34,34)

\put(17.5,29.5){$\scriptstyle x_2+\ga_2$}

\put(14,10){\circle{1}}
\put(28,10){\circle{1}}
\put(34,10){\circle{1}}
\put(48,10){\circle{1}}

\put(18,16){\circle{1}}
\put(24,16){\circle*{1}}
\put(38,16){\circle{1}}
\put(44,16){\circle{1}}

\put(14,22){\circle{1}}
\put(28,22){\circle*{1}}
\put(34,22){\circle{1}}
\put(48,22){\circle{1}}

\put(18,28){\circle{1}}
\put(24,28){\circle{1}}
\put(38,28){\circle{1}}
\put(44,28){\circle{1}}

\put(14,34){\circle{1}}
\put(28,34){\circle{1}}
\put(34,34){\circle{1}}
\put(48,34){\circle{1}}

\put(18,40){\circle{1}}
\put(24,40){\circle{1}}
\put(38,40){\circle{1}}
\put(44,40){\circle{1}}

\put(28,10){\line(1,0){6.00}}
\put(18,16){\line(1,0){6.00}}
\put(38,16){\line(1,0){6.00}}

\put(28,22){\line(1,0){6.00}}
\put(18,28){\line(1,0){6.00}}
\put(38,28){\line(1,0){6.00}}

\put(28,34){\line(1,0){6.00}}
\put(18,40){\line(1,0){6.00}}
\put(38,40){\line(1,0){6.00}}

\put(14,10){\line(2,3){4.00}}

\put(34,10){\line(2,3){4.00}}
\put(24,16){\line(2,3){4.00}}
\put(44,16){\line(2,3){4.00}}

\put(14,22){\line(2,3){4.00}}
\put(34,22){\line(2,3){4.00}}
\put(24,28){\line(2,3){4.00}}
\put(44,28){\line(2,3){4.00}}

\put(14,34){\line(2,3){4.00}}
\put(34,34){\line(2,3){4.00}}

\put(28,10){\line(-2,3){4.00}}
\put(48,10){\line(-2,3){4.00}}
\put(38,16){\line(-2,3){4.00}}
\put(18,16){\line(-2,3){4.00}}

\put(28,22){\line(-2,3){4.00}}
\put(48,22){\line(-2,3){4.00}}
\put(38,28){\line(-2,3){4.00}}
\put(18,28){\line(-2,3){4.00}}

\put(28,34){\line(-2,3){4.00}}
\put(48,34){\line(-2,3){4.00}}

\put(30,18){$\scriptstyle \ga_1$}
\put(20.5,22){$\scriptstyle \ga_2$}

\bezier{20}(24,28)(29,31)(34,34)
\bezier{20}(34,22)(34,28)(34,34)

\put(24,16){\vector(0,1){12.0}}
\put(33,21.3){\vector(3,2){0.5}}

\qbezier(24,16)(29,19)(34,22)

\put(35,21.5){$\scriptstyle x_2+\ga_1$}
\put(30,35){$\scriptstyle x_2+\ga_1+\ga_2$}

\put(21.5,14.0){$\scriptstyle x_2$}

\put(30,28){$\Omega$}
\put(32,22){\vector(1,0){0.5}}
\put(26.4,19.6){\vector(2,3){0.5}}
\put(26,25){\vector(-2,3){0.5}}

\put(31,23){$\scriptstyle \mathbf{e}_2$}
\put(25.5,26){$\scriptstyle \mathbf{e}_3$}
\put(24.0,20.5){$\scriptstyle \mathbf{e}_1$}
\put(27.5,23){$\scriptstyle x_1$}

\put(14,10){\line(-1,0){3.00}}
\put(14,22){\line(-1,0){3.00}}
\put(14,34){\line(-1,0){3.00}}
\put(48,10){\line(1,0){3.00}}
\put(48,22){\line(1,0){3.00}}
\put(48,34){\line(1,0){3.00}}

\bezier{50}(14,10)(15,8.5)(16,7)
\bezier{50}(34,10)(35,8.5)(36,7)
\bezier{50}(28,10)(27,8.5)(26,7)
\bezier{50}(48,10)(47,8.5)(46,7)

\bezier{50}(18,40)(17,41.5)(16,43)
\bezier{50}(38,40)(37,41.5)(36,43)
\bezier{50}(24,40)(25,41.5)(26,43)
\bezier{50}(44,40)(45,41.5)(46,43)

\put(65,10){\circle*{1}}
\put(72,20){\circle*{1}}

\put(65,30){\circle{1}}
\put(85,40){\circle{1}}
\put(85,20){\circle{1}}

\put(65,10){\vector(0,1){20.0}}
\put(65,10){\vector(2,1){20.0}}

\multiput(85,20)(0,7){3}{\line(0,1){4}}
\put(65,30){\line(2,1){4.0}}
\put(72,33.5){\line(2,1){4.0}}
\put(79,37){\line(2,1){4.0}}

\qbezier(72,20)(78.5,20.0)(85,20)
\qbezier(72,20)(68.5,15.0)(65,10)
\qbezier(72,20)(68.5,25.0)(65,30)

\put(79,20){\vector(1,0){1.0}}
\put(68.6,15){\vector(3,4){1.0}}
\put(69,24.3){\vector(-3,4){1.0}}

\put(61,8.0){$x_2$}
\put(72.5,21){$x_1$}
\put(86,19.0){$x_2$}
\put(60.5,31.0){$x_2$}
\put(83.5,42.0){$x_2$}
\put(75,13){$\ga_1$}
\put(60,20.0){$\ga_2$}

\bezier{40}(65.5,10.25)(65.5,20.25)(65.5,30.25)
\bezier{40}(66.0,10.50)(66.0,20.50)(66.0,30.50)
\bezier{40}(66.5,10.75)(66.5,20.75)(66.5,30.75)
\bezier{40}(67.0,11.00)(67.0,21.00)(67.0,31.00)
\bezier{40}(67.5,11.25)(67.5,21.25)(67.5,31.25)
\bezier{40}(68.0,11.50)(68.0,21.50)(68.0,31.50)
\bezier{40}(68.5,11.75)(68.5,21.75)(68.5,31.75)
\bezier{40}(69.0,12.00)(69.0,22.00)(69.0,32.00)
\bezier{40}(69.5,12.25)(69.5,22.25)(69.5,32.25)
\bezier{40}(70.0,12.50)(70.0,22.50)(70.0,32.50)
\bezier{40}(70.5,12.75)(70.5,22.75)(70.5,32.75)
\bezier{40}(71.0,13.00)(71.0,23.00)(71.0,33.00)
\bezier{40}(71.5,13.25)(71.5,23.25)(71.5,33.25)
\bezier{40}(72.0,13.50)(72.0,23.50)(72.0,33.50)
\bezier{40}(72.5,13.75)(72.5,23.75)(72.5,33.75)
\bezier{40}(73.0,14.00)(73.0,24.00)(73.0,34.00)
\bezier{40}(73.5,14.25)(73.5,24.25)(73.5,34.25)
\bezier{40}(74.0,14.50)(74.0,24.50)(74.0,34.50)
\bezier{40}(74.5,14.75)(74.5,24.75)(74.5,34.75)
\bezier{40}(75.0,15.00)(75.0,25.00)(75.0,35.00)
\bezier{40}(75.5,15.25)(75.5,25.25)(75.5,35.25)
\bezier{40}(76.0,15.50)(76.0,25.50)(76.0,35.50)
\bezier{40}(76.5,15.75)(76.5,25.75)(76.5,35.75)
\bezier{40}(77.0,16.00)(77.0,26.00)(77.0,36.00)
\bezier{40}(77.5,16.25)(77.5,26.25)(77.5,36.25)
\bezier{40}(78.0,16.50)(78.0,26.50)(78.0,36.50)
\bezier{40}(78.5,16.75)(78.5,26.75)(78.5,36.75)
\bezier{40}(79.0,17.00)(79.0,27.00)(79.0,37.00)
\bezier{40}(79.5,17.25)(79.5,27.25)(79.5,37.25)
\bezier{40}(80.0,17.50)(80.0,27.50)(80.0,37.50)
\bezier{40}(80.5,17.75)(80.5,27.75)(80.5,37.75)
\bezier{40}(81.0,18.00)(81.0,28.00)(81.0,38.00)
\bezier{40}(81.5,18.25)(81.5,28.25)(81.5,38.25)
\bezier{40}(82.0,18.50)(82.0,28.50)(82.0,38.50)
\bezier{40}(82.5,18.75)(82.5,28.75)(82.5,38.75)
\bezier{40}(83.0,19.00)(83.0,29.00)(83.0,39.00)
\bezier{40}(83.5,19.25)(83.5,29.25)(83.5,39.25)
\bezier{40}(84.0,19.50)(84.0,29.50)(84.0,39.50)
\bezier{40}(84.5,19.75)(84.5,29.75)(84.5,39.75)
\bezier{40}(85.0,20.00)(85.0,30.00)(85.0,40.00)
\put(54,5){\emph{b})}
\put(7,38){$\bG$}
\put(65,38){$\bG_*$}
\put(74.3,28.5){$\Omega$}

\put(66,17.2){$\be_1$}
\put(79,21.2){$\be_2$}
\put(68,26.6){$\be_3$}

\put(105,10){\circle*{1}}
\put(117,26){\circle*{1}}
\put(102,7.0){$x_2$}
\put(117,27){$x_1$}
\bezier{200}(105,10)(111,18)(117,26)
\bezier{200}(105,10)(102,27)(117,26)
\bezier{200}(105,10)(121,10)(117,26)
\put(111.0,18.0){\vector(2,3){1.0}}
\put(107.0,23.0){\vector(-1,-1){1.0}}
\put(117.7,17.0){\vector(-1,-2){1.0}}
\put(116,12){$\scriptstyle(1,0)$}
\put(101,24.0){$\scriptstyle(0,1)$}
\put(109.5,15){$\scriptstyle(0,0)$}

\put(118.5,15){$\be_2$}
\put(101.5,21){$\be_3$}
\put(107.0,18.5){$\be_1$}

\end{picture}
\vspace{-0.2cm}
\caption{\footnotesize \emph{a}) The hexagonal lattice $\bG$; $\ga_1,\ga_2$ are the periods of $\bG$; the fundamental cell $\Omega$ is shaded; \; \emph{b})~the fundamental graph $\bG_*$; the indices are shown near the corresponding edges.} \label{ff.0.3}
\end{figure}
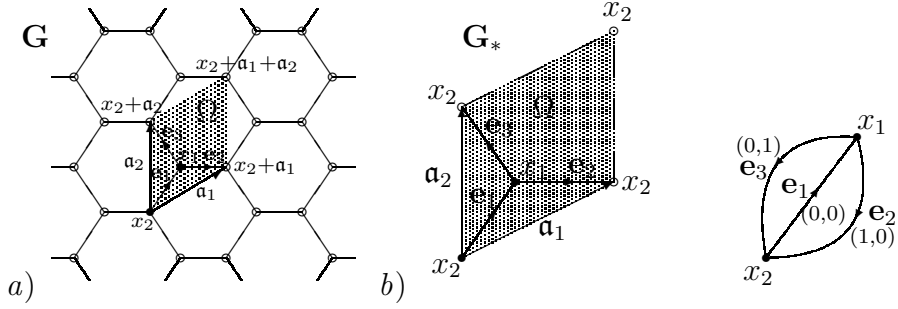

\setlength{\unitlength}{1.0mm}
\begin{figure}[h]
\centering
\unitlength 1mm 
\linethickness{0.4pt}
\ifx\plotpoint\undefined\newsavebox{\plotpoint}\fi 

\begin{picture}(70,55)(0,0)
\bezier{25}(21,15)(26,25)(31,35)
\bezier{25}(22,15)(27,25)(32,35)
\bezier{25}(23,15)(28,25)(33,35)
\bezier{25}(24,15)(29,25)(34,35)
\bezier{25}(25,15)(30,25)(35,35)
\bezier{25}(26,15)(31,25)(36,35)
\bezier{25}(27,15)(32,25)(37,35)
\bezier{25}(28,15)(33,25)(38,35)
\bezier{25}(29,15)(34,25)(39,35)
\bezier{25}(30,15)(35,25)(40,35)
\bezier{25}(31,15)(36,25)(41,35)
\bezier{25}(32,15)(37,25)(42,35)
\bezier{25}(33,15)(38,25)(43,35)
\bezier{25}(34,15)(39,25)(44,35)
\bezier{25}(35,15)(40,25)(45,35)
\bezier{25}(36,15)(41,25)(46,35)
\bezier{25}(37,15)(42,25)(47,35)
\bezier{25}(38,15)(43,25)(48,35)
\bezier{25}(39,15)(44,25)(49,35)

\put(10.0,51){$\bK$}
\put(34.0,25){$\Omega$}
\put(9.0,20){$\scriptstyle(-1,0)$}
\put(32.0,20){$\scriptstyle(0,0)$}
\put(52.0,20){$\scriptstyle(1,0)$}

\put(21.0,44){$\scriptstyle(-1,1)$}
\put(42.0,44){$\scriptstyle(0,1)$}
\put(62.0,44){$\scriptstyle(1,1)$}

\put(-3,15){\line(1,0){66.0}}
\put(7,35){\line(1,0){66.0}}
\put(17,55){\line(1,0){66.0}}
\put(-1.5,12){\line(1,2){23.0}}
\put(18.5,12){\line(1,2){23.0}}
\put(38.5,12){\line(1,2){23.0}}
\put(58.5,12){\line(1,2){23.0}}

\put(11.5,12){\line(-1,2){8.0}}
\put(31.5,12){\line(-1,2){18.0}}
\put(51.5,12){\line(-1,2){23.0}}
\put(66.5,22){\line(-1,2){18.0}}
\put(76.5,42){\line(-1,2){8.0}}

\put(20,15){\vector(1,0){20.0}}
\put(20,15){\vector(1,2){10.0}}

\put(30,15){\vector(-1,0){10.0}}
\put(20,15){\vector(1,2){5.0}}
\put(25,25){\vector(1,-2){5.0}}
\put(24,18){$\scriptstyle\f_1$}

\put(50,35){\vector(-1,0){10.0}}
\put(40,35){\vector(1,-2){5.0}}
\put(45,25){\vector(1,2){5.0}}
\put(43.5,30){$\scriptstyle\f_2$}

\put(0,15){\circle{1}}
\put(10,15){\circle{1}}
\put(20,15){\circle*{1.5}}
\put(30,15){\circle*{1.5}}
\put(40,15){\circle{1}}
\put(50,15){\circle{1}}
\put(60,15){\circle{1}}

\put(5,25){\circle{1}}
\put(25,25){\circle*{1.5}}
\put(45,25){\circle{1}}
\put(65,25){\circle{1}}

\put(10,35){\circle{1}}
\put(20,35){\circle{1}}
\put(30,35){\circle{1}}
\put(40,35){\circle{1}}
\put(50,35){\circle{1}}
\put(60,35){\circle{1}}
\put(70,35){\circle{1}}

\put(15,45){\circle{1}}
\put(35,45){\circle{1}}
\put(55,45){\circle{1}}
\put(75,45){\circle{1}}

\put(20,55){\circle{1}}
\put(30,55){\circle{1}}
\put(40,55){\circle{1}}
\put(50,55){\circle{1}}
\put(60,55){\circle{1}}
\put(70,55){\circle{1}}
\put(80,55){\circle{1}}

\put(24,13){$\scriptstyle \be_1$}
\put(19,20){$\scriptstyle \be_2$}
\put(28,20){$\scriptstyle \be_3$}
\put(43.5,33){$\scriptstyle\be_4$}
\put(39,30){$\scriptstyle\be_5$}
\put(48.5,30){$\scriptstyle\be_6$}

\put(34.0,13){$\scriptstyle \ga_1$}
\put(24.5,31.0){$\scriptstyle \ga_2$}
\put(20.5,13){$\scriptstyle x_1$}
\put(11.5,33.0){$\scriptstyle x_3-\ga_1+\ga_2$}
\put(22.0,36.0){$\scriptstyle x_1+\ga_2$}
\put(39.5,36.0){$\scriptstyle x_3+\ga_2$}
\put(40.0,13){$\scriptstyle x_1+\ga_1$}
\put(21.0,25.5){$\scriptstyle x_2$}
\put(46.0,24.5){$\scriptstyle x_2+\ga_1$}
\put(50.0,33.0){$\scriptstyle x_1+\ga_1+\ga_2$}
\put(29.5,16.5){$\scriptstyle x_3$}
\put(-5,7.0){\emph{a})}
\end{picture}
\begin{picture}(90,50)(0,0)
\bezier{30}(11,15)(17.5,28)(24,41)
\bezier{30}(12,15)(18.5,28)(25,41)
\bezier{30}(13,15)(19.5,28)(26,41)
\bezier{30}(14,15)(20.5,28)(27,41)
\bezier{30}(15,15)(21.5,28)(28,41)
\bezier{30}(16,15)(22.5,28)(29,41)
\bezier{30}(17,15)(23.5,28)(30,41)
\bezier{30}(18,15)(24.5,28)(31,41)
\bezier{30}(19,15)(25.5,28)(32,41)
\bezier{30}(20,15)(26.5,28)(33,41)
\bezier{30}(21,15)(27.5,28)(34,41)
\bezier{30}(22,15)(28.5,28)(35,41)
\bezier{30}(23,15)(29.5,28)(36,41)
\bezier{30}(24,15)(30.5,28)(37,41)
\bezier{30}(25,15)(31.5,28)(38,41)
\bezier{30}(26,15)(32.5,28)(39,41)
\bezier{30}(27,15)(33.5,28)(40,41)
\bezier{30}(28,15)(34.5,28)(41,41)
\bezier{30}(29,15)(35.5,28)(42,41)
\bezier{30}(30,15)(36.5,28)(43,41)
\bezier{30}(31,15)(37.5,28)(44,41)
\bezier{30}(32,15)(38.5,28)(45,41)
\bezier{30}(33,15)(39.5,28)(46,41)
\bezier{30}(34,15)(40.5,28)(47,41)
\bezier{30}(35,15)(41.5,28)(48,41)

\put(27,27){$\Omega$}
\put(22,48){$\bK_*$}
\multiput(24,41)(5,0){5}{\line(1,0){3}}
\multiput(37,17)(3,6){4}{\line(1,2){2}}

\put(10,15){\vector(1,0){26.0}}
\put(10,15){\vector(1,2){13.0}}
\bezier{30}(23,41)(36,41)(49,41)
\bezier{30}(36,15)(42.5,28)(49,41)
\put(10,15){\circle*{1.5}}
\put(23,15){\circle*{1.5}}
\put(36,15){\circle{1}}

\put(23,15){\vector(-1,0){13}}
\put(10,15){\vector(1,2){6.5}}
\put(16.5,28){\vector(1,-2){6.5}}

\put(36,41){\vector(1,-2){6.5}}
\put(37,41){\vector(-1,0){1}}
\put(48.3,40){\vector(1,2){0.5}}

\put(23,15){\line(-1,2){6.5}}
\put(42.5,28){\line(-1,2){6.5}}

\put(16.5,28){\circle*{1.5}}
\put(42.5,28){\circle{1}}

\put(23,41){\circle{1}}
\put(36,41){\circle{1}}
\put(49,41){\circle{1}}

\put(26.7,12.0){$\ga_1$}
\put(15,34.5){$\ga_2$}
\put(8.0,12.0){$x_1$}
\put(21.0,42.7){$x_1$}
\put(32.0,42.7){$x_3$}
\put(34.0,12.0){$x_1$}
\put(11.5,28.0){$x_2$}
\put(37.5,28.0){$x_2$}
\put(46.0,42.7){$x_1$}
\put(20.0,12.0){$x_3$}
\put(0,7.0){\emph{b})}

\put(40,42.7){$\be_4$}
\put(15.0,16){$\be_1$}
\put(20.5,21){$\be_3$}
\put(8.5,21){$\be_2$}
\put(47.0,34){$\be_6$}
\put(35.0,34){$\be_5$}

\put(55,15){\line(1,0){26.0}}

\put(55,15){\line(1,2){13.0}}

\put(80.8,29.5){\vector(-1,2){1.0}}
\put(56.2,31.5){\vector(-1,-2){1.0}}
\put(68,9.0){\vector(1,0){1.0}}

\put(81,15){\line(-1,2){13.0}}

\put(55,15){\circle*{1.5}}
\put(50.5,13.5){$x_1$}
\put(81,15){\circle*{1.5}}
\put(82,13.5){$x_3$}
\put(68,41){\circle*{1.5}}
\put(68.0,42){$x_2$}
\bezier{200}(55,15)(68,3)(81,15)
\bezier{200}(55,15)(50,33)(68,41)
\bezier{200}(81,15)(86,33)(68,41)

\put(69.0,16.0){$\scriptstyle(0,0)$}
\put(60.5,24.0){$\scriptstyle(0,0)$}
\put(69.5,24.0){$\scriptstyle(0,0)$}

\put(64.0,5.0){$\scriptstyle(-1,0)$}
\put(53.0,37.0){$\scriptstyle(0,1)$}
\put(75.0,37.0){$\scriptstyle(1,-1)$}
\put(61.5,16){$\be_1$}
\put(66.5,10.5){$\be_4$}
\put(62.5,28){$\be_2$}
\put(70,28){$\be_3$}
\put(51.0,30){$\be_6$}
\put(81.5,30){$\be_5$}

\put(61.0,27.0){\vector(1,2){1.0}}
\put(74.2,28.5){\vector(1,-2){1.0}}
\put(68,15.0){\vector(-1,0){1.0}}

\end{picture}

\vspace{-0.7cm} \caption{\footnotesize \emph{a}) The Kagome lattice $\bK$; $\ga_1,\ga_2$ are the periods of $\bK$; the fundamental cell $\Omega$ is shaded; the shifted copies of $\Omega$ have the coordinates $(-1,0)$, $(0,0)$, $(1,0)$, $(-1,1)$ and so on; the influence of the magnetic field on the spectrum of $\D_\a$ on $\bK$ is completely determined by two fluxes $\phi_1$ and $\phi_2$ through the cycles $\bc_1=(\be_1,\be_2,\be_3)$ and $\bc_2=(\be_4,\be_5,\be_6)$ of $\bK$;\, \emph{b})~ the fundamental graph $\bK_*$; the indices are shown near the corresponding edges.}
\label{FKaL}
\end{figure}

\subsection{Edge indices} \lb{Sedin} We define the important notion of an {\it edge index} which was introduced in \cite{KS14}. This notion allows one to consider, instead of a periodic graph, the finite fundamental graph with edges labeled by some integer vectors called indices.

For each $\bx\in\R^d$, we denote by $\bx_\A\in\R^d$ the coordinate vector of $\bx$ with respect to the basis  $\A=\{\ga_1,\ldots,\ga_d\}$ of the lattice~$\G$, i.e.,
\[\lb{cola}
\bx_\A=(x_1,\ldots,x_d), \qqq \textrm{where} \qq \bx=
\textstyle\sum\limits_{s=1}^dx_s\ga_s.
\]

For any vertex $x\in\cV$ of a $\G$-periodic graph $\cG$ the
following  unique representation holds true:
\[\lb{Dv}
x=x_0+[x], \qq \textrm{where}\qq x_0\in\cV\cap\Omega,\qquad [x]\in\G,
\]
and $\Omega$ is the fundamental cell of the lattice $\G$ defined by
\er{fuce}. In other words, each vertex $x$ can be obtained from a
vertex $x_0$ from the fundamental cell $\Omega$ by a shift by a vector $[x]\in\G$. For any
oriented edge $\be=(x,y)\in\cA$ of the periodic graph $\cG$ we define the \emph{edge index}
$\t(\be)$ as the vector of the lattice $\Z^d$ given by
\[
\lb{in}
\t(\be)=[y]_\A-[x]_\A\in\Z^d,
\]
where $[x]\in\G$ is defined by \er{Dv} and the vector
$[x]_\A\in\Z^d$  is given by \er{cola}.

Due to the periodicity of the graph $\cG$, the edge indices $\t(\be)$, $\be\in\cA$, satisfy
\[\lb{Gpe}
\t(\be+\ga)=\t(\be),\qqq \forall\, (\be,\ga)\in\cA \ts\G.
\]
On the set $\cA$ of all oriented edges of the $\G$-periodic graph $\cG$ we define the surjection
\[\lb{sur}
\gf:\cA\rightarrow\cA_*=\cA/\G,
\]
which maps each $\be\in\cA$ to its equivalence class $\be_*=\gf(\be)$ which is an oriented edge of the fundamental graph $\cG_*$.

For any oriented edge $\be_*\in\cA_*$ of the fundamental graph $\cG_*=(\cV_*,\cA_*)$ we define the \emph{edge index}  $\t(\be_*)\in\Z^d$ by
\[
\lb{dco}
\t(\be_*)=\t(\be) \qq \textrm{ for some $\be\in\cA$ \; such that }  \; \be_*=\gf(\be),
\]
where $\gf$ is defined by \er{sur}. In other words, edge indices of
the fundamental graph $\cG_*$ are induced by edge indices of the
periodic graph~$\cG$. Due to \er{Gpe}, the edge index $\t(\be_*)$ is uniquely determined by \er{dco}. From this definition it also follows that
\[\lb{inin0}
\t(\ul\be\,)=-\t(\be), \qqq \forall\,\be\in\cA_*,
\]
where $\ul\be$ is the inverse edge of $\be$.

\begin{remark}
\emph{i}) Indices of edges of a $\G$-periodic graph $\cG$ depend on the choice of the embedding of $\cG$ into $\R^d$ and on the choice of basis $\A=\{\ga_1,\ldots,\ga_d\}$ of the lattice $\G$. But when the embedding $\cG\ss\R^d$ and the basis $\A$ are fixed, the edge indices of $\cG$ are uniquely defined.

\emph{ii}) We can give the following simple geometric interpretation of edge indices. We define the \emph{coordinates} of a shifted copy $\Omega+\ga$, $\ga\in\G$, of the fundamental cell $\Omega$ as the coordinate vector of $\ga$ with respect to the basis $\A$. For example, the cell $\Omega+\ga_1$ has the coordinates $(1,0,0,\ldots,0)$, the cell $\Omega+\ga_2$ has the coordinates $(0,1,0,\ldots,0)$ and so on (see Fig.\ref{FKaL}\emph{a}). The coordinates of the fundamental cell $\Omega$ is $(0,\ldots,0)$. Then the edge index is equal to the difference between the coordinates of the cells, containing the terminal and initial vertices of the edge. Thus, in some sense, the edge index shows how long the edge is. For example, for the graph $\bK$ shown in Fig.\ref{FKaL}\emph{a}, the index of the edge $\be_5=(x_3+\ga_2,x_2+\ga_1)$ equals $(1,-1)$, since this edge starts in the cell with coordinates $(0,1)$ and ends at the cell with coordinates $(1,0)$. Edges connecting vertices from the same cell have zero indices.
\end{remark}

\begin{example1} We calculate edge indices for the periodic and fundamental graphs from Example \ref{ExPG}, using the definitions \er{in}, \er{dco}.

\emph{i}) The fundamental graph $\bS_*=(\cV_*,\cA_*)$ of the square lattice $\bS=(\cV,\cA)$ consists of one vertex $x$ and two loop edges $\be_1,\be_2$ at this vertex (Fig.\ref{FSqL}\emph{b}). The index of the loop $\be_1\in\cA_*$ is $(1,0)$, since the corresponding edge $\be_1=(x,x+\ga_1)\in\cA$ of the square lattice $\bS$ has the index
$$
\t(\be_1)=[x+\ga_1]_\A-[x]_\A=(1,0)-(0,0)=(1,0).
$$
Similarly, the index of the loop $\be_2\in\cA_*$ is $(0,1)$, since the corresponding edge $\be_2=(x,x+\ga_2)\in\cA$ of $\bS$ has the index
$$
\t(\be_2)=[x+\ga_2]_\A-[x]_\A=(0,1)-(0,0)=(0,1).
$$
Thus, the edge indices of the fundamental graph $\bS_*$ of the square lattice $\bS$ are given by
\[\lb{insl}
\t(\be_1)=(1,0), \qqq \t(\be_2)=(0,1).
\]

\emph{ii}) The fundamental graph $\bG_*=(\cV_*,\cA_*)$ of the hexagonal lattice $\bG=(\cV,\cA)$ consists of two vertices $x_1,x_2$ and three multiple edges  $\be_1,\be_2,\be_3$ connecting these vertices (Fig.\ref{ff.0.3}\emph{b}). The index of the edge $\be_1=(x_2,x_1)\in\cA_*$ is $(0,0)$, since the corresponding edge $\be_1=(x_2,x_1)\in\cA$ of the hexagonal lattice $\bG$ has the index
$$
\t(\be_1)=[x_1]_\A-[x_2]_\A=(0,0)-(0,0)=(0,0).
$$
Similarly, the indices of the edges $\be_2,\be_3\in\cA_*$ is $(1,0)$ and $(0,1)$, respectively, since the corresponding edges $\be_2=(x_1,x_2+\ga_1)\in\cA$ and $\be_3=(x_1,x_2+\ga_2)\in\cA$ of $\bG$ have the indices
$$
\begin{array}{l}
\t(\be_2)=[x_2+\ga_1]_\A-[x_1]_\A=(1,0)-(0,0)=(1,0),\\[4pt]
\t(\be_3)=[x_2+\ga_2]_\A-[x_1]_\A=(0,1)-(0,0)=(0,1).
\end{array}
$$
Thus, the edge indices of the fundamental graph $\bG_*$ of the hexagonal lattice $\bG$ are given by
\[\lb{ingl}
\t(\be_1)=(0,0), \qqq \t(\be_2)=(1,0), \qqq \t(\be_3)=(0,1).
\]
Note that, by \er{inin0}, the inverse edges $\ul\be_1,\ul\be_2,\ul\be_3\in\cA_*$ have indices
$$
\t(\ul\be_1)=(0,0), \qqq \t(\ul\be_2)=(-1,0), \qqq \t(\ul\be_3)=(0,-1).
$$

\emph{iii}) The fundamental graph $\bK_*=(\cV_*,\cA_*)$ of the Kagome lattice $\bK=(\cV,\cA)$ consists of three vertices $x_1,x_2,x_3$ and six edges $\be_1,\ldots,\be_6$ (Fig.~\ref{FKaL}\emph{b}). The index of the edge $\be_1=(x_3,x_1)\in\cA_*$ is $(0,0)$, since the corresponding edge $\be_1=(x_3,x_1)\in\cA$ of the Kagome lattice $\bK$ has the index
$$
\t(\be_1)=[x_1]_\A-[x_3]_\A=(0,0)-(0,0)=(0,0).
$$
The index of the edge $\be_5=(x_3,x_2)\in\cA_*$ is $(1,-1)$, since the corresponding edge $\be_5=(x_3+\ga_2,x_2+\ga_1)\in\cA$ of $\bK$ has the index
$$
\t(\be_5)=[x_2+\ga_1]_\A-[x_3+\ga_2]_\A=
(1,0)-(0,1)=(1,-1).
$$
Using similar calculations for the remaining edges of $\bK_*$, we obtain the edge indices of the fundamental graph $\bK_*$ of the Kagome lattice $\bK$:
\[\lb{inKl}
\begin{aligned}
&\t(\be_1)=\t(\be_2)=\t(\be_3)=(0,0), \\
&\t(\be_4)=(-1,0), \qq \t(\be_5)=(1,-1),\qq \t(\be_6)=(0,1).
\end{aligned}
\]
\end{example1}

\subsection{Magnetic Schr\"odinger operators on graphs.}
Let $\ell^2(\cV)$ be the Hilbert space of all square summable
functions  $f:\cV\to \C$ equipped with the norm
$$
\|f\|^2_{\ell^2(\cV)}=\sum_{x\in\cV}|f_x|^2<\infty.
$$

We consider the magnetic Schr\"odinger operator $H_\a$ acting on $\ell^2(\cV)$ and given by
\[\lb{Sh}
H_\a=-\D_\a+V,\qqq \D_\a=\vk-A_\a,
\]
where $V$ is a real electric potential and $\vk$ is a degree potential given by
\[\lb{DMO}
(Vf)_x=V_xf_x, \qqq (\vk f)_x=\vk_xf_x, \qqq f\in\ell^2(\cV), \qqq x\in\cV,
\]
$\vk_x$ is the degree of the vertex $x$; $\D_\a$ is the combinatorial magnetic Laplacian, and $A_\a$ is the magnetic \emph{adjacency} operator having the form
\[
\lb{ALO}
(A_\a f)_x=\sum_{\be=(x,y)\in\cA}e^{i\a(\be)}f_y, \qqq f\in\ell^2(\cV), \qqq x\in\cV.
\]
Here $\a:\cA\ra\R$ is a \emph{magnetic potential on $\cG=(\cV,\cA)$}, i.e., it satisfies:
\[\lb{mpp}
\a(\ul\be\,)=-\a(\be),\qqq \forall\,\be\in\cA,
\]
where $\ul\be$ is the inverse edge of $\be$. The sum in \er{ALO} is taken over all edges from $\cA$ starting at the vertex $x$. It is known that the magnetic Schr\"odinger operator $H_\a$ is a bounded self-adjoint operator on $\ell^2(\cV)$ (see, e.g., \cite{HS99a}) and the spectra of $A_\a$ and $\D_\a$ satisfy
\[\lb{spAD}
\s(A_\a)\subseteq[-\vk_+,\vk_+],\qqq \s(\D_\a)\subseteq[0,2\vk_+],\qqq\where \qqq \vk_+=\max\limits_{x\in\cV_*}\vk_x.
\]

We assume that the magnetic potential $\a$ and the
electric potential $V$ are $\G$-periodic, i.e., they satisfy
\[\lb{Gpps}
\a(\be+\ga)=\a(\be),\qqq V_{x+\ga}=V_x,\qqq \forall\,(x,\be,\ga)\in\cV\ts\cA\ts\G.
\]
The periodicity of the magnetic potential $\a$ guarantees a band structure of the spectrum and the absence of Cantor spectrum.

The first identity in \er{Gpps} allows us to define uniquely the \emph{magnetic potential} $\a$ on the fundamental graph $\cG_*=(\cV_*,\cA_*)$ which is induced by the magnetic potential $\a$ on the periodic graph $\cG=(\cV,\cA)$:
\[\lb{mpf}
\a(\be_*)=\a(\be)\qq \textrm{ for some $\be\in\cA$ \; such that }  \; \be_*=\gf(\be), \qqq \be_*\in\cA_*,
\]
where $\gf$ is defined by \er{sur}. From \er{mpp} and \er{mpf} it follows that
\[\lb{inin}
\a(\ul\be\,)=-\a(\be), \qqq \forall\,\be\in\cA_*.
\]

\begin{remark} \emph{i}) If $\a=0$, then $\D_0$ is just the usual combinatorial Laplacian $\D$ without magnetic potentials:
$$
(\D f)_x=\sum_{(x,y)\in\cA}(f_x-f_y), \qqq \forall\, f\in\ell^2(\cV),\qqq \forall\, x\in\cV.
$$

\emph{ii}) If $\cG$ is a regular graph of degree $\vk_+$, i.e., all vertices of $\cG$ have the same degree $\vk_+$, then the magnetic Laplacian $\D_\a$ has the the form
\[\lb{suD}
\D_\a=\vk_+I-A_\a,
\]
where $I$ is the identity operator, and $A_\a$ is the magnetic adjacency operator  given by \er{ALO}. Thus, the operators $-\D_\a$ and $A_\a$ on a regular graph differ only by a shift.

\emph{iii})  The case of periodic graphs in uniform and periodic magnetic fields with rational fluxes can also be reduced to the Laplacians with periodic magnetic potentials. Indeed, for a $\G$-invariant magnetic field $\mB$, the magnetic potential $\a$ is \emph{weakly $\G$-invariant}, i.e., for any $\ga\in\G$ there exists a function $\varphi_\ga:\cV\to\R$ such that
\[\lb{MPot}
\a(\be+\ga)=\a(\be)+\varphi_\ga(u)-\varphi_\ga(v),\qqq \forall\,\be=(u,v)\in\cA,
\]
for more details see \cite{S94}. If the (appropriately defined) fluxes $\f$ of the magnetic field $\mB$ satisfy $\frac\f{2\pi}=\frac pq$ for some relatively prime numbers $p\in\Z$ and $q\in\N$, then there exists a gauge transformation which reduces the Laplacian $\D_\a$ with the \emph{weakly $\G$-invariant} magnetic potential $\a$ to the Laplacian $\D_{\wt\a}$ with a magnetic potential $\wt\a$ which is periodic with periods $q\ga_1,\ldots,q\ga_{d-1},\ga_d$. Thus, enlargement of the lattice $\G$ permits reduction of the problem to a periodic magnetic potential by increasing the size of the fundamental cell in each direction $\ga_1,\ldots,\ga_{d-1}$ by $q$ times. Note that the electric potential $V$ is invariant under the gauge transformation.
\end{remark}

\begin{example1}
The \emph{Harper operator} $\D_\f$ is the discrete magnetic Laplacian on the square lattice $\bS=(\Z^2,\cA)$ in a uniform magnetic field $\mB$ orthogonal to the lattice and with flux $\f$ through the unit cell of the lattice. It acts on $f\in\ell^2(\Z^2)$ and has the form
\[
\lb{haop}
\begin{array}{l}
(\D_\f f)_x=e^{-i\f{x_2\/2}}f_{x+\ga_1}+e^{i\f{x_2\/2}}f_{x-\ga_1}+
e^{i\f{x_1\/2}}f_{x+\ga_2}+e^{-i\f{x_1\/2}}f_{x-\ga_2},\\[10pt]
x=(x_1,x_2)\in\Z^2,\qqq \ga_1=(1,0),\qq \ga_2=(0,1).
\end{array}
\]
The magnetic potential $\a:\cA\to\R$ of the uniform magnetic field $\mB$ on edges of $\bS$ is given by
\[\lb{abe}
\a(\be)=\left\{
\begin{array}{rl}
  -\,{\f x_2\/2}\,, & \textrm{if } \be=(x,x+\ga_1) \\[4pt]
  {\f x_1\/2}\,, & \textrm{if } \be=(x,x+\ga_2)
\end{array}\right.,\qqq \be\in\cA,
\]
and $\a$ is not periodic (see Fig.\ref{ff.0.1}\emph{a}).

Let ${\f\/2\pi}={p\/q}$\,, where $p\in\Z$ and $q\in\N$ are relatively prime.
We define the magnetic potential $\wt\a:\cA\to\R$ by
\[\lb{mpsl*}
\wt\a(\be)=\left\{
\begin{array}{cl}
  0, & \textrm{if } \be=(x,x+\ga_1) \\[6pt]
  \f x_1, & \textrm{if } \be=(x,x+\ga_2)
\end{array}\right.,\qqq \be\in\cA,
\]
see Fig.\ref{ff.0.1}\emph{b}. The potential $\wt\a$ satisfies
$$
\wt\a(\be+q\ga_1)=\left\{
\begin{array}{cl}
  \wt\a(\be), & \textrm{if } \be=(x,x+\ga_1) \\[6pt]
  \wt\a(\be)+2\pi p, & \textrm{if } \be=(x,x+\ga_2)
\end{array}\right.,\qqq \wt\a(\be+\ga_2)=\wt\a(\be),\qqq \forall\,\be\in\cA.
$$
Note that we may consider a magnetic potential $\a$ modulo $2\pi$, since it appears in the magnetic adjacency operator $A_\a$ via the magnetic phase $e^{i\a(\be)}$, $\be\in\cA$, see \er{ALO}.
Thus, the magnetic potential $\wt\a$ (modulo $2\pi$) is periodic with periods $q\ga_1,\ga_2$.

The magnetic potentials $\a$ and $\wt\a$ given by \er{abe} and \er{mpsl*}, respectively, have the same flux through each cycle of the square lattice $\bS$ (Fig.\ref{ff.0.1}\emph{a},\emph{b}). Then the Laplacian $\D_\f$ with the non-periodic potential $\a$ is gauge equivalent to the Laplacian $\wt\D_\f$ with the periodic potential $\wt\a$ with periods $q\ga_1,\ga_2$. The extended fundamental graph $\wt\bS_*=\bS/\wt\G$, where $\wt\G$ is the lattice with basis $q\ga_1,\ga_2$, is shown in Fig.\ref{ff.0.1}\emph{d} (when $q=3$).
\end{example1}

\setlength{\unitlength}{1.0mm}
\begin{figure}[h]
\centering
\unitlength 1.4mm 
\linethickness{0.4pt}
\ifx\plotpoint\undefined\newsavebox{\plotpoint}\fi 
\begin{picture}(130,45)(0,0)

\put(5,10){\line(1,0){40.00}}
\put(5,20){\line(1,0){40.00}}
\put(5,30){\line(1,0){40.00}}
\put(5,40){\line(1,0){40.00}}
\put(20,5){\line(0,1){40.00}}
\put(20,20){\vector(1,0){10.00}}
\put(20,20){\vector(0,1){10.00}}

\put(24,24){$\Omega$}

\put(14.0,11.0){$\scriptstyle \frac{\f}2$}
\put(24.0,11.0){$\scriptstyle \frac{\f}2$}
\put(34.0,11.0){$\scriptstyle \frac{\f}2$}

\put(14.5,20.3){$\scriptstyle 0$}
\put(24.5,20.3){$\scriptstyle 0$}
\put(34.5,20.3){$\scriptstyle 0$}

\put(13.0,31.0){$\scriptstyle -\frac{\f}2$}
\put(23.0,31.0){$\scriptstyle -\frac{\f}2$}
\put(33.0,31.0){$\scriptstyle -\frac{\f}2$}

\put(13.0,40.3){$\scriptstyle -\f$}
\put(23.0,40.3){$\scriptstyle -\f$}
\put(33.0,40.3){$\scriptstyle -\f$}

\put(10.5,14){$\scriptstyle -\frac{\f}2$}
\put(20.5,14){$\scriptstyle 0$}
\put(30.2,14){$\scriptstyle \frac{\f}2$}
\put(40.2,14){$\scriptstyle \f$}\

\put(10.5,24){$\scriptstyle -\frac{\f}2$}
\put(20.5,24){$\scriptstyle 0$}
\put(30.2,24){$\scriptstyle \frac{\f}2$}
\put(40.2,24){$\scriptstyle \f$}

\put(10.5,34){$\scriptstyle -\frac{\f}2$}
\put(20.5,34){$\scriptstyle 0$}
\put(30.2,34){$\scriptstyle \frac{\f}2$}
\put(40.2,34){$\scriptstyle \f$}

\put(27.5,18.3){$\scriptstyle \ga_1$}
\put(17.5,28.2){$\scriptstyle \ga_2$}

\put(5,10){\vector(1,0){5.00}}
\put(5,20){\vector(1,0){5.00}}
\put(5,30){\vector(1,0){5.00}}
\put(5,40){\vector(1,0){5.00}}

\put(10,10){\vector(1,0){10.00}}
\put(10,20){\vector(1,0){10.00}}
\put(10,30){\vector(1,0){10.00}}
\put(10,40){\vector(1,0){10.00}}

\put(20,10){\vector(1,0){10.00}}
\put(20,20){\vector(1,0){10.00}}
\put(20,30){\vector(1,0){10.00}}
\put(20,40){\vector(1,0){10.00}}

\put(30,10){\vector(1,0){10.00}}
\put(30,20){\vector(1,0){10.00}}
\put(30,30){\vector(1,0){10.00}}
\put(30,40){\vector(1,0){10.00}}

\put(10,10){\vector(0,1){10.00}}
\put(20,10){\vector(0,1){10.00}}
\put(30,10){\vector(0,1){10.00}}
\put(40,10){\vector(0,1){10.00}}

\put(10,20){\vector(0,1){10.00}}
\put(20,20){\vector(0,1){10.00}}
\put(30,20){\vector(0,1){10.00}}
\put(40,20){\vector(0,1){10.00}}

\put(10,30){\vector(0,1){10.00}}
\put(20,30){\vector(0,1){10.00}}
\put(30,30){\vector(0,1){10.00}}
\put(40,30){\vector(0,1){10.00}}

\put(10,5){\vector(0,1){5.00}}
\put(20,5){\vector(0,1){5.00}}
\put(30,5){\vector(0,1){5.00}}
\put(40,5){\vector(0,1){5.00}}

\put(10,5){\line(0,1){40.00}}
\put(30,5){\line(0,1){40.00}}
\put(40,5){\line(0,1){40.00}}

\bezier{20}(20.5,20)(20.5,25)(20.5,30)
\bezier{20}(21.0,20)(21.0,25)(21.0,30)
\bezier{20}(21.5,20)(21.5,25)(21.5,30)
\bezier{20}(22.0,20)(22.0,25)(22.0,30)
\bezier{20}(22.5,20)(22.5,25)(22.5,30)
\bezier{20}(23.0,20)(23.0,25)(23.0,30)
\bezier{20}(23.5,20)(23.5,25)(23.5,30)
\bezier{20}(24.0,20)(24.0,25)(24.0,30)
\bezier{20}(24.5,20)(24.5,25)(24.5,30)
\bezier{20}(25.0,20)(25.0,25)(25.0,30)
\bezier{20}(25.5,20)(25.5,25)(25.5,30)
\bezier{20}(26.0,20)(26.0,25)(26.0,30)
\bezier{20}(26.5,20)(26.5,25)(26.5,30)
\bezier{20}(27.0,20)(27.0,25)(27.0,30)
\bezier{20}(27.5,20)(27.5,25)(27.5,30)
\bezier{20}(28.0,20)(28.0,25)(28.0,30)
\bezier{20}(28.5,20)(28.5,25)(28.5,30)
\bezier{20}(29.0,20)(29.0,25)(29.0,30)
\bezier{20}(29.5,20)(29.5,25)(29.5,30)
\put(10,10){\circle{1}}
\put(20,10){\circle{1}}
\put(30,10){\circle{1}}
\put(40,10){\circle{1}}

\put(10,20){\circle{1}}
\put(20,20){\circle{1}}
\put(30,20){\circle{1}}
\put(40,20){\circle{1}}

\put(10,30){\circle{1}}
\put(20,30){\circle{1}}
\put(30,30){\circle{1}}
\put(40,30){\circle{1}}

\put(10,40){\circle{1}}
\put(20,40){\circle{1}}
\put(30,40){\circle{1}}
\put(40,40){\circle{1}}

\put(10.5,8.0){$\scriptstyle(-1,-1)$}
\put(20.5,8.0){$\scriptstyle(0,-1)$}
\put(30.5,8.0){$\scriptstyle(1,-1)$}
\put(40.5,8.0){$\scriptstyle(2,-1)$}
\put(10.5,18.0){$\scriptstyle(-1,0)$}
\put(20.5,18.0){$\scriptstyle(0,0)$}
\put(30.5,18.0){$\scriptstyle(1,0)$}
\put(40.5,18.0){$\scriptstyle(2,0)$}
\put(10.5,28.0){$\scriptstyle(-1,1)$}
\put(20.5,28.0){$\scriptstyle(0,1)$}
\put(30.5,28.0){$\scriptstyle(1,1)$}
\put(40.5,28.0){$\scriptstyle(2,1)$}
\put(10.5,38.0){$\scriptstyle(-1,2)$}
\put(20.5,38.0){$\scriptstyle(0,2)$}
\put(30.5,38.0){$\scriptstyle(1,2)$}
\put(40.5,38.0){$\scriptstyle(2,2)$}
\put(3,42){$\bS,\a$}

\put(3,5){\emph{a})}
\linethickness{1.8pt}
\put(65,10){\line(1,0){40.00}}
\put(65,20){\line(1,0){40.00}}
\put(65,30){\line(1,0){40.00}}
\put(65,40){\line(1,0){40.00}}
\put(80,5){\line(0,1){40.00}}
\put(80,20){\vector(1,0){10.00}}
\put(80,20){\vector(0,1){10.00}}
\linethickness{0.4pt}

\put(84,23.5){$\Omega$}

\put(70,24){$\scriptstyle-\f$}
\put(90.5,24){$\scriptstyle \f$}
\put(100.4,24){$\scriptstyle 2\f$}

\put(70,14){$\scriptstyle-\f$}
\put(90.5,14){$\scriptstyle \f$}
\put(100.4,14){$\scriptstyle 2\f$}

\put(70,34){$\scriptstyle-\f$}
\put(90.5,34){$\scriptstyle \f$}
\put(100.4,34){$\scriptstyle 2\f$}

\put(87.5,18.3){$\scriptstyle \ga_1$}
\put(77.5,28.2){$\scriptstyle \ga_2$}

\put(70,10){\vector(0,1){10.00}}
\put(90,10){\vector(0,1){10.00}}
\put(100,10){\vector(0,1){10.00}}

\put(70,20){\vector(0,1){10.00}}
\put(90,20){\vector(0,1){10.00}}
\put(100,20){\vector(0,1){10.00}}

\put(70,30){\vector(0,1){10.00}}
\put(90,30){\vector(0,1){10.00}}
\put(100,30){\vector(0,1){10.00}}

\put(70,5){\vector(0,1){5.00}}
\put(90,5){\vector(0,1){5.00}}
\put(100,5){\vector(0,1){5.00}}

\put(70,5){\line(0,1){40.00}}
\put(90,5){\line(0,1){40.00}}
\put(100,5){\line(0,1){40.00}}

\bezier{20}(80.5,20)(80.5,25)(80.5,30)
\bezier{20}(81.0,20)(81.0,25)(81.0,30)
\bezier{20}(81.5,20)(81.5,25)(81.5,30)
\bezier{20}(82.0,20)(82.0,25)(82.0,30)
\bezier{20}(82.5,20)(82.5,25)(82.5,30)
\bezier{20}(83.0,20)(83.0,25)(83.0,30)
\bezier{20}(83.5,20)(83.5,25)(83.5,30)
\bezier{20}(84.0,20)(84.0,25)(84.0,30)
\bezier{20}(84.5,20)(84.5,25)(84.5,30)
\bezier{20}(85.0,20)(85.0,25)(85.0,30)
\bezier{20}(85.5,20)(85.5,25)(85.5,30)
\bezier{20}(86.0,20)(86.0,25)(86.0,30)
\bezier{20}(86.5,20)(86.5,25)(86.5,30)
\bezier{20}(87.0,20)(87.0,25)(87.0,30)
\bezier{20}(87.5,20)(87.5,25)(87.5,30)
\bezier{20}(88.0,20)(88.0,25)(88.0,30)
\bezier{20}(88.5,20)(88.5,25)(88.5,30)
\bezier{20}(89.0,20)(89.0,25)(89.0,30)
\bezier{20}(89.5,20)(89.5,25)(89.5,30)
\put(70,10){\circle{1}}
\put(80,10){\circle{1}}
\put(90,10){\circle{1}}
\put(100,10){\circle{1}}

\put(70,20){\circle{1}}
\put(80,20){\circle{1}}
\put(90,20){\circle{1}}
\put(100,20){\circle{1}}

\put(70,30){\circle{1}}
\put(80,30){\circle{1}}
\put(90,30){\circle{1}}
\put(100,30){\circle{1}}

\put(70,40){\circle{1}}
\put(80,40){\circle{1}}
\put(90,40){\circle{1}}
\put(100,40){\circle{1}}

\put(70.5,8.0){$\scriptstyle(-1,-1)$}
\put(80.5,8.0){$\scriptstyle(0,-1)$}
\put(90.5,8.0){$\scriptstyle(1,-1)$}
\put(100.5,8.0){$\scriptstyle(2,-1)$}
\put(70.5,18.0){$\scriptstyle(-1,0)$}
\put(80.5,18.0){$\scriptstyle(0,0)$}
\put(90.5,18.0){$\scriptstyle(1,0)$}
\put(100.5,18.0){$\scriptstyle(2,0)$}
\put(70.5,28.0){$\scriptstyle(-1,1)$}
\put(80.5,28.0){$\scriptstyle(0,1)$}
\put(90.5,28.0){$\scriptstyle(1,1)$}
\put(100.5,28.0){$\scriptstyle(2,1)$}
\put(70.5,38.0){$\scriptstyle(-1,2)$}
\put(80.5,38.0){$\scriptstyle(0,2)$}
\put(90.5,38.0){$\scriptstyle(1,2)$}
\put(100.5,38.0){$\scriptstyle(2,2)$}
\put(63,42){$\bS,\wt\a$}

\put(63,5){\emph{b})}
\end{picture}

\begin{picture}(130,20)(0,0)
\put(24,9){$\Omega$}
\put(24,19){$\bS_*$}
\put(20,5){\line(1,0){10.00}}
\put(20,5){\line(0,1){10.00}}

\bezier{20}(20.5,5)(20.5,10)(20.5,15)
\bezier{20}(21.0,5)(21.0,10)(21.0,15)
\bezier{20}(21.5,5)(21.5,10)(21.5,15)
\bezier{20}(22.0,5)(22.0,10)(22.0,15)
\bezier{20}(22.5,5)(22.5,10)(22.5,15)
\bezier{20}(23.0,5)(23.0,10)(23.0,15)
\bezier{20}(23.5,5)(23.5,10)(23.5,15)
\bezier{20}(24.0,5)(24.0,10)(24.0,15)
\bezier{20}(24.5,5)(24.5,10)(24.5,15)
\bezier{20}(25.0,5)(25.0,10)(25.0,15)
\bezier{20}(25.5,5)(25.5,10)(25.5,15)
\bezier{20}(26.0,5)(26.0,10)(26.0,15)
\bezier{20}(26.5,5)(26.5,10)(26.5,15)
\bezier{20}(27.0,5)(27.0,10)(27.0,15)
\bezier{20}(27.5,5)(27.5,10)(27.5,15)
\bezier{20}(28.0,5)(28.0,10)(28.0,15)
\bezier{20}(28.5,5)(28.5,10)(28.5,15)
\bezier{20}(29.0,5)(29.0,10)(29.0,15)
\bezier{20}(29.5,5)(29.5,10)(29.5,15)
\bezier{20}(30.0,5)(30.0,10)(30.0,15)

\put(20,5){\circle*{1}}
\put(20,15){\circle{1}}
\put(30,5){\circle{1}}
\put(30,15){\circle{1}}

\put(20,5){\vector(1,0){10.00}}
\put(20,5){\vector(0,1){10.00}}

\put(18,2.5){$x_1$}
\put(29,2.5){$x_1$}
\put(29,16.5){$x_1$}
\put(18,16.5){$x_1$}
\put(24,2.5){$\ga_1$}
\put(17,10){$\ga_2$}

\multiput(20,15)(4,0){3}{\line(1,0){2}}
\multiput(30,5)(0,4){3}{\line(0,1){2}}

\put(10,3){\emph{c})}

\put(74,9){$\Omega$}
\put(80.5,9){$\f$}
\put(90.5,9){$2\f$}
\put(84,19){$\wt\bS_*$}
\put(70,5){\line(1,0){30.00}}
\put(70,5){\line(0,1){10.00}}
\put(80,5){\line(0,1){10.00}}
\put(90,5){\line(0,1){10.00}}

\bezier{20}(70.5,5)(70.5,10)(70.5,15)
\bezier{20}(71.0,5)(71.0,10)(71.0,15)
\bezier{20}(71.5,5)(71.5,10)(71.5,15)
\bezier{20}(72.0,5)(72.0,10)(72.0,15)
\bezier{20}(72.5,5)(72.5,10)(72.5,15)
\bezier{20}(73.0,5)(73.0,10)(73.0,15)
\bezier{20}(73.5,5)(73.5,10)(73.5,15)
\bezier{20}(74.0,5)(74.0,10)(74.0,15)
\bezier{20}(74.5,5)(74.5,10)(74.5,15)
\bezier{20}(75.0,5)(75.0,10)(75.0,15)
\bezier{20}(75.5,5)(75.5,10)(75.5,15)
\bezier{20}(76.0,5)(76.0,10)(76.0,15)
\bezier{20}(76.5,5)(76.5,10)(76.5,15)
\bezier{20}(77.0,5)(77.0,10)(77.0,15)
\bezier{20}(77.5,5)(77.5,10)(77.5,15)
\bezier{20}(78.0,5)(78.0,10)(78.0,15)
\bezier{20}(78.5,5)(78.5,10)(78.5,15)
\bezier{20}(79.0,5)(79.0,10)(79.0,15)
\bezier{20}(79.5,5)(79.5,10)(79.5,15)
\bezier{20}(80.0,5)(80.0,10)(80.0,15)

\put(70,5){\circle*{1}}
\put(80,5){\circle*{1}}
\put(90,5){\circle*{1}}
\put(100,5){\circle{1}}
\put(70,15){\circle{1}}
\put(80,15){\circle{1}}
\put(90,15){\circle{1}}
\put(100,15){\circle{1}}

\put(70,5){\vector(1,0){30.00}}
\put(70,5){\vector(1,0){10.00}}
\put(70,5){\vector(0,1){10.00}}
\put(80,5){\vector(0,1){10.00}}
\put(90,5){\vector(0,1){10.00}}

\put(68,2.5){$x_1$}
\put(79,2.5){$x_2$}
\put(89,2.5){$x_3$}
\put(100,2.5){$x_1$}
\put(68,16.5){$x_1$}
\put(79,16.5){$x_2$}
\put(89,16.5){$x_3$}
\put(100,16.5){$x_1$}
\put(74,2.5){$\ga_1$}
\put(95,2.5){$3\ga_1$}
\put(67,10){$\ga_2$}

\multiput(70,15)(4,0){8}{\line(1,0){2}}
\multiput(100,5)(0,4){3}{\line(0,1){2}}

\put(60,3){\emph{d})}

\end{picture}

\caption{\footnotesize  \emph{a}) -- \emph{b}) The square lattice $\bS=(\Z^2,\cA)$; $\ga_1,\ga_2$ are the periods of $\bS$; the magnetic potentials $\a(\be)$ and $\wt\a(\be)$ are shown near the corresponding edges $\be\in\cA$; the magnetic potential $\wt\a$ on the bold edges is zero; \, \emph{c}) the minimal fundamental graph $\bS_*$ of $\bS$;\, \emph{d}) the extended fundamental graph $\wt\bS_*=\bS/\wt\G$, where $\wt\G$ is the lattice with basis $3\ga_1,\ga_2$.}
\label{ff.0.1}
\end{figure}

For other examples of such reductions of planar periodic graphs in uniform and periodic magnetic fields with rational fluxes to a periodic problem with extended fundamental graphs see \cite{BBR,BHJ19,HKeR16} and references therein.

\subsection{Spectrum of magnetic Schr\"odinger operators.}
The spectral analysis of periodic operators is usually based on the Floquet-Bloch theory (\cite{RS78}, Section XIII.16). We briefly describe the Floquet-Bloch technique for periodic discrete graphs (for more details see, e.g., \cite{HS99a} or \cite{KS17}). This technique allows one to reduce the study of the spectrum of the magnetic Schr\"odinger operator $H_\a$ on a periodic graph $\cG$ to a family of spectral problems on the finite fundamental graph $\cG_*$.

Let $\T^d=\R^d/(2\pi\Z)^d=[-\pi,\pi]^d$ be the $d$-dimensional torus. We introduce the Hilbert space
\[\lb{did}
\mH=L^2\Big(\T^d,\frac{dk}{(2\pi)^d}\,;\ell^2(\cV_*)\Big)
\]
of all square integrable functions on $\T^d$ with values in $\ell^2(\cV_*)$
equipped with the norm
$$
\|g\|^2_{\mH}=\int_{\T^d}\|g(k)\|_{\ell^2(\cV_*)}^2\frac{dk}{(2\pi)^d}\,, \qqq g\in\mH.
$$
The parameter $k\in\T^d$ is called the \emph{quasimomentum} (the name comes from the solid state physics).

We identify the vertices of the periodic graph $\cG$ from the fundamental cell $\Omega$ with the corresponding vertices of the fundamental graph $\cG_*=(\cV_*,\cA_*)$ and define the unitary operator $U:\ell^2(\cV)\to \mH$ (the discrete \emph{Floquet} transform) by
\[\lb{DFT}
(Uf)_x(k)=\sum\limits_{\ga\in\G}e^{-i\lan[\ga]_\A,k\ran }
f_{x+\ga}, \qqq \forall\,(k,x)\in \T^d\ts\cV_*,
\]
where $\lan \cdot,\cdot\ran$ denotes the standard inner product in $\R^d$, and $[\ga]_\A\in\Z^d$ is the coordinate vector of $\ga\in\G$ with respect to the basis $\A=\{\ga_1,\ldots,\ga_d\}$ of the lattice~$\G$.

After the Floquet transform the magnetic Schr\"odinger operator $H_\a$ on a periodic graph $\cG$ becomes an operator of multiplication in $\mH$:
\[\lb{did1}
(UH_\a U^{-1}g)(k)=H_\a(k)g(k),\qqq \forall\, g\in\mH, \qq \forall\, k\in\T^d,
\]
and the spectrum $\s(H_\a)$ of $H_\a$ satisfies
\[\lb{specH0}
\s(H_\a)=\bigcup_{k\in\T^d}\s\big(H_\a(k)\big).
\]
For each $k\in\T^d$ the Floquet (or \emph{fiber}) operator $H_\a(k)$ on $\ell^2(\cV_*)$ is given by
\[\label{Hvt'}
H_\a(k)=-\D_\a(k)+V, \qqq \D_\a(k)=\vk-A_\a(k).
\]
Here $V$ and $\vk$ are the electric and degree potentials on $\ell^2(\cV_*)$; $\D_\a(k)$ is the fiber magnetic Laplacian, and $A_\a(k)$ is the fiber magnetic adjacency operator having the form \cite{KS17}
\[
\label{fado}
\big(A_\a(k)f\big)_x=\sum_{\be=(x,y)\in\cA_*}e^{i(\a(\be)+\lan\t(\be),\,k\ran)}f_y,
\qqq f\in\ell^2(\cV_*),\qqq x\in \cV_*,
\]
where $\t(\be)$ is the index of the edge $\be\in\cA_*$ defined by \er{in}, \er{dco}.

\begin{remark}\lb{RFD}
\emph{i}) The Floquet transform $U:\ell^2(\cV)\to \mH$ defined by \er{DFT} is a Fourier transform with respect to the action of the lattice $\G$ on the vertex set $\cV$.

\emph{ii}) The identities \er{did}, \er{did1} are the \emph{direct integral decomposition} of the magnetic Schr\"odin\-ger operator $H_\a$:
$$
\mH=\int_{\T^{d}}^{\os}\ell^2(\cV_*)\,{dk \/(2\pi)^d}\,, \qqq UH_\a U^{-1}=\int^\oplus_{\T^d}H_\a(k){dk\/(2\pi)^d}\,, \qqq
\T^d=\R^d/(2\pi\Z)^d.
$$

\emph{iii}) We can consider the fiber operator $H_\a(k)$, $k\in
\T^d$, as the magnetic Schr\"odinger operator defined by \er{Sh} -- \er{ALO} with the magnetic potential $\wt\a(\be)=\a(\be)+\lan\t(\be),k\ran$, $\be\in\cA_*$, on the fundamental graph $\cG_*=(\cV_*,\cA_*)$. The fiber operator $H_\a(0)$ is just the magnetic Schr\"odinger operator with the magnetic potential $\a$ on $\cG_*$.
\end{remark}

Let $\#M$ denote the number of elements in a set $M$. Each fiber operator $H_\a(k)$, $k\in\T^{d}$, acts on the space $\ell^2(\cV_*)=\C^\n$, $\n=\#\cV_*$, and has $\n$ real eigenvalues $\l_{\a,j}(k)$, $j=1,\ldots,\n$, labeled in non-decreasing order by
$$
\l_{\a,1}(k)\leq\l_{\a,2}(k)\leq\ldots\leq\l_{\a,\nu}(k), \qqq
\forall\,k\in\T^{d},
$$
counting multiplicities. Since the operator $H_\a(k)$ is self-adjoint and analytic in $k\in\T^d$, each $\l_{\a,j}(\cdot)$ is continuous and piecewise real analytic function on the torus $\T^{d}$ (see, e.g., \cite{Ka95}) and defines the \emph{spectral band} $\s_j(H_\a)$ as the range of $\l_{\a,j}(\cdot)$:
\[\lb{ban.1H}
\begin{array}{l}
\displaystyle\s_j(H_\a)=[\l_{\a,j}^-,\l_{\a,j}^+]=\l_{\a,j}(\T^{d}),\qqq j\in\N_\n, \qqq \N_\n=\{1,\ldots,\n\},\\[6pt]
\displaystyle\textrm{where}\qqq \l_{\a,j}^-=\min_{k\in\T^d}\l_{\a,j}(k),\qqq \l_{\a,j}^+=\max_{k\in\T^d}\l_{\a,j}(k).
\end{array}
\]
Some of $\l_{\a,j}(\cdot)$ may be constant, i.e., $\l_{\a,j}(\cdot)=\L_j=\const$, on some subset of $\T^d$ of positive Lebesgue measure. In this case the magnetic Schr\"odinger operator $H_\a$ on $\cG$ has the eigenvalue $\L_j$ of infinite multiplicity. We call $\{\L_j\}$ a \emph{flat band}. Thus, the spectrum of the magnetic Schr\"odinger operator $H_\a$ on the periodic graph $\cG$ has the form
\[\lb{specH}
\s(H_\a)=\bigcup_{k\in\T^d}\s\big(H_\a(k)\big)=
\bigcup_{j=1}^{\nu}\s_j(H_\a)=\s_{ac}(H_\a)\cup \s_{fb}(H_\a),
\]
where $\s_{ac}(H_\a)$ is the absolutely continuous spectrum, which is a
union of non-degenerate bands from \er{ban.1H}, and $\s_{fb}(H_\a)$ is
the set of all flat bands (eigenvalues of infinite multiplicity).

\medskip

The paper is organized as follows. In Section \ref{Sec2} we formulate our main results:

\smallskip

$\bu$ We determine trace formulas for the magnetic Schr\"odinger operators $H_\a$ on periodic graphs $\cG$ (Theorem \ref{TPG}). The traces of the fiber operators $H_\a(k)$ are expressed as finite Fourier series of the quasimomentum $k$. The coefficients of the Fourier series are given in terms of the magnetic fluxes, electric potentials and cycles in the fundamental graph $\cG_*$ from some specific cycle sets.

$\bu$ We apply the trace formulas to obtain the lower estimates of the total bandwidth for the magnetic Schr\"odinger operator $H_\a$ in terms of geometric parameters of the graph, magnetic fluxes and electric potentials (Theorem \ref{TVECS}) and formulate sufficient conditions on the magnetic fluxes and geometry of the periodic graph $\cG$ under which the spectrum of $H_\a$ on $\cG$ has at least one non-degenerate band (Theorem \ref{TEsFK}).

Section \ref{Sec3} is devoted to the trace formulas. Here we prove Theorem \ref{TPG} and Corollary \ref{TCo0} about trace formulas for the magnetic Schr\"odinger operators and compare traces of the Schr\"odinger operators with and without magnetic fields (Corollary \ref{CCrFA}). In this section we also present the trace formulas for the magnetic adjacency operator $A_\a$ defined by \er{ALO} (Theorem \ref{TFao}). In Section \ref{Sec4} we prove Theorems \ref{TVECS} and \ref{TEsFK} about the lower bounds of the total bandwidth for the magnetic Schr\"odinger operators.

\section{Main results}
\setcounter{equation}{0}
\lb{Sec2}
\subsection{Magnetic fluxes, cycle indices and cycle sets}\lb{SoC} In order to
formulate trace formulas for the magnetic Schr\"odinger operators we need some notation. An \emph{oriented path} (or, for short, just a \emph{path}) $\bp$ in a graph $\cG=(\cV,\cA)$ is a sequence of oriented edges
\[\lb{depa}
\bp=(\be_1,\be_2,\ldots,\be_n), \qqq \textrm{where} \qq
\be_s=(x_{s-1},x_s)\in\cA,\qq s\in \N_n,
\]
for some vertices $x_0,x_1,\ldots,x_n\in\cV$. The vertices $x_0$ and $x_n$ are called the \emph{initial} and \emph{terminal} vertices of the path $\bp$, respectively. The \emph{reverse} of the path $\bp$ given by \er{depa} is the path $\ul\bp=(\ul\be_n,\ldots,\ul\be_1)$. If $x_0=x_n$, then the path $\bp$ is
called a \emph{cycle}. The number $n$ of edges in a cycle $\bc$ is
called the \emph{length} of $\bc$ and is denoted by $|\bc|$, i.e.,
$|\bc|=n$. A cycle $\bc=(\be_1,\ldots,\be_n)$ has a \emph{back-tracking}, if $\be_{s+1}=\ul\be_s$ for some $s\in\N_n$ ($\be_{n+1}$ is understood as $\be_1$).

\begin{remark} A path $\bp$ is uniquely defined by the sequence of
its oriented edges $(\be_1,\be_2,\ldots,\be_n)$. The sequence of
its vertices $(x_0,x_1,\ldots,x_n)$ does not uniquely define~$\bp$,
since multiple edges are allowed in the graph $\cG$.
\end{remark}

Let $\cC$ be the set of all cycles of the fundamental graph $\cG_*=(\cV_*,\cA_*)$. For any cycle $\bc\in\cC$ we define

$\bu$ the \emph{cycle index} $\t(\bc)\in\Z^d$ by
\[\lb{cyin}
\t(\bc)=\sum\limits_{\be\in\bc}\t(\be),  \qqq \bc\in\cC,
\]
where $\t(\be)\in\Z^d$ is the index of the edge $\be\in\cA_*$ defined by \er{in}, \er{dco};

$\bu$  the \emph{flux} $\a(\bc)\in[-\pi,\pi]$ of the magnetic potential
 $\a$ through the cycle $\bc$ by
\[\lb{cyin2}
\a(\bc)=\bigg(\sum\limits_{\be\in\bc}\a(\be)\bigg)\hspace{-3mm} \mod 2\pi,
 \qqq \qqq \bc\in\cC.
\]
From these definitions and the identities \er{inin0}, \er{inin} it follows that
\[\lb{cyin1}
\t(\ul\bc\,)=-\t(\bc), \qqq \a(\ul\bc\,)=-\a(\bc), \qqq  \,\forall\,\bc\in\cC.
\]

\begin{remark}\label{Re22}
Edge indices depend on the choice of the embedding
of the periodic graph $\cG$ into $\R^d$.
Cycle indices \emph{do not} depend on this choice. This property
of cycle indices has the following simple explanation. Any cycle $\bc$ in the fundamental graph $\cG_*$ is obtained by factorization of a path in the periodic graph $\cG=(\cV,\cA)$ connecting some $\G$-equivalent vertices $x\in\cV$ and $x+\ga\in\cV$, $\ga\in\G$. Furthermore, the index of the cycle $\bc$ is equal to $\mm=(m_j)_{j=1}^d\in\Z^d$, where $\ga=m_1\ga_1+\ldots+m_d\ga_d$, and $\{\ga_1,\ldots,\ga_d\}$ is the basis of the lattice $\G$ and, therefore, does not depend on
the choice of the embedding. In particular, $\t(\bc)=0$ if and only if the cycle $\bc$ in $\cG_*$ corresponds to a cycle in~$\cG$.
\end{remark}

For the set  $\cC$ of all cycles of the fundamental graph $\cG_*$  we define the following subsets of $\cC$ which will be used throughout this paper:\\
$\bu$  $\cC_n$ is the set of all cycles of length $n$:
\[\lb{cNn}
\cC_{n}=\{\bc\in\cC: |\bc|=n \};
\]
$\bu$  $\cC_n^\mm$ is the set of all cycles of length $n$ and with index $\mm$:
\[
\lb{cNnm} \cC_n^\mm=\{\bc\in\cC: |\bc|=n \;\textrm{and}\;\t(\bc)=\mm\}\ss\cC_n;
\]
$\bu$  $\cC_n^+$ is the set of all cycles of length $n$ with non-zero indices:
\[
\lb{cNn+} \cC_n^+=\{\bc\in\cC: |\bc|=n \;\textrm{and}\;\t(\bc)\neq0\}\ss\cC_n;
\]
$\bu$  $\cC_n^{odd}$ is the set of all cycles of length $n$ with odd
sum of index components:
\[
\lb{cNno} \cC_n^{odd}=\{\bc\in\cC: |\bc|=n \;\textrm{and}\;\lan\t(\bc),\1\ran\textrm{ is odd}\}\ss\cC_n,
\]
where $\1=(1,\ldots,1)\in\R^d$. Here and below\\
$\bu$  $|\bc|$ is the length  of the cycle $\bc$;\\
$\bu$  $\t(\bc)\in\Z^d$ is the index of the cycle $\bc$ defined by \er{cyin}.

\begin{remark}\lb{Rccy}
The sets $\cC_n,\cC_n^\mm,\cC_n^+,\cC_n^{odd}$ are finite and include the corresponding cycles with back-tracking parts.
\end{remark}

\begin{example1} We obtain the cycle sets $\cC_n,\cC_n^\mm,\cC_n^+,\cC_n^{odd}$ (for some $n$ and $m$) for the fundamental graphs from Example~\ref{ExPG}.\emph{i},\emph{ii}.

\emph{i}) The fundamental graph $\bS_*$ of the square lattice $\bS$ labeled by edge indices is shown in Fig.\ref{FSqL}\emph{b}. All cycles of length one in $\bS_*$ are the loops $\be_1,\be_2$ with indices $\t(\be_1)=(1,0)$, $\t(\be_2)=(0,1)$, and their reverse loops $\ul\be_1,\ul\be_2$. Thus, $\cC_1=\{\be_1,\be_2,\ul\be_1,\ul\be_2\}$,
$$
\cC_1^{(1,0)}=\{\be_1\},\qqq \cC_1^{(0,1)}=\{\be_2\},\qqq
\cC_1^{(-1,0)}=\{\ul\be_1\},\qqq \cC_1^{(0,-1)}=\{\ul\be_2\},
$$
and $\cC_1^+=\cC_1^{odd}=\cC_1$.

\emph{ii}) The fundamental graph $\bG_*=(\cV_*,\cA_*)$ of the hexagonal lattice $\bG=(\cV,\cA)$ labeled by edge indices is shown in Fig.\ref{ff.0.3}\emph{b}. There are no cycles of length one (i.e., loops) and length three in $\bG_*$. All cycles of length two without back-tracking in $\bS_*$ are the cycles
$$
\bc_1=(\be_1,\be_2),\qqq \bc_2=(\be_1,\be_3), \qqq \bc_3=(\be_2,\ul\be_3),
$$
their permutations
$$
\wt\bc_1=(\be_2,\be_1),\qqq \wt\bc_2=(\be_3,\be_1), \qqq \wt\bc_3=(\ul\be_3,\be_2),
$$
and their reverse cycles $\ul\bc_s,\ul{\wt\bc}_s$, $s=1,2,3$. The indices of these cycles are given by
$$
\begin{array}{l}
\t(\bc_1)=\t(\wt\bc_1)=\t(\be_1)+\t(\be_2)=(1,0), \\[4pt]
\t(\bc_2)=\t(\wt\bc_2)=\t(\be_1)+\t(\be_3)=(0,1),\\[4pt] \t(\bc_3)=\t(\wt\bc_3)=\t(\be_2)+\t(\ul\be_3)=\t(\be_2)-\t(\be_3)=(1,-1),
\end{array}
$$
$\t(\ul\bc_s)=-\t(\bc_s)$, $\t(\ul{\wt\bc}_s)=-\t(\wt\bc_s)$, $s=1,2,3$.
Each oriented edge $\be\in\cA_*$ gives rise to a backtracking cycle $\bc_\be=(\be,\ul\be\,)$ of length two with zero index. Thus,
$$
\begin{array}{l}
\cC_2^{(0,0)}=\{\bc_\be \mid \be\in\cA_*\},\\[4pt]
\cC_2^{(1,0)}=\{\bc_1,\wt\bc_1\},\qqq\qq \cC_2^{(0,1)}=\{\bc_2,\wt\bc_2\},\qqq\qq \cC_2^{(1,-1)}=\{\bc_3,\wt\bc_3\},\\[4pt]
\cC_2^+=\{\bc_s,\wt\bc_s,\ul\bc_s,\ul{\wt\bc}_s \mid s=1,2,3\},\qqq\qq \cC_2^{odd}=\{\bc_s,\wt\bc_s,\ul\bc_s,\ul{\wt\bc}_s \mid s=1,2\}.
\end{array}
$$
\end{example1}

\subsection{Trace formulas for magnetic Schr\"odinger operators} In order to determine trace formulas for the magnetic Schr\"odinger operators $H_\a=-\D_\a+V$, we need some modification of the fundamental graph $\cG_*$. At each vertex $x$ of $\cG_*$ we add a loop with weight $v_x=V_x-\vk_x$. This modification allows us to replace the fiber magnetic Schr\"odinger operator $H_\a(k)$ acting on the fundamental graph $\cG_*$ by a fiber weighted magnetic adjacency operator on the modified fundamental graph.

More precisely, at each vertex $x$ of $\cG_*=(\cV_*,\cA_*)$ we add a loop  $\be_x$ with index $\t(\be_x)=0$ and magnetic potential $\a(\be_x)=0$, and consider the modified fundamental graph $\wt\cG_*=(\cV_*,\wt\cA_*)$, where
\[\lb{wtAs}
\wt\cA_*=\cA_*\cup\{\be_x\}_{x\in\cV_*}, \qqq \t(\be_x)=0, \qqq \a(\be_x)=0.
\]

Let $\wt\cC$ be the set of all cycles in $\wt\cG_*$. For each cycle
$$
\bc=(\be_1,\ldots,\be_n)\in\wt\cC, \qqq  \be_s=(x_s,x_{s+1})\in\wt\cA_*, \qq s\in \N_n, \qq x_{n+1}=x_1,
$$
we define the \emph{weight} $\o(\bc)$ by
\[\lb{Wcy} \o(\bc)=\o(\be_1)\ldots\o(\be_n),
\]
where $\o(\be)$ has the form
\[\lb{webe}
\o(\be)=\left\{
\begin{array}{cl}
1,  & \qq \textrm{if} \qq  \be\in\cA_* \\[2pt]
v_x, & \qq \textrm{if} \qq \be=\be_x
\end{array}\right.,\qqq v_x=V_x-\vk_x,
\]
and $\vk_x$ is the degree of the vertex $x\in\cV_*$. From this definition it follows that
\[\lb{omin}
\o(\ul\bc\,)=\o(\bc), \qqq \forall\,\bc\in\wt\cC.
\]

\begin{remark}
Note that
\[\lb{oc1C}
\o(\bc)=1\qqq \textrm{for each cycle}\qqq \bc\in\cC.
\]
\end{remark}

Let the cycle sets $\wt\cC_n$ and $\wt\cC_n^\mm$ for the modified fundamental graph $\wt\cG_*$ be defined as the corresponding cycle sets $\cC_n$ and $\cC_n^\mm$ for the fundamental graph $\cG_*$, see \er{cNn} and \er{cNnm}.

We formulate trace formulas for the fiber magnetic Schr\"odinger operators.
Let
\[\lb{tapl}
\t_+=\max\limits_{\be\in\cA_*}\|\t(\be)\|,
\]
where $\t(\be)\in\Z^d$ is the index of the edge $\be$ defined by \er{in}, \er{dco}, and $\|\cdot\|$ is the standard norm in $\R^d$. Recall that
\begin{itemize}
  \item $\t(\bc)$ is the index of the cycle $\bc$ defined by \er{cyin};
  \item $\a(\bc)$ is the flux of the magnetic potential $\a$ through the cycle $\bc$ defined by \er{cyin2};
  \item $\o(\bc)$ is the weight of the cycle $\bc$ defined by \er{Wcy}.
\end{itemize}

\begin{theorem}\lb{TPG}
Let $H_\a(k)$, $k\in\T^d$, be the fiber magnetic Schr\"odinger operator defined by \er{Hvt'} -- \er{fado} on the fundamental graph $\cG_*=(\cV_*,\cA_*)$. Then for each $n\in\N$

i) The trace of $H_\a^n(k)$ has the form
\[\lb{TrnH}
\Tr H_\a^n(k)=\sum_{j=1}^\n\l_{\a,j}^n(k)=\cT_{\a,n}(k), \qqq \n=\#\cV_*,
\]
where
\[
\lb{deTnk}
\begin{aligned}
&\cT_{\a,n}(k)=\sum_{\bc\in\wt\cC_n}\o(\bc)\,e^{-i\f(\bc,k)}=
\sum_{\bc\in\wt\cC_n}\o(\bc)\cos\f(\bc,k), \\ &\f(\bc,k)=\a(\bc)+\lan\t(\bc),k\ran,
\end{aligned}
\]
or, in the form of Fourier series,
\[\lb{FsTrH}
\cT_{\a,n}(k)=\sum_{\mm\in\Z^d\atop \|\mm\|\leq n\t_+}
\cT_{\a,n,\mm}e^{-i\lan\mm,k\ran},\qqq
\cT_{\a,n,\mm}=\sum_{\bc\in\wt\cC_n^\mm}\o(\bc)e^{-i\a(\bc)},
\]
$\t_+$ is defined by \er{tapl}.

ii) The trace of $H_\a^n(k)$ satisfies
\[\lb{ITrH}
\frac1{(2\pi)^d}\int_{\T^d}\Tr H_\a^n(k)dk=\cT_{\a,n,0},
\]
where
\[\lb{IVFu}
\cT_{\a,n,0}=\sum_{\bc\in\wt\cC_n^0}\o(\bc)e^{-i\a(\bc)}=
\sum_{\bc\in\wt\cC_n^0}\o(\bc)\cos\a(\bc).
\]
\end{theorem}

\begin{remark}
\emph{i}) The formulas \er{TrnH} -- \er{FsTrH} and \er{ITrH}, \er{IVFu} are \emph{trace formulas}, where the traces of the fiber operators are expressed in terms of geometric parameters of the graph (vertex degrees, cycle
indices and lengths), magnetic fluxes and electric potentials. The explicit form for the identities \er{TrnH}, \er{deTnk} and \er{ITrH}, \er{IVFu} when $n=1$ and $n=2$ is given in Proposition~\ref{TrTO}.

\emph{ii}) The functional $\cT_{\a,n}(k)$ is a finite Fourier series of the quasimomentum $k\in\T^d$, and $\o(\bc)$ are polynomials of degree $\le n$, $n=|\bc|$,  with respect to the electric potential $V$.

\emph{iii}) We sometimes write $\cT_{\a,n}(k,V),\cT_{\a,n,\mm}(V),\o(\bc,V)$ instead of $\cT_{\a,n}(k),\cT_{\a,n,\mm},\o(\bc)$, when several electric potentials $V$ are considered.

\emph{iv}) The trace formulas for the magnetic Laplacian $-\D_\a=A_\a-\vk$ are given by the same identities \er{TrnH} -- \er{FsTrH} and \er{ITrH}, \er{IVFu}, where $\cT_{\a,n}(k)=\cT_{\a,n}(k,0)$ and $\cT_{\a,n,0}=\cT_{\a,n,0}(0)$.

\emph{v}) The traces $\Tr H_\a^n(k)$, $(n,k)\in\N\ts\T^d$, given by \er{TrnH}, \er{deTnk}, depend on the magnetic fluxes $\a(\bc)$ through the cycles $\bc$ of the fundamental graph $\cG_*$. Therefore, the effect of the magnetic field on the traces $\Tr H_\a^n(k)$ is fully determined by $\b$ fluxes through a chosen set of basis cycles of $\cG_*$. Here and below
$$
\b=\#\cE_*-\#\cV_*+1,
$$
is the Betti number of the fundamental graph $\cG_*=(\cV_*,\cE_*)$, i.e., the number of independent cycles in $\cG_*$. Using a shift of the quasimomentum $k\in\T^d$, we can reduce the number of these independent fluxes to $\b-d$, where $d$ is the dimension of the periodic graph (see Corollary \ref{TCo0} and Example \ref{EKL}). In particular, if $\b=d$, then the traces $\Tr H_\a^n(k)$, $n\in\N$, and, consequently, the spectrum of the magnetic Schr\"odinger operator $H_\a$ do not depend on the magnetic fluxes.

\emph{vi}) For $\a=0$ the trace formulas \er{TrnH} -- \er{FsTrH} and \er{ITrH}, \er{IVFu} coincide with trace formulas for the Schr\"odinger operators $H_0$ without magnetic fields presented in \cite{KS22}.
\end{remark}

\subsection{Estimates of the total bandwidth for magnetic Schr\"odinger operators} \lb{S2.3} We apply the obtained trace formulas \er{TrnH} -- \er{IVFu} to estimate the total bandwidth for the magnetic Schr\"odinger operator $H_\a$. In order to present our results we describe known total bandwidth estimates for $H_\a$. There are only \emph{upper} estimates determined in \cite{KS17,KS18}. We do not know lower estimates.

Define the total bandwidth $\gS(H_\a)$ of the operator $H_\a$ by
$$
\gS(H_\a)= \sum_{j=1}^{\n}|\s_j(H_\a)|, \qqq \n=\#\cV_*,
$$
where $\s_j(H_\a)$ is the spectral band of $H_\a$ defined by \er{ban.1H}.

In \cite{KS17} the authors estimated the total bandwidth $\gS(H_\a)$ in terms of the Betti number $\b$ of the fundamental graph $\cG_*=(\cV_*,\cE_*)$:
\[
\lb{esLm1} \gS(H_\a)\le 4\b,\qqq \b=\#\cE_*-\#\cV_*+1.
\]
In \cite{KS18} the estimate \er{esLm1} was improved:
\[
\lb{ues1} \gS(H_\a)\leq 4\cI.
\]
Here $\cI$ is an invariant of the periodic graph $\cG$ introduced in \cite{KS20}:
\[\lb{dIm}
\textstyle\cI=\frac12\,\min\limits_{\cG\ss\R^d}\#\cB_*, \qqq \cB_*=\supp\t,
\]
i.e., $\cI$ is the minimal number of the fundamental graph edges $\be\in\cE_*$ with non-zero indices $\t(\be)$, where the minimum is taken over all embeddings of $\cG$ into $\R^d$.

\begin{remark}
\emph{i}) The invariant $\cI$ defined by \er{dIm} satisfies $\cI\leq\b$ and the difference between the Betti number $\b$ and the invariant $\cI$ may be any non-negative integer number (for specific graphs). For more details see Theorem 2.1.\emph{v} and Proposition 2.2 in \cite{KS20}.

\emph{ii}) The upper estimates \er{esLm1} and \er{ues1} of the total bandwidth for $H_\a=-\D_\a+V$ do not depend on the magnetic potential $\a$ and the electric potential $V$.

\emph{iii}) Estimates for the total bandwidth of Schr\"odinger operators  without magnetic fields were discussed in \cite{KS14,KS19,KS20,KS22a}. In particular, in \cite{KS22a} it was shown that
\[
\lb{CEC1} \frac{2d_*}{v_*^{\n-1}}\leq \gS(H_0)\leq 4\cI,\qqq d_*=\left\{\begin{array}{cl}
              d, & \textrm{if $d$ is even} \\
              d+1, & \textrm{if $d$ is odd}
            \end{array}\right.,
\]
where $v_*$ and $\cI$ are defined in \er{vsdv} and \er{dIm}, respectively.
\end{remark}

Now we present lower estimates of the total bandwidth for the magnetic Schr\"odinger operators $H_\a=-\D_\a+V$. Let the cycle sets $\wt\cC_n^+$ and $\wt\cC_n^{odd}$ for the modified fundamental graph $\wt\cG_*$ be defined as the corresponding cycle sets $\cC_n^+$ and $\cC_n^{odd}$ for the fundamental graph $\cG_*$, see \er{cNn+} and \er{cNno}.

We define the diameter of the electric potential $V$ by
$$
\diam V=\max_{x\in\cV_*}V_x-\min_{x\in\cV_*}V_x.
$$
Recall that
\begin{itemize}
  \item $\a(\bc)$ is the flux of the magnetic potential $\a$ through the cycle $\bc$ defined by \er{cyin2};
  \item $\o(\bc)$ is the weight of the cycle $\bc$ defined by \er{Wcy}.
\end{itemize}

\begin{theorem}\lb{TVECS}
Let $H_\a=-\D_\a+V$ be the Schr\"odinger operator defined by \er{Sh} -- \er{ALO}  with a periodic magnetic potential $\a$ and a periodic electric potential $V$ on a periodic graph $\cG$. Then its total bandwidth $\gS(H_\a)$ satisfies
\[\lb{leHa}
\gS(H_\a)\geq\frac{\max\big\{B_{\a,n}^+,2B_{\a,n}^{odd}\big\}}{nv_*^{n-1}}\,,\qqq \forall\, n\in\N,
\]
where
\[\lb{Banpo}
B_{\a,n}^+=\Big|\sum_{\bc\in\wt\cC_n^+}\o(\bc)\cos\a(\bc)\Big|\,,\qqq
B_{\a,n}^{odd}=\Big|\sum_{\bc\in\wt\cC_n^{odd}}\o(\bc)\cos\a(\bc)\Big|\,,
\]
and \[\lb{vsdv}
v_*=\diam V+\vk_+,\qqq \vk_+=\max\limits_{x\in\cV_*}\vk_x,
\]
$\vk_x$ is the degree of the vertex $x\in\cV_*$.

In particular, if $n$ is the length of the shortest cycle with non-zero index in the fundamental graph $\cG_*$, then
\[\lb{Banp1}
B_{\a,n}^+=\Big|\sum_{\bc\in\cC_{n}^+}\cos\a(\bc)\Big|\,,\qqq
B_{\a,n}^{odd}=\Big|\sum_{\bc\in\cC_{n}^{odd}}\cos\a(\bc)\Big|.
\]
\end{theorem}

\begin{remark}
\emph{i}) Since $\wt\cC_{n}^{odd}\subseteq\wt\cC_{n}^+$, each term of the sum for $B_{\a,n}^{odd}$, see \er{Banpo}, is also in the sum for $B_{\a,n}^+$, $n\in\N$. But, the constant $2B_{\a,n}^{odd}$ may be greater than $B_{\a,n}^+$.

\emph{ii}) For $\a=0$ the estimate \er{leHa}, \er{Banp1} coincides with the estimate for the total bandwidth of the Schr\"odinger operator $H_0$ without magnetic fields obtained in \cite{KS22a}:
\[\lb{esa0}
\mathfrak{S}(H_0)\geq\frac{\max\big\{\cN_n^+,2\cN_n^{odd}\big\}}
{nv_*^{n-1}}\,,\qqq \cN_n^+=\#\cC_n^+,\qqq \cN_n^{odd}=\#\cC_n^{odd},
\]
where $n$ is the length of the shortest cycle with non-zero index in $\cG_*$, and $\cC_n^+$, $\cC_n^{odd}$ are defined in \er{cNn+}, \er{cNno}. For the Schr\"odinger operator $H_0=-\D_0+V$ with a real $\n$-periodic potential $V$ on the one-dimensional lattice $\Z$ the estimate \er{esa0} coincides with the Last's estimate
$$
\gS(H_0)\geq\frac{4}{v_*^{\n-1}}\,,\qqq \textrm{where}\qqq v_*=2+\diam V,
$$
see \cite{L92}, which is sharp.

\emph{iii}) The problem of estimating the measure of the spectrum for the Harper operator $\D_\f$, see \er{haop}, was studied by several authors \cite{AMS90,HK95,JKK22,L94}. For irrational ${\f\/2\pi}$ the spectrum $\s(\D_\f)$ of $\D_\f$ is a zero measure Cantor set. For rational ${\f\/2\pi}=\frac pq$, where $p,q$ are coprime positive integers, the Lebesgue measure $|\s(\D_\f)|$ of $\s(\D_\f)$ satisfies
\[\lb{esHa}
\textstyle \frac{2(\sqrt{5}+1)}q<|\s(\D_\f)|<\frac{4\pi}q\,,
\]
see Lemma 1 in \cite{L94} and Theorem 1 in \cite{JKK22}. The proof of these estimates was based on the decomposition of the Harper operator $\D_\f$ into the direct integral in $k$ of the one-dimensional Schr\"odinger operator with a specific electric potential given by
$$
(H_{\f,k}f)_n=f_{n+1}+f_{n-1}+2\cos(\f n+k)f_n, \qq f\in\ell^2(\Z),\qq \f,k\in[0,2\pi).
$$
Similar results for the Harper model on the hexagonal lattice were obtained in \cite{BHJ19}. We do not know results about the Lebesgue measure of the spectrum of the Harper type models for periodic graphs of dimension $d\geq3$. For these cases the problem is not reduced to the one-dimensional case.
\end{remark}

For some magnetic potentials $\a$ the constants $B_{\a,n}^+$ and $B_{\a,n}^{odd}$ defined by \er{Banpo} may be zero, see \er{coBB}. Thus, in general, Theorem \ref{TVECS} does not guarantee that the absolutely continuous spectrum of the magnetic Schr\"odinger operator $H_\a$ is not empty. In the next theorem we formulate one more lower bound for the total bandwidth of $H_\a$ and give some simple sufficient conditions on the magnetic fluxes and geometry of the periodic graph $\cG$ under which the spectrum of $H_\a$ on $\cG$ has at least one non-degenerate band.

\begin{theorem}\lb{TEsFK}
Let $H_\a=-\D_\a+V$ be the Schr\"odinger operator defined by \er{Sh} -- \er{ALO} with a periodic magnetic potential $\a$ and a periodic electric potential $V$ on a periodic graph $\cG$. Then its total bandwidth $\gS(H_\a)$ satisfies
\[\lb{EsFS}
\gS(H_\a)\geq\frac{2\big|\cT_{\a,n,\mm}\big|}{nv_*^{n-1}}\,,\qqq \forall\, (n,\mm)\in\N\ts(\Z^d\sm\{0\}),
\]
where $\cT_{\a,n,\mm}$ are the Fourier coefficients defined in \er{FsTrH}, and $v_*$ is given by \er{vsdv}.

In particular, let $\mm\in\Z^d\sm\{0\}$ and $\frac1q\,\mm\not\in\Z^d$ for any integer $q\geq2$, and let $n(\mm)$ be the length of the shortest cycle with index $\mm$ in the fundamental graph $\cG_*$. Assume that $\bc_1,\ldots,\bc_p$ (and their cyclic edge permutations) are all cycles of length $n(\mm)$ and with index $\mm$ in $\cG_*$, then
\[\lb{EsFS1}
\big|\cT_{\a,n(\mm),\mm}\big|=
n(\mm)\Big|1+\sum_{j=2}^pe^{-i\a(\bc_j-\bc_1)}\Big|\,.
\]
Moreover,

i) if $p=1$, then $\big|\cT_{\a,n(\mm),\mm}\big|=n(\mm)>0$;

ii) if $p=2$, then $\big|\cT_{\a,n(\mm),\mm}\big|=2n(\mm)\big|\cos\frac{\a(\bc_2-\bc_1)}2\big|$;

iii) if $p\geq2$ and $\a_+=\max_{2\leq j\leq p}\big|\a(\bc_j-\bc_1)\big|<\pi/2$, then
$$
\big|\cT_{\a,n(\mm),\mm}\big|\geq n(\mm)(1+(p-1)\cos\a_+)>0.
$$
In the cases i) and iii) and in the case ii) for $\a(\bc_2-\bc_1)\neq\pi$ the magnetic Schr\"odinger operator $H_\a$ has at least one non-degenerate band.
\end{theorem}

\begin{remark}
\emph{i}) Let $\bc$ be a cycle in the fundamental graph $\cG_*$ with index $\mm\in\Z^d\sm\{0\}$. The condition $\frac1q\,\mm\not\in\Z^d$ for any integer $q\geq2$ means that $\bc$ is a \emph{prime} cycle, i.e., it is not obtained by repeating $q$-times any cycle in $\cG_*$.

\emph{ii}) In contrast to the upper estimates \er{esLm1} and \er{ues1}, the lower estimates \er{leHa} and \er{EsFS} of the total bandwidth for $H_\a$ do depend on the magnetic potential~$\a$.

\emph{iii}) For fixed $n\in\N$ the lower bound \er{leHa} is expressed in terms of \emph{all} cycles of length $n$ with non-zero indices in the modified fundamental graph $\wt\cG_*$, whereas the lower bound \er{EsFS} uses only cycles of length $n$ \emph{with some fixed index} $\mm\in\Z^d\sm\{0\}$. Thus, in general, the estimate \er{leHa} is better than \er{EsFS}. But the last one gives simple sufficient conditions under which the spectrum of the magnetic Schr\"odinger operator $H_\a$ has an absolutely continuous component, see Theorem~\ref{TEsFK}.\emph{i} -- \emph{iii}.

\emph{iv}) We think that the condition $p=1$ or $p=2$  for some $\mm\in\Z^d\sm\{0\}$ and $n(\mm)$ (see items \emph{i}, \emph{ii} of Theorem~\ref{TEsFK}) is quite general (at least for "minimal"\, fundamental graphs). Thus, in general case, \er{EsFS} gives a non-zero lower bound for the total bandwidth of $H_\a$.

\emph{v}) The estimates for the magnetic Laplacian $-\D_\a=A_\a-\vk$ are given by the same inequalities \er{leHa} -- \er{vsdv} and \er{EsFS}, where $\o(\bc)=\o(\bc,0)$, $\cT_{\a,n,\mm}=\cT_{\a,n,\mm}(0)$ and $v_*=\vk_+$ (since the electric potential $V=0$).

\emph{vi}) All spectral bands of $H_\a$ are degenerate, i.e., $\gS(H_\a)=0$, iff $\cT_{\a,n,\mm}=0$ for all $(n,\mm)\in\N_\n\ts\big(\Z^d\sm\{0\}\big)$, see Lemma \ref{Lcond}.
\end{remark}

\subsection{An example and flat band spectrum} \lb{Exe} We apply the obtained estimates \er{leHa} and \er{EsFS} of the total bandwidth to the Laplacian $\D_\a$ with a periodic magnetic potential $\a$ on the simple $\Z$-periodic graph $\cG$ shown in Fig.\ref{FEx1}\emph{a}. We also show that there exists a magnetic potential $\a$ such that the absolutely continuous spectrum of $\D_\a$ on $\cG$ is empty.

\begin{figure}[h]
\centering
\unitlength 1mm 
\linethickness{0.4pt}
\ifx\plotpoint\undefined\newsavebox{\plotpoint}\fi 
\begin{picture}(100,30)(0,0)

\put(40.0,21.5){$\be_5$}
\put(54.5,21.5){$\be_4$}
\put(41.5,7){$\be_2$}
\put(55.0,7.0){$\be_3$}
\put(37,12.0){$x_1$}
\put(59,12.0){$x_4$}
\put(50,2.0){$x_3$}
\put(49.5,26.0){$x_2$}
\put(67.5,11.5){$\be_1$}
\put(49,14.0){$\bc_0$}

\put(45,10){\vector(1,-1){1.0}}
\put(55,10){\vector(1,1){1.0}}
\put(45,20){\vector(-1,-1){1.0}}
\put(55,20){\vector(-1,1){1.0}}
\put(68,14.9){\vector(1,0){1.0}}

\put(0,15){\line(1,0){5.0}}

\put(15,5){\circle{1}}
\put(5,15){\circle{1}}
\put(25,15){\circle{1}}
\put(15,25){\circle{1}}
\put(15,5){\line(-1,1){10.0}}
\put(15,5){\line(1,1){10.0}}
\put(15,25){\line(-1,-1){10.0}}
\put(15,25){\line(1,-1){10.0}}

\put(25,15){\line(1,0){15.0}}

\put(50,5){\circle*{1}}
\put(40,15){\circle*{1}}
\put(60,15){\circle*{1}}
\put(50,25){\circle*{1}}
\put(50,5){\line(-1,1){10.0}}
\put(50,5){\line(1,1){10.0}}
\put(50,25){\line(-1,-1){10.0}}
\put(50,25){\line(1,-1){10.0}}

\put(60,15){\line(1,0){15.0}}

\put(85,5){\circle{1}}
\put(75,15){\circle{1}}
\put(95,15){\circle{1}}
\put(85,25){\circle{1}}
\put(85,5){\line(-1,1){10.0}}
\put(85,5){\line(1,1){10.0}}
\put(85,25){\line(-1,-1){10.0}}
\put(85,25){\line(1,-1){10.0}}

\put(95,15){\line(1,0){5.0}}

\put(0,3){\emph{a})}
\put(5,22){$\cG$}

\end{picture}\hspace{10mm}
\begin{picture}(30,30)(0,0)
\put(15,5){\circle*{1}}
\put(5,15){\circle*{1}}
\put(25,15){\circle*{1}}
\put(15,25){\circle*{1}}
\put(10,10){\vector(1,-1){1.0}}
\put(20,10){\vector(1,1){1.0}}
\put(15,5){\line(-1,1){10.0}}
\put(15,5){\line(1,1){10.0}}
\put(15,25){\line(-1,-1){10.0}}
\put(15,25){\line(1,-1){10.0}}
\put(10,20){\vector(-1,-1){1.0}}
\put(20,20){\vector(-1,1){1.0}}
\put(5,15){\line(1,0){20.0}}
\put(15.0,14.9){\vector(-1,0){1.0}}
\put(0,3){\emph{b})}
\put(0,24){$\cG_*$}
\put(6,21.5){$\be_5$}
\put(19.5,21.5){$\be_4$}
\put(7.5,7){$\be_2$}
\put(19.5,7.0){$\be_3$}
\put(2,12.0){$x_1$}
\put(25,12.0){$x_4$}
\put(15,2.0){$x_3$}
\put(14.5,26.0){$x_2$}
\put(15,12.0){$\be_1$}
\end{picture}
\caption{\footnotesize  \emph{a}) A $\Z$-periodic graph $\cG$; the influence of the magnetic field on the spectrum of $\D_\a$ on $\cG$ is completely determined by the flux $\phi:=\a(\bc_0)$ through the cycle $\bc_0=(\be_2,\be_3,\be_4,\be_5)$ of $\cG$; \quad \emph{b}) the fundamental graph  $\cG_*$.} \label{FEx1}
\end{figure}
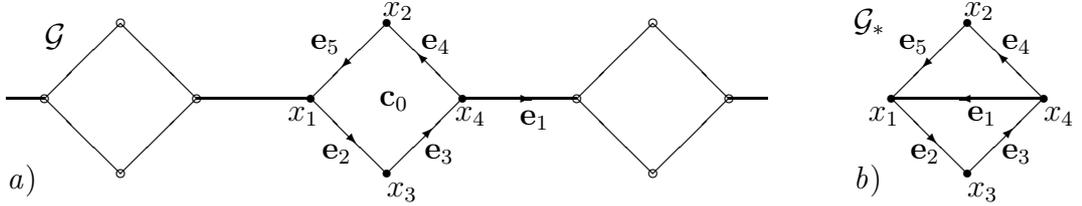

\begin{example}\lb{Exa}
Let $\D_\a$ be the Laplacian with a periodic magnetic potential $\a$ on the $\Z$-periodic graph $\cG$ shown in Fig.\ref{FEx1}a. Then

i) The total bandwidth $\gS(\D_\a)$ of $\D_\a$ satisfies
\[\lb{esti}
\textstyle\max\big\{\frac89\,\cos^2\frac{\phi}2\,,\,\frac{4}{9}\,
\big|\cos\frac{\phi}2\big|\big\}
\leq\gS(\D_\a)\leq4,
\]
where $\phi:=\alpha(\mathbf{c}_0)$ is the flux through the single (up to periodicity) cycle $\mathbf{c}_0=(\mathbf{e}_2,\mathbf{e}_3,\mathbf{e}_4,\mathbf{e}_5)$ of the periodic graph $\mathcal{G}$, see Fig.\ref{FEx1}a. The estimates \er{esti} are sharp.

ii) If $\phi=\pi$, then the spectrum of $\D_\a$ is given by
$$
\s(\D_\a)=\s_{fb}(\D_\a)=\{2\pm\sqrt{2}\,,3\pm\sqrt{3}\,\},
$$
i.e., $\s_{ac}(\D_\a)=\varnothing$.
\end{example}

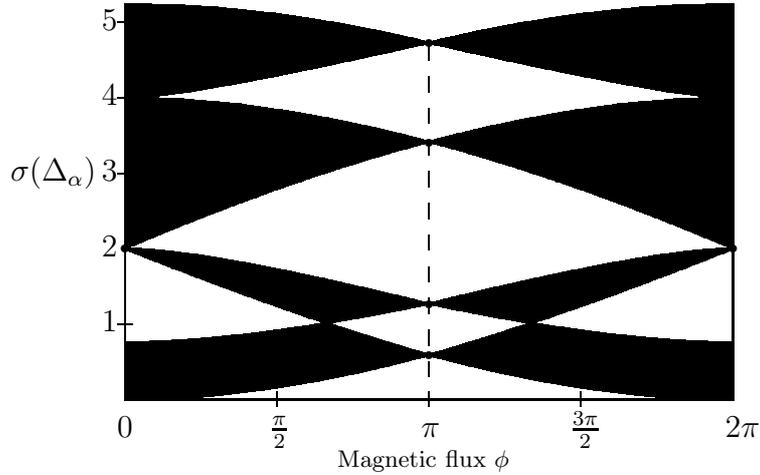
\begin{figure}[h]
\centering
\unitlength 1mm 
\linethickness{0.4pt}
\ifx\plotpoint\undefined\newsavebox{\plotpoint}\fi 
\begin{picture}(80,55)(0,0)

\put(20,-1){\line(0,1){2}}
\put(40,-1){\line(0,1){2}}
\put(60,-1){\line(0,1){2}}
\put(0,0){\line(1,0){80.0}}
\put(-1,10){\line(1,0){2.0}}
\put(-1,30){\line(1,0){2.0}}
\put(-1,40){\line(1,0){2.0}}
\put(-1,50){\line(1,0){2.0}}
\put(0,0){\line(0,1){52.4}}
\put(80,0){\line(0,1){52.4}}
\put(0,20){\circle*{1}}
\put(80,20){\circle*{1}}
\put(40,5.9){\circle*{1}}
\put(40,12.7){\circle*{1}}
\put(40,34.1){\circle*{1}}
\put(40,47.3){\circle*{1}}
\multiput(40,3)(0,4.2){11}{\line(0,1){2.0}}
\put(-1,-5){$0$}
\put(19,-5){$\frac\pi2$}
\put(39,-5){$\pi$}
\put(58.5,-5){$\frac{3\pi}2$}
\put(79,-5){$2\pi$}
\put(-3,9){$1$}
\put(-3,19){$2$}
\put(-3,29){$3$}
\put(-3,39){$4$}
\put(-3,49){$5$}
\put(-15,29){$\s(\D_\a)$}
\put(28,-9){\scriptsize Magnetic flux $\phi$}
\bezier{600}(0,0.0)(18.4,0.0)(40,5.9)
\bezier{600}(80,20.0)(57.6,11.0)(40,5.9)
\bezier{600}(80,0.0)(61.6,0.0)(40,5.9)
\bezier{600}(0,20.0)(22.4,11.0)(40,5.9)

\bezier{600}(0,7.6)(16,7.6)(40,12.7)
\bezier{600}(80,20.0)(69.6,20)(40,12.7)
\bezier{600}(80,7.6)(64,7.6)(40,12.7)
\bezier{600}(0,20.0)(10.4,20)(40,12.7)

\bezier{900}(0,20)(46.4,40)(80,40)
\bezier{900}(0,40)(33.6,40)(80,20)

\bezier{600}(0,40)(8.0,40)(40,47.3)
\bezier{600}(80,52.4)(60.8,52.4)(40,47.3)
\bezier{600}(80,40)(72.0,40)(40,47.3)
\bezier{600}(0,52.4)(19.2,52.4)(40,47.3)

\linethickness{1.5pt}
\put(0.3,0.0){\line(0,1){7.6}}
\put(0.5,0.0){\line(0,1){7.6}}
\put(1.0,0.0){\line(0,1){7.6}}
\put(1.5,0.0){\line(0,1){7.6}}
\put(2.0,0.0){\line(0,1){7.6}}
\put(2.5,0.0){\line(0,1){7.6}}
\put(3.0,0.0){\line(0,1){7.6}}
\put(3.5,0.0){\line(0,1){7.6}}
\put(4.0,0.0){\line(0,1){7.6}}
\put(4.5,0.0){\line(0,1){7.6}}
\put(5.0,0.0){\line(0,1){7.7}}
\put(5.5,0.0){\line(0,1){7.7}}
\put(6.0,0.2){\line(0,1){7.6}}
\put(6.5,0.2){\line(0,1){7.6}}
\put(7.0,0.2){\line(0,1){7.6}}
\put(7.5,0.2){\line(0,1){7.6}}
\put(8.0,0.2){\line(0,1){7.6}}
\put(8.5,0.2){\line(0,1){7.6}}
\put(9.0,0.3){\line(0,1){7.6}}
\put(9.5,0.4){\line(0,1){7.6}}
\put(10.0,0.4){\line(0,1){7.6}}
\put(10.5,0.4){\line(0,1){7.6}}
\put(11.0,0.5){\line(0,1){7.6}}
\put(11.5,0.5){\line(0,1){7.6}}
\put(12.0,0.5){\line(0,1){7.6}}
\put(12.5,0.6){\line(0,1){7.6}}
\put(13.0,0.7){\line(0,1){7.6}}
\put(13.5,0.7){\line(0,1){7.6}}
\put(14.0,0.8){\line(0,1){7.6}}
\put(14.5,0.8){\line(0,1){7.6}}
\put(15.0,0.9){\line(0,1){7.6}}
\put(15.5,1.0){\line(0,1){7.5}}
\put(16.0,1.1){\line(0,1){7.6}}
\put(16.5,1.1){\line(0,1){7.6}}
\put(17.0,1.2){\line(0,1){7.6}}
\put(17.5,1.3){\line(0,1){7.50}}
\put(18.0,1.4){\line(0,1){7.4}}
\put(18.5,1.4){\line(0,1){7.5}}
\put(19.0,1.5){\line(0,1){7.5}}
\put(19.5,1.5){\line(0,1){7.6}}
\put(20.0,1.6){\line(0,1){7.6}}
\put(20.5,1.6){\line(0,1){7.6}}
\put(21.0,1.7){\line(0,1){7.6}}
\put(21.5,1.9){\line(0,1){7.4}}
\put(22.0,1.9){\line(0,1){7.5}}
\put(22.5,2.0){\line(0,1){7.5}}
\put(23.0,2.1){\line(0,1){7.5}}
\put(23.5,2.2){\line(0,1){7.5}}
\put(24.0,2.2){\line(0,1){7.5}}
\put(24.5,2.3){\line(0,1){7.5}}
\put(25.0,2.4){\line(0,1){7.5}}
\put(25.5,2.5){\line(0,1){7.5}}
\put(26.0,2.6){\line(0,1){7.4}}
\put(26.5,2.7){\line(0,1){7.4}}
\put(27.0,2.8){\line(0,1){7.2}}
\put(27.5,2.9){\line(0,1){6.9}}
\put(28.0,3.0){\line(0,1){6.6}}
\put(28.5,3.2){\line(0,1){6.3}}
\put(29.0,3.2){\line(0,1){6.1}}
\put(29.5,3.3){\line(0,1){5.8}}
\put(30.0,3.4){\line(0,1){5.6}}
\put(30.5,3.6){\line(0,1){5.2}}
\put(31.0,3.6){\line(0,1){5.0}}
\put(31.5,3.8){\line(0,1){4.7}}
\put(32.0,3.9){\line(0,1){4.4}}
\put(32.5,4.0){\line(0,1){4.2}}
\put(33.0,4.1){\line(0,1){3.9}}
\put(33.5,4.2){\line(0,1){3.7}}
\put(34.0,4.4){\line(0,1){3.3}}
\put(34.5,4.5){\line(0,1){3.1}}
\put(35.0,4.6){\line(0,1){2.8}}
\put(35.5,4.8){\line(0,1){2.4}}
\put(36.0,4.9){\line(0,1){2.2}}
\put(36.5,5.0){\line(0,1){1.9}}
\put(37.0,5.1){\line(0,1){1.7}}
\put(37.5,5.3){\line(0,1){1.3}}
\put(38.0,5.4){\line(0,1){1.1}}
\put(38.5,5.5){\line(0,1){0.8}}
\put(39.0,5.6){\line(0,1){0.6}}
\put(39.5,5.7){\line(0,1){0.4}}

\put(79.7,0.0){\line(0,1){7.6}}
\put(79.5,0.0){\line(0,1){7.6}}
\put(79.0,0.0){\line(0,1){7.6}}
\put(78.5,0.0){\line(0,1){7.6}}
\put(78.0,0.0){\line(0,1){7.6}}
\put(77.5,0.0){\line(0,1){7.6}}
\put(77.0,0.0){\line(0,1){7.6}}
\put(76.5,0.0){\line(0,1){7.6}}
\put(76.0,0.0){\line(0,1){7.6}}
\put(75.5,0.0){\line(0,1){7.6}}
\put(75.0,0.0){\line(0,1){7.7}}
\put(74.5,0.0){\line(0,1){7.7}}
\put(74.0,0.2){\line(0,1){7.6}}
\put(73.5,0.2){\line(0,1){7.6}}
\put(73.0,0.2){\line(0,1){7.6}}
\put(72.5,0.2){\line(0,1){7.6}}
\put(72.0,0.2){\line(0,1){7.6}}
\put(71.5,0.2){\line(0,1){7.6}}
\put(71.0,0.3){\line(0,1){7.6}}
\put(70.5,0.4){\line(0,1){7.6}}
\put(70.0,0.4){\line(0,1){7.6}}
\put(69.5,0.4){\line(0,1){7.6}}
\put(69.0,0.5){\line(0,1){7.6}}
\put(68.5,0.5){\line(0,1){7.6}}
\put(68.0,0.5){\line(0,1){7.6}}
\put(67.5,0.6){\line(0,1){7.6}}
\put(67.0,0.7){\line(0,1){7.6}}
\put(66.5,0.7){\line(0,1){7.6}}
\put(66.0,0.8){\line(0,1){7.6}}
\put(65.5,0.8){\line(0,1){7.6}}
\put(65.0,0.9){\line(0,1){7.6}}
\put(64.5,1.0){\line(0,1){7.5}}
\put(64.0,1.1){\line(0,1){7.6}}
\put(63.5,1.1){\line(0,1){7.6}}
\put(63.0,1.2){\line(0,1){7.6}}
\put(62.5,1.3){\line(0,1){7.50}}
\put(62.0,1.4){\line(0,1){7.4}}
\put(61.5,1.4){\line(0,1){7.5}}
\put(61.0,1.5){\line(0,1){7.5}}
\put(60.5,1.5){\line(0,1){7.6}}
\put(60.0,1.6){\line(0,1){7.6}}
\put(59.5,1.6){\line(0,1){7.6}}
\put(59.0,1.7){\line(0,1){7.6}}
\put(58.5,1.9){\line(0,1){7.4}}
\put(58.0,1.9){\line(0,1){7.5}}
\put(57.5,2.0){\line(0,1){7.5}}
\put(57.0,2.1){\line(0,1){7.5}}
\put(56.5,2.2){\line(0,1){7.5}}
\put(56.0,2.2){\line(0,1){7.5}}
\put(55.5,2.3){\line(0,1){7.5}}
\put(55.0,2.4){\line(0,1){7.5}}
\put(54.5,2.5){\line(0,1){7.5}}
\put(54.0,2.6){\line(0,1){7.4}}
\put(53.5,2.7){\line(0,1){7.4}}
\put(53.0,2.8){\line(0,1){7.2}}
\put(52.5,2.9){\line(0,1){6.9}}
\put(52.0,3.0){\line(0,1){6.6}}
\put(51.5,3.2){\line(0,1){6.3}}
\put(51.0,3.2){\line(0,1){6.1}}
\put(50.5,3.3){\line(0,1){5.8}}
\put(50.0,3.4){\line(0,1){5.6}}
\put(49.5,3.6){\line(0,1){5.2}}
\put(49.0,3.6){\line(0,1){5.0}}
\put(48.5,3.8){\line(0,1){4.7}}
\put(48.0,3.9){\line(0,1){4.4}}
\put(47.5,4.0){\line(0,1){4.2}}
\put(47.0,4.1){\line(0,1){3.9}}
\put(46.5,4.2){\line(0,1){3.7}}
\put(46.0,4.4){\line(0,1){3.3}}
\put(45.5,4.5){\line(0,1){3.1}}
\put(45.0,4.6){\line(0,1){2.8}}
\put(44.5,4.8){\line(0,1){2.4}}
\put(44.0,4.9){\line(0,1){2.2}}
\put(43.5,5.0){\line(0,1){1.9}}
\put(43.0,5.1){\line(0,1){1.7}}
\put(42.5,5.3){\line(0,1){1.3}}
\put(42.0,5.4){\line(0,1){1.1}}
\put(41.5,5.5){\line(0,1){0.8}}
\put(41.0,5.6){\line(0,1){0.6}}
\put(40.5,5.7){\line(0,1){0.4}}

\put(0.5,19.7){\line(0,1){0.3}}
\put(1.0,19.6){\line(0,1){0.4}}
\put(1.5,19.4){\line(0,1){0.6}}
\put(2.0,19.2){\line(0,1){0.8}}
\put(2.5,19.0){\line(0,1){1.0}}
\put(3.0,18.8){\line(0,1){1.1}}
\put(3.5,18.6){\line(0,1){1.2}}
\put(4.0,18.4){\line(0,1){1.4}}
\put(4.5,18.2){\line(0,1){1.5}}
\put(5.0,18.0){\line(0,1){1.7}}
\put(5.5,17.8){\line(0,1){1.8}}
\put(6.0,17.6){\line(0,1){1.9}}
\put(6.5,17.4){\line(0,1){2.1}}
\put(7.0,17.2){\line(0,1){2.3}}
\put(7.5,17.0){\line(0,1){2.4}}
\put(8.0,16.9){\line(0,1){2.5}}
\put(8.5,16.7){\line(0,1){2.6}}
\put(9.0,16.5){\line(0,1){2.7}}
\put(9.5,16.3){\line(0,1){2.8}}
\put(10.0,16.1){\line(0,1){2.9}}
\put(10.5,15.9){\line(0,1){3.1}}
\put(11.0,15.7){\line(0,1){3.1}}
\put(11.5,15.5){\line(0,1){3.3}}
\put(12.0,15.3){\line(0,1){3.4}}
\put(12.5,15.1){\line(0,1){3.5}}
\put(13.0,15.0){\line(0,1){3.6}}
\put(13.5,14.8){\line(0,1){3.7}}
\put(14.0,14.6){\line(0,1){3.8}}
\put(14.5,14.4){\line(0,1){3.9}}
\put(15.0,14.2){\line(0,1){4.0}}
\put(15.5,14.0){\line(0,1){4.1}}
\put(16.0,13.8){\line(0,1){4.2}}
\put(16.5,13.7){\line(0,1){4.3}}
\put(17.0,13.5){\line(0,1){4.4}}
\put(17.5,13.3){\line(0,1){4.5}}
\put(18.0,13.1){\line(0,1){4.6}}
\put(18.5,12.9){\line(0,1){4.6}}
\put(19.0,12.8){\line(0,1){4.6}}
\put(19.5,12.6){\line(0,1){4.7}}
\put(20.0,12.4){\line(0,1){4.8}}
\put(20.5,12.2){\line(0,1){4.9}}
\put(21.0,12.1){\line(0,1){4.9}}
\put(21.5,11.9){\line(0,1){5.0}}
\put(22.0,11.7){\line(0,1){5.1}}
\put(22.5,11.5){\line(0,1){5.2}}
\put(23.0,11.4){\line(0,1){5.2}}
\put(23.5,11.2){\line(0,1){5.3}}
\put(24.0,11.0){\line(0,1){5.4}}
\put(24.5,10.8){\line(0,1){5.5}}
\put(25.0,10.7){\line(0,1){5.5}}
\put(25.5,10.5){\line(0,1){5.6}}
\put(26.0,10.3){\line(0,1){5.7}}
\put(26.5,10.2){\line(0,1){5.7}}
\put(27.0,10.3){\line(0,1){5.5}}
\put(27.5,10.4){\line(0,1){5.3}}
\put(28.0,10.4){\line(0,1){5.2}}
\put(28.5,10.5){\line(0,1){4.9}}
\put(29.0,10.6){\line(0,1){4.7}}
\put(29.5,10.7){\line(0,1){4.5}}
\put(30.0,10.8){\line(0,1){4.3}}
\put(30.5,10.9){\line(0,1){4.1}}
\put(31.0,10.9){\line(0,1){3.9}}
\put(31.5,11.0){\line(0,1){3.7}}
\put(32.0,11.1){\line(0,1){3.5}}
\put(32.5,11.2){\line(0,1){3.3}}
\put(33.0,11.3){\line(0,1){3.1}}
\put(33.5,11.4){\line(0,1){2.9}}
\put(34.0,11.5){\line(0,1){2.7}}
\put(34.5,11.6){\line(0,1){2.4}}
\put(35.0,11.7){\line(0,1){2.2}}
\put(35.5,11.8){\line(0,1){2.0}}
\put(36.0,11.9){\line(0,1){1.8}}
\put(36.5,12.0){\line(0,1){1.6}}
\put(37.0,12.1){\line(0,1){1.3}}
\put(37.5,12.2){\line(0,1){1.1}}
\put(38.0,12.3){\line(0,1){0.9}}
\put(38.5,12.4){\line(0,1){0.7}}
\put(39.0,12.5){\line(0,1){0.5}}
\put(39.5,12.6){\line(0,1){0.3}}

\put(79.5,19.7){\line(0,1){0.3}}
\put(79.0,19.6){\line(0,1){0.4}}
\put(78.5,19.4){\line(0,1){0.6}}
\put(78.0,19.2){\line(0,1){0.8}}
\put(77.5,19.0){\line(0,1){1.0}}
\put(77.0,18.8){\line(0,1){1.1}}
\put(76.5,18.6){\line(0,1){1.2}}
\put(76.0,18.4){\line(0,1){1.4}}
\put(75.5,18.2){\line(0,1){1.5}}
\put(75.0,18.0){\line(0,1){1.7}}
\put(74.5,17.8){\line(0,1){1.8}}
\put(74.0,17.6){\line(0,1){1.9}}
\put(73.5,17.4){\line(0,1){2.1}}
\put(73.0,17.2){\line(0,1){2.3}}
\put(72.5,17.0){\line(0,1){2.4}}
\put(72.0,16.9){\line(0,1){2.5}}
\put(71.5,16.7){\line(0,1){2.6}}
\put(71.0,16.5){\line(0,1){2.7}}
\put(70.5,16.3){\line(0,1){2.8}}
\put(70.0,16.1){\line(0,1){2.9}}
\put(69.5,15.9){\line(0,1){3.1}}
\put(69.0,15.7){\line(0,1){3.1}}
\put(68.5,15.5){\line(0,1){3.3}}
\put(68.0,15.3){\line(0,1){3.4}}
\put(67.5,15.1){\line(0,1){3.5}}
\put(67.0,15.0){\line(0,1){3.6}}
\put(66.5,14.8){\line(0,1){3.7}}
\put(66.0,14.6){\line(0,1){3.8}}
\put(65.5,14.4){\line(0,1){3.9}}
\put(65.0,14.2){\line(0,1){4.0}}
\put(64.5,14.0){\line(0,1){4.1}}
\put(64.0,13.8){\line(0,1){4.2}}
\put(63.5,13.7){\line(0,1){4.3}}
\put(63.0,13.5){\line(0,1){4.4}}
\put(62.5,13.3){\line(0,1){4.5}}
\put(62.0,13.1){\line(0,1){4.6}}
\put(61.5,12.9){\line(0,1){4.6}}
\put(61.0,12.8){\line(0,1){4.6}}
\put(60.5,12.6){\line(0,1){4.7}}
\put(60.0,12.4){\line(0,1){4.8}}
\put(59.5,12.2){\line(0,1){4.9}}
\put(59.0,12.1){\line(0,1){4.9}}
\put(58.5,11.9){\line(0,1){5.0}}
\put(58.0,11.7){\line(0,1){5.1}}
\put(57.5,11.5){\line(0,1){5.2}}
\put(57.0,11.4){\line(0,1){5.2}}
\put(56.5,11.2){\line(0,1){5.3}}
\put(56.0,11.0){\line(0,1){5.4}}
\put(55.5,10.8){\line(0,1){5.5}}
\put(55.0,10.7){\line(0,1){5.5}}
\put(54.5,10.5){\line(0,1){5.6}}
\put(54.0,10.3){\line(0,1){5.7}}
\put(53.5,10.2){\line(0,1){5.7}}
\put(53.0,10.3){\line(0,1){5.5}}
\put(52.5,10.4){\line(0,1){5.3}}
\put(52.0,10.4){\line(0,1){5.2}}
\put(51.5,10.5){\line(0,1){4.9}}
\put(51.0,10.6){\line(0,1){4.7}}
\put(50.5,10.7){\line(0,1){4.5}}
\put(50.0,10.8){\line(0,1){4.3}}
\put(49.5,10.9){\line(0,1){4.1}}
\put(49.0,10.9){\line(0,1){3.9}}
\put(48.5,11.0){\line(0,1){3.7}}
\put(48.0,11.1){\line(0,1){3.5}}
\put(47.5,11.2){\line(0,1){3.3}}
\put(47.0,11.3){\line(0,1){3.1}}
\put(46.5,11.4){\line(0,1){2.9}}
\put(46.0,11.5){\line(0,1){2.7}}
\put(45.5,11.6){\line(0,1){2.4}}
\put(45.0,11.7){\line(0,1){2.2}}
\put(44.5,11.8){\line(0,1){2.0}}
\put(44.0,11.9){\line(0,1){1.8}}
\put(43.5,12.0){\line(0,1){1.6}}
\put(43.0,12.1){\line(0,1){1.3}}
\put(42.5,12.2){\line(0,1){1.1}}
\put(42.0,12.3){\line(0,1){0.9}}
\put(41.5,12.4){\line(0,1){0.7}}
\put(41.0,12.5){\line(0,1){0.5}}
\put(40.5,12.6){\line(0,1){0.3}}

\put(0.3,20.2){\line(0,1){19.8}}
\put(0.5,20.3){\line(0,1){19.7}}
\put(1.0,20.4){\line(0,1){19.5}}
\put(1.5,20.6){\line(0,1){19.4}}
\put(2.0,20.9){\line(0,1){19.1}}
\put(2.5,21.1){\line(0,1){18.8}}
\put(3.0,21.3){\line(0,1){18.6}}
\put(3.5,21.5){\line(0,1){18.4}}
\put(4.0,21.7){\line(0,1){18.2}}
\put(4.5,21.9){\line(0,1){18.0}}
\put(5.0,22.1){\line(0,1){17.8}}
\put(5.5,22.3){\line(0,1){17.6}}
\put(6.0,22.5){\line(0,1){17.4}}
\put(6.5,22.7){\line(0,1){17.2}}
\put(7.0,22.9){\line(0,1){17.0}}
\put(7.5,23.1){\line(0,1){16.8}}
\put(8.0,23.3){\line(0,1){16.5}}
\put(8.5,23.6){\line(0,1){16.2}}
\put(9.0,23.7){\line(0,1){15.9}}
\put(9.5,23.9){\line(0,1){15.7}}
\put(10.0,24.1){\line(0,1){15.5}}
\put(10.5,24.3){\line(0,1){15.3}}
\put(11.0,24.5){\line(0,1){15.1}}
\put(11.5,24.7){\line(0,1){14.9}}
\put(12.0,24.9){\line(0,1){14.6}}
\put(12.5,25.1){\line(0,1){14.3}}
\put(13.0,25.3){\line(0,1){14.0}}
\put(13.5,25.5){\line(0,1){13.7}}
\put(14.0,25.7){\line(0,1){13.5}}
\put(14.5,25.9){\line(0,1){13.3}}
\put(15.0,26.1){\line(0,1){13.1}}
\put(15.5,26.3){\line(0,1){12.8}}
\put(16.0,26.4){\line(0,1){12.5}}
\put(16.5,26.7){\line(0,1){12.2}}
\put(17.0,26.8){\line(0,1){12.0}}
\put(17.5,27.0){\line(0,1){11.7}}
\put(18.0,27.2){\line(0,1){11.5}}
\put(18.5,27.4){\line(0,1){11.3}}
\put(19.0,27.5){\line(0,1){11.1}}
\put(19.5,27.7){\line(0,1){10.8}}
\put(20.0,27.9){\line(0,1){10.5}}
\put(20.5,28.1){\line(0,1){10.2}}
\put(21.0,28.3){\line(0,1){9.9}}
\put(21.5,28.4){\line(0,1){9.7}}
\put(22.0,28.6){\line(0,1){9.5}}
\put(22.5,28.8){\line(0,1){9.2}}
\put(23.0,28.9){\line(0,1){9.0}}
\put(23.5,29.1){\line(0,1){8.8}}
\put(24.0,29.3){\line(0,1){8.5}}
\put(24.5,29.5){\line(0,1){8.2}}
\put(25.0,29.7){\line(0,1){7.9}}
\put(25.5,29.8){\line(0,1){7.6}}
\put(26.0,30.0){\line(0,1){7.4}}
\put(26.5,30.2){\line(0,1){7.0}}
\put(27.0,30.3){\line(0,1){6.9}}
\put(27.5,30.5){\line(0,1){6.6}}
\put(28.0,30.6){\line(0,1){6.4}}
\put(28.5,30.8){\line(0,1){6.1}}
\put(29.0,31.0){\line(0,1){5.8}}
\put(29.5,31.1){\line(0,1){5.6}}
\put(30.0,31.3){\line(0,1){5.3}}
\put(30.5,31.4){\line(0,1){5.1}}
\put(31.0,31.6){\line(0,1){4.8}}
\put(31.5,31.8){\line(0,1){4.5}}
\put(32.0,31.9){\line(0,1){4.3}}
\put(32.5,32.0){\line(0,1){4.0}}
\put(33.0,32.2){\line(0,1){3.7}}
\put(33.5,32.4){\line(0,1){3.4}}
\put(34.0,32.5){\line(0,1){3.2}}
\put(34.5,32.6){\line(0,1){2.9}}
\put(35.0,32.8){\line(0,1){2.6}}
\put(35.5,32.9){\line(0,1){2.4}}
\put(36.0,33.1){\line(0,1){2.1}}
\put(36.5,33.2){\line(0,1){1.9}}
\put(37.0,33.4){\line(0,1){1.5}}
\put(37.5,33.5){\line(0,1){1.3}}
\put(38.0,33.6){\line(0,1){1.1}}
\put(38.5,33.7){\line(0,1){0.9}}
\put(39.0,33.9){\line(0,1){0.5}}
\put(39.5,34.0){\line(0,1){0.3}}

\put(79.7,20.2){\line(0,1){19.8}}
\put(79.5,20.3){\line(0,1){19.7}}
\put(79.0,20.4){\line(0,1){19.5}}
\put(78.5,20.6){\line(0,1){19.4}}
\put(78.0,20.9){\line(0,1){19.1}}
\put(77.5,21.1){\line(0,1){18.8}}
\put(77.0,21.3){\line(0,1){18.6}}
\put(76.5,21.5){\line(0,1){18.4}}
\put(76.0,21.7){\line(0,1){18.2}}
\put(75.5,21.9){\line(0,1){18.0}}
\put(75.0,22.1){\line(0,1){17.8}}
\put(74.5,22.3){\line(0,1){17.6}}
\put(74.0,22.5){\line(0,1){17.4}}
\put(73.5,22.7){\line(0,1){17.2}}
\put(73.0,22.9){\line(0,1){17.0}}
\put(72.5,23.1){\line(0,1){16.8}}
\put(72.0,23.3){\line(0,1){16.5}}
\put(71.5,23.6){\line(0,1){16.2}}
\put(71.0,23.7){\line(0,1){15.9}}
\put(70.5,23.9){\line(0,1){15.7}}
\put(70.0,24.1){\line(0,1){15.5}}
\put(69.5,24.3){\line(0,1){15.3}}
\put(69.0,24.5){\line(0,1){15.1}}
\put(68.5,24.7){\line(0,1){14.9}}
\put(68.0,24.9){\line(0,1){14.6}}
\put(67.5,25.1){\line(0,1){14.3}}
\put(67.0,25.3){\line(0,1){14.0}}
\put(66.5,25.5){\line(0,1){13.7}}
\put(66.0,25.7){\line(0,1){13.5}}
\put(65.5,25.9){\line(0,1){13.3}}
\put(65.0,26.1){\line(0,1){13.1}}
\put(64.5,26.3){\line(0,1){12.8}}
\put(64.0,26.4){\line(0,1){12.5}}
\put(63.5,26.7){\line(0,1){12.2}}
\put(63.0,26.8){\line(0,1){12.0}}
\put(62.5,27.0){\line(0,1){11.7}}
\put(62.0,27.2){\line(0,1){11.5}}
\put(61.5,27.4){\line(0,1){11.3}}
\put(61.0,27.5){\line(0,1){11.1}}
\put(60.5,27.7){\line(0,1){10.8}}
\put(60.0,27.9){\line(0,1){10.5}}
\put(59.5,28.1){\line(0,1){10.2}}
\put(59.0,28.3){\line(0,1){9.9}}
\put(58.5,28.4){\line(0,1){9.7}}
\put(58.0,28.6){\line(0,1){9.5}}
\put(57.5,28.8){\line(0,1){9.2}}
\put(57.0,28.9){\line(0,1){9.0}}
\put(56.5,29.1){\line(0,1){8.8}}
\put(56.0,29.3){\line(0,1){8.5}}
\put(55.5,29.5){\line(0,1){8.2}}
\put(55.0,29.7){\line(0,1){7.9}}
\put(54.5,29.8){\line(0,1){7.6}}
\put(54.0,30.0){\line(0,1){7.4}}
\put(53.5,30.2){\line(0,1){7.0}}
\put(53.0,30.3){\line(0,1){6.9}}
\put(52.5,30.5){\line(0,1){6.6}}
\put(52.0,30.6){\line(0,1){6.4}}
\put(51.5,30.8){\line(0,1){6.1}}
\put(51.0,31.0){\line(0,1){5.8}}
\put(50.5,31.1){\line(0,1){5.6}}
\put(50.0,31.3){\line(0,1){5.3}}
\put(49.5,31.4){\line(0,1){5.1}}
\put(49.0,31.6){\line(0,1){4.8}}
\put(48.5,31.8){\line(0,1){4.5}}
\put(48.0,31.9){\line(0,1){4.3}}
\put(47.5,32.0){\line(0,1){4.0}}
\put(47.0,32.2){\line(0,1){3.7}}
\put(46.5,32.4){\line(0,1){3.4}}
\put(46.0,32.5){\line(0,1){3.2}}
\put(45.5,32.6){\line(0,1){2.9}}
\put(45.0,32.8){\line(0,1){2.6}}
\put(44.5,32.9){\line(0,1){2.4}}
\put(44.0,33.1){\line(0,1){2.1}}
\put(43.5,33.2){\line(0,1){1.9}}
\put(43.0,33.4){\line(0,1){1.5}}
\put(42.5,33.5){\line(0,1){1.3}}
\put(42.0,33.6){\line(0,1){1.1}}
\put(41.5,33.7){\line(0,1){0.9}}
\put(41.0,33.9){\line(0,1){0.5}}
\put(40.5,34.0){\line(0,1){0.3}}

\put(0.3,40.1){\line(0,1){12.2}}
\put(0.5,40.1){\line(0,1){12.2}}
\put(1.0,40.1){\line(0,1){12.2}}
\put(1.5,40.1){\line(0,1){12.2}}
\put(2.0,40.2){\line(0,1){12.2}}
\put(2.5,40.2){\line(0,1){12.2}}
\put(3.0,40.2){\line(0,1){12.1}}
\put(3.5,40.3){\line(0,1){12.0}}
\put(4.0,40.3){\line(0,1){12.0}}
\put(4.5,40.4){\line(0,1){11.9}}
\put(5.0,40.5){\line(0,1){11.8}}
\put(5.5,40.5){\line(0,1){11.8}}
\put(6.0,40.6){\line(0,1){11.6}}
\put(6.5,40.6){\line(0,1){11.6}}
\put(7.0,40.7){\line(0,1){11.5}}
\put(7.5,40.7){\line(0,1){11.5}}
\put(8.0,40.8){\line(0,1){11.4}}
\put(8.5,40.8){\line(0,1){11.3}}
\put(9.0,40.9){\line(0,1){11.2}}
\put(9.5,41.0){\line(0,1){11.1}}
\put(10.0,41.1){\line(0,1){11.0}}
\put(10.5,41.2){\line(0,1){10.9}}
\put(11.0,41.3){\line(0,1){10.8}}
\put(11.5,41.4){\line(0,1){10.6}}
\put(12.0,41.5){\line(0,1){10.5}}
\put(12.5,41.6){\line(0,1){10.4}}
\put(13.0,41.7){\line(0,1){10.2}}
\put(13.5,41.8){\line(0,1){10.0}}
\put(14.0,41.8){\line(0,1){10.0}}
\put(14.5,41.9){\line(0,1){9.7}}
\put(15.0,42.0){\line(0,1){9.6}}
\put(15.5,42.1){\line(0,1){9.5}}
\put(16.0,42.2){\line(0,1){9.4}}
\put(16.5,42.3){\line(0,1){9.3}}
\put(17.0,42.4){\line(0,1){9.1}}
\put(17.5,42.5){\line(0,1){8.9}}
\put(18.0,42.6){\line(0,1){8.8}}
\put(18.5,42.7){\line(0,1){8.6}}
\put(19.0,42.8){\line(0,1){8.4}}
\put(19.5,42.9){\line(0,1){8.3}}
\put(20.0,43.0){\line(0,1){8.2}}
\put(20.5,43.1){\line(0,1){8.0}}
\put(21.0,43.2){\line(0,1){7.8}}
\put(21.5,43.3){\line(0,1){7.6}}
\put(22.0,43.4){\line(0,1){7.4}}
\put(22.5,43.5){\line(0,1){7.2}}
\put(23.0,43.6){\line(0,1){7.1}}
\put(23.5,43.7){\line(0,1){6.9}}
\put(24.0,43.8){\line(0,1){6.7}}
\put(24.5,43.9){\line(0,1){6.5}}
\put(25.0,44.0){\line(0,1){6.3}}
\put(25.5,44.1){\line(0,1){6.1}}
\put(26.0,44.2){\line(0,1){5.9}}
\put(26.5,44.3){\line(0,1){5.8}}
\put(27.0,44.4){\line(0,1){5.6}}
\put(27.5,44.5){\line(0,1){5.5}}
\put(28.0,44.6){\line(0,1){5.3}}
\put(28.5,44.7){\line(0,1){5.1}}
\put(29.0,44.8){\line(0,1){4.9}}
\put(29.5,45.0){\line(0,1){4.6}}
\put(30.0,45.1){\line(0,1){4.4}}
\put(30.5,45.2){\line(0,1){4.2}}
\put(31.0,45.3){\line(0,1){4.0}}
\put(31.5,45.4){\line(0,1){3.8}}
\put(32.0,45.5){\line(0,1){3.6}}
\put(32.5,45.6){\line(0,1){3.4}}
\put(33.0,45.7){\line(0,1){3.2}}
\put(33.5,45.8){\line(0,1){3.0}}
\put(34.0,45.9){\line(0,1){2.8}}
\put(34.5,46.1){\line(0,1){2.5}}
\put(35.0,46.2){\line(0,1){2.3}}
\put(35.5,46.3){\line(0,1){2.1}}
\put(36.0,46.4){\line(0,1){1.9}}
\put(36.5,46.5){\line(0,1){1.6}}
\put(37.0,46.6){\line(0,1){1.4}}
\put(37.5,46.7){\line(0,1){1.2}}
\put(38.0,46.8){\line(0,1){1.0}}
\put(38.5,46.9){\line(0,1){0.8}}
\put(39.0,47.1){\line(0,1){0.5}}
\put(39.5,47.2){\line(0,1){0.3}}

\put(79.7,40.1){\line(0,1){12.2}}
\put(79.5,40.1){\line(0,1){12.2}}
\put(79.0,40.1){\line(0,1){12.2}}
\put(78.5,40.1){\line(0,1){12.2}}
\put(78.0,40.2){\line(0,1){12.2}}
\put(77.5,40.2){\line(0,1){12.2}}
\put(77.0,40.2){\line(0,1){12.1}}
\put(76.5,40.3){\line(0,1){12.0}}
\put(76.0,40.3){\line(0,1){12.0}}
\put(75.5,40.4){\line(0,1){11.9}}
\put(75.0,40.5){\line(0,1){11.8}}
\put(74.5,40.5){\line(0,1){11.8}}
\put(74.0,40.6){\line(0,1){11.6}}
\put(73.5,40.6){\line(0,1){11.6}}
\put(73.0,40.7){\line(0,1){11.5}}
\put(72.5,40.7){\line(0,1){11.5}}
\put(72.0,40.8){\line(0,1){11.4}}
\put(71.5,40.8){\line(0,1){11.3}}
\put(71.0,40.9){\line(0,1){11.2}}
\put(70.5,41.0){\line(0,1){11.1}}
\put(70.0,41.1){\line(0,1){11.0}}
\put(69.5,41.2){\line(0,1){10.9}}
\put(69.0,41.3){\line(0,1){10.8}}
\put(68.5,41.4){\line(0,1){10.6}}
\put(68.0,41.5){\line(0,1){10.5}}
\put(67.5,41.6){\line(0,1){10.4}}
\put(67.0,41.7){\line(0,1){10.2}}
\put(66.5,41.8){\line(0,1){10.0}}
\put(66.0,41.8){\line(0,1){10.0}}
\put(65.5,41.9){\line(0,1){9.7}}
\put(65.0,42.0){\line(0,1){9.6}}
\put(64.5,42.1){\line(0,1){9.5}}
\put(64.0,42.2){\line(0,1){9.4}}
\put(63.5,42.3){\line(0,1){9.3}}
\put(63.0,42.4){\line(0,1){9.1}}
\put(62.5,42.5){\line(0,1){8.9}}
\put(62.0,42.6){\line(0,1){8.8}}
\put(61.5,42.7){\line(0,1){8.6}}
\put(61.0,42.8){\line(0,1){8.4}}
\put(60.5,42.9){\line(0,1){8.3}}
\put(60.0,43.0){\line(0,1){8.2}}
\put(59.5,43.1){\line(0,1){8.0}}
\put(59.0,43.2){\line(0,1){7.8}}
\put(58.5,43.3){\line(0,1){7.6}}
\put(58.0,43.4){\line(0,1){7.4}}
\put(57.5,43.5){\line(0,1){7.2}}
\put(57.0,43.6){\line(0,1){7.1}}
\put(56.5,43.7){\line(0,1){6.9}}
\put(56.0,43.8){\line(0,1){6.7}}
\put(55.5,43.9){\line(0,1){6.5}}
\put(55.0,44.0){\line(0,1){6.3}}
\put(54.5,44.1){\line(0,1){6.1}}
\put(54.0,44.2){\line(0,1){5.9}}
\put(53.5,44.3){\line(0,1){5.8}}
\put(53.0,44.4){\line(0,1){5.6}}
\put(52.5,44.5){\line(0,1){5.5}}
\put(52.0,44.6){\line(0,1){5.3}}
\put(51.5,44.7){\line(0,1){5.1}}
\put(51.0,44.8){\line(0,1){4.9}}
\put(50.5,45.0){\line(0,1){4.6}}
\put(50.0,45.1){\line(0,1){4.4}}
\put(49.5,45.2){\line(0,1){4.2}}
\put(49.0,45.3){\line(0,1){4.0}}
\put(48.5,45.4){\line(0,1){3.8}}
\put(48.0,45.5){\line(0,1){3.6}}
\put(47.5,45.6){\line(0,1){3.4}}
\put(47.0,45.7){\line(0,1){3.2}}
\put(46.5,45.8){\line(0,1){3.0}}
\put(46.0,45.9){\line(0,1){2.8}}
\put(45.5,46.1){\line(0,1){2.5}}
\put(45.0,46.2){\line(0,1){2.3}}
\put(44.5,46.3){\line(0,1){2.1}}
\put(44.0,46.4){\line(0,1){1.9}}
\put(43.5,46.5){\line(0,1){1.6}}
\put(43.0,46.6){\line(0,1){1.4}}
\put(42.5,46.7){\line(0,1){1.2}}
\put(42.0,46.8){\line(0,1){1.0}}
\put(41.5,46.9){\line(0,1){0.8}}
\put(41.0,47.1){\line(0,1){0.5}}
\put(40.5,47.2){\line(0,1){0.3}}

\end{picture}
\vspace{8mm}
\caption{\footnotesize  Dependence of the spectrum $\s(\D_\a)$ on the magnetic flux $\phi:=\alpha(\mathbf{c}_0)$ through the cycle $\mathbf{c}_0=(\mathbf{e}_2,\mathbf{e}_3,\mathbf{e}_4,\mathbf{e}_5)$ of the periodic graph~$\mathcal{G}$.} \label{FSD}
\end{figure}

\begin{remark}\lb{rExa}
\emph{i}) The influence of the magnetic field on the spectrum of the magnetic Laplacian $\D_\a$ on $\cG$ is completely determined by the flux $\phi:=\a(\bc_0)$ through the single (up to periodicity) cycle $\bc_0=(\be_2,\be_3,\be_4,\be_5)$ of the periodic graph $\cG$, see Fig.~\ref{FEx1}\emph{a} and also Remark \ref{Re1}.\emph{i}. Dependence of the spectrum of $\D_\a$ on the magnetic flux $\phi\in[0,2\pi]$ is shown in Fig.~\ref{FSD}. The magnetic flux $\phi$ is plotted on the horizontal axis. The vertical axis represents the spectrum of $\D_\a$. The spectral bands are the black vertical intervals which appear as the intersection of the black region with a line $\phi=\const$ and the spectral gaps are the white vertical intervals.

\emph{ii}) In \er{esti}, the first lower bound
$\frac89\,\cos^2\frac{\phi}2\leq\gS(\D_\a)$ is obtained by using the estimate \er{leHa}, and for the second one $\frac{4}{9}\,
\big|\cos\frac{\phi}2\big|\leq\gS(\D_\a)$ we use the estimate \er{EsFS}.
If $\phi\in[0,\frac{2\pi}3)$ or $\phi\in(\frac{4\pi}3,2\pi]$, then the first lower bound is better than the second one. For $\phi\in(\frac{2\pi}3,\frac{4\pi}3)$, the better estimate is the second one. If $\phi=\frac{2\pi}3$ or $\phi=\frac{4\pi}3$, then the both estimates are the same.

\emph{iii}) If the magnetic flux $\phi=0$, then the spectrum of $\D_\a$ is given by
$$
\begin{aligned}
&\s(\D_\a)=\s_{ac}(\D_\a)\cup\s_{fb}(\D_\a), \qqq \textrm{where}\\
&\s_{ac}(\D_\a)=\s_1\cup\s_3\cup\s_4,\qqq \s_{fb}(\D_\a)=\s_2=\{2\},\\
&\s_1=[0,3-\sqrt{5}\;],\qqq \s_3=[2,4],\qqq \s_4=[4,3+\sqrt{5}\;],
\end{aligned}
$$
see Fig.~\ref{FSD}, and the total bandwidth $\gS(\D_\a)=4$.
\end{remark}

It is known (see Theorem 2.1.\emph{ii} in \cite{KS19}) that the first spectral band of the Schr\"odinger operator $H_0$ without magnetic fields is non-degenerate. Thus, the spectrum of $H_0$ always has an absolutely continuous component. Example \ref{Exa}.\emph{ii} shows that, in contrast to the non-magnetic case, all spectral bands of the magnetic Schr\"odinger  operator $H_\a$ may be flat, i.e., the absolutely continuous spectrum of $H_\a$ may be empty. In the following proposition we prove that this is a quite rare situation.

\begin{proposition}\lb{Tfmp}
Let $H_{t\a}=-\D_{t\a}+V$ be the Schr\"odinger operator defined  by \er{Sh} -- \er{ALO} with a periodic magnetic potential $t\a$ and a periodic electric potential $V$ on a periodic graph $\cG$, where $t$ is a coupling constant. Then the absolutely continuous spectrum of $H_{t\a}$ is not empty for all except finitely many $t\in[0,1]$.
\end{proposition}

\begin{remark}
\emph{i}) In order to prove Proposition \ref{Tfmp} we use the result from \cite{KS19} that the first spectral band of the Schr\"odinger operator $H_0$ without magnetic field is not flat and classical function theory.

\emph{ii}) Proposition \ref{Tfmp} also holds with a similar proof for weighted magnetic Laplacians and Schr\"odinger operators, in particular, for the normalized ones.
\end{remark}

\section{Trace formulas for magnetic Schr\"odinger operators}
\setcounter{equation}{0}
\lb{Sec3}

\subsection{Trace formulas for magnetic adjacency operators} In this section we discuss trace formulas for the magnetic adjacency operator $A_\a$ defined by \er{ALO}. In the standard orthonormal basis of $\ell^2(\cV_*)=\C^\n$, $\n=\#\cV_*$, the $\n\ts\n$ matrix $A_\a(k)=\big(A_{\a,xy}(k)\big)_{x,y\in\cV_*}$ of the fiber magnetic adjacency
operator $A_\a(k)$ given by \er{fado} has the form
\[\lb{fnml0}
A_{\a,xy}(k)=\sum_{\be=(x,y)\in\cA_*}e^{-i(\a(\be)+\lan\t(\be),k\ran)},\qqq k\in\T^d.
\]
The diagonal entries of the matrix $A_\a(k)$ can also be written as
$$
A_{\a,xx}(k)=\sum\limits_{\be=(x,x)\in\cA_*}
\cos\big(\a(\be)+\lan\t(\be),k\ran\big), \qqq \forall\,x\in\cV_*.
$$

The eigenvalues of the fiber magnetic adjacency operator $A_\a(k)$ will be denoted  by $\l^o_{\a,j}(k)$, $j\in\N_\n$. The spectral bands $\s_j(A_\a)$, $j\in\N_{\n}$, for the magnetic adjacency operator $A_\a$ have the form
$$
\begin{array}{l}
\s_j(A_\a)=[\l_{\a,j}^{o-},\l_{\a,j}^{o+}]=\l_{\a,j}^o(\T^d),\\[8pt]
\displaystyle\textrm{where}\qqq \l_{\a,j}^{o-}=\min_{k\in\T^d}\l^o_{\a,j}(k),\qqq \l_{\a,j}^{o+}=\max_{k\in\T^d}\l^o_{\a,j}(k).
\end{array}
$$

We express the traces for the $n$-th power of the fiber magnetic adjacency operator in terms of magnetic fluxes through cycles of the fundamental graph from the sets $\cC_n$ and $\cC_n^\mm$ defined by \er{cNn} and \er{cNnm}. First, we recall some simple properties of the numbers $\cN_n=\#\cC_n$ and $\cN_n^\mm=\#\cC_n^\mm$ (see Proposition 2.1.\emph{i} in \cite{KS22}).

\begin{proposition}\lb{spcN}
The numbers $\cN_n=\#\cC_n$ and $\cN_n^\mm=\#\cC_n^\mm$ satisfy
$$
\cN_n\leq\n\vk^n_+,\qqq \forall\,n\in\N, \qqq \n=\#\cV_*,
$$
$$
\cN_{2n}^0\geq\#\cA_*,\qqq\cN_n^\mm=\cN_n^{-\mm},\qqq \forall\,(n,\mm)\in\N\ts\Z^d,
$$
\[\lb{cNe0}
\cN_n^\mm=0, \qqq \textrm{if}\qqq \|\mm\|>n\t_+,\qqq\textrm{where}\qqq \t_+=\max\limits_{\be\in\cA_*}\|\t(\be)\|,
\]
$\t(\be)$ is the index of the edge $\be\in\cA_*$  defined by \er{in}, \er{dco}; $\|\cdot\|$ is the standard norm in $\R^d$, and $\vk_+$ is given in \er{spAD}.

If there are no loops and multiple edges in the periodic graph $\cG$, then
$\cN_2^0=\#\cA_*$.
\end{proposition}

Now we describe trace formulas for the magnetic adjacency operator $A_\a$.

\begin{theorem}
\lb{TFao} Let $A_\a(k)$, $k\in\T^d$, be the fiber magnetic adjacency operator defined by \er{fado} on the fundamental graph $\cG_*=(\cV_*,\cA_*)$. Then for each $n\in\N$

i) The trace of $A_\a^n(k)$ has the form
\[\lb{TrAo}
\begin{aligned}
&\Tr A_\a^n(k)=\sum_{j=1}^\n\big(\l_{\a,j}^o(k)\big)^n=
\sum_{\bc\in\cC_n}e^{-i\f(\bc,k)}=
\sum_{\bc\in\cC_n}\cos\f(\bc,k), \\
&\f(\bc,k)=\a(\bc)+\lan\t(\bc),k\ran,\qqq \n=\#\cV_*,
\end{aligned}
\]
or, in the form of Fourier series,
\[\lb{FTank}
\Tr A_\a^n(k)=\sum_{\mm\in\Z^d\atop \|\mm\|\leq n\t_+}
T_{\a,n,\mm}\,e^{-i\lan\mm,k\ran},\qqq
T_{\a,n,\mm}=\sum_{\bc\in\cC_n^\mm}e^{-i\a(\bc)},
\]
where $\t_+$ is defined in \er{tapl}.

ii) The trace of $A_\a^n(k)$ satisfies
\[\lb{ITrA}
\frac1{(2\pi)^d}\int_{\T^d}\Tr A_\a^n(k)dk=T_{\a,n,0}, \qqq
T_{\a,n,0}=\sum_{\bc\in\cC_n^0}\cos \a(\bc).
\]
\end{theorem}

\begin{remark}
\emph{i}) The traces $\Tr A_\a^n(k)$ are finite Fourier series of the quasimomentum $k\in\T^d$.

\emph{ii}) For $k=0$ the trace formulas \er{TrAo} have the form
$$
\Tr A_\a^n(0)=\sum_{\bc\in\cC_n}\cos\a(\bc),
$$
and coincide with the trace formulas for the magnetic adjacency operator $A_\a(0)$ on a finite graph $\cG_*$ presented in \cite{OGS09}. The traces for the $n$-th power of $A_\a(0)$ are expressed in terms of the magnetic fluxes $\a(\bc)$ through cycles $\bc$ of length $n$ in the graph $\cG_*$.
\end{remark}

We need the following simple identity
\[\lb{inex}
\frac1{(2\pi)^d}\int_{\T^d}e^{i\lan\mm,k\ran}dk=
\frac1{(2\pi)^d}\int_{\T^d}\cos\lan\mm,k\ran\,dk=\d_{\mm,0},\qqq \forall\,\mm\in\Z^d,
\]
where $\d_{\mm,0}$ is the Kronecker delta.

\medskip

\no\textbf{Proof of Theorem \ref{TFao}.} \emph{i}) Using \er{fnml0} and applying the definitions \er{cyin} and \er{cyin2} of the cycle index $\t(\bc)$ and the magnetic flux $\a(\bc)$, we obtain that for each $n\in\N$
\begin{multline*}
\Tr A_\a^n(k)=\sum_{x_1,\ldots,x_n\in\cV_*}A_{\a,x_1x_2}(k)A_{\a,x_2x_3}(k)
\ldots A_{\a,x_{n-1}x_n}(k)A_{\a,x_nx_1}(k)\\
=\sum_{x_1,\ldots,x_n\in\cV_*}\sum\limits_{\be_1,\ldots,\be_n\in\cA_*}
e^{-i(\a(\be_1)+\a(\be_2)+\ldots+\a(\be_n)+
\lan\t(\be_1)+\t(\be_2)+\ldots+\t(\be_n),k\ran)}
=\sum_{\bc\in\cC_n}e^{-i(\a(\bc)+\lan\t(\bc),k\ran)},
\end{multline*}
where $\be_s=(x_s,x_{s+1})$, $s\in\N_n$, $x_{n+1}=x_1$. Since for each cycle $\bc\in\cC_n$ there exists a reverse cycle $\ul\bc\in\cC_n$, and, due to \er{cyin1}, $\a(\ul\bc\,)=-\a(\bc)$ and $\t(\ul\bc\,)=-\t(\bc)$, we have \er{TrAo}.

The identity \er{TrAo} is a finite Fourier series for the $2\pi\Z^d$-periodic function $\Tr A_\a^n(k)$, since $\t(\bc)\in\Z^d$. We rewrite this Fourier series in the standard form
\[\lb{FSAa}
\Tr A_\a^n(k)=\sum_{\bc\in\cC_n}e^{-i(\a(\bc)+\lan\t(\bc),k\ran)}=
\sum_{\mm\in\Z^d}\sum_{\bc\in\cC_n^\mm}e^{-i\a(\bc)}e^{-i\lan\mm,k\ran}
=\sum\limits_{\mm\in\Z^d}T_{\a,n,\mm}e^{-i\lan\mm,k\ran},
\]
where the coefficients $T_{\a,n,\mm}$ are given in \er{FTank}. Due to \er{cNe0}, $T_{\a,n,\mm}=0$ for all $\mm\in\Z^d$ such that $\|\mm\|>n\t_+$, i.e., the sum in RHS of \er{FSAa} is finite. Thus, the identity \er{FSAa} can be written in the form \er{FTank}.

\emph{ii})  Integrating \er{FTank} over $k\in\T^d$ and using \er{inex}, we obtain the first identity in \er{ITrA}, where
$$
T_{\a,n,0}=\sum_{\bc\in\cC_n^0}e^{-i\a(\bc)}.
$$
Finally, using that for each cycle $\bc\in\cC_n^0$ there exists a reverse cycle $\ul\bc\in\cC_n^0$, and $\a(\ul\bc\,)=-\a(\bc)$, we get the second formula in \er{ITrA}. \qq $\BBox$

\subsection{Trace formulas for magnetic Schr\"odinger operators}
In the standard orthonormal basis of  $\ell^2(\cV_*)=\C^\n$, $\n=\#\cV_*$, the $\n\ts\n$ matrix $H_\a(k)=\big(H_{\a,xy}(k)\big)_{x,y\in\cV_*}$ of the fiber magnetic Schr\"odinger operator $H_\a(k)$ defined by \er{Hvt'} -- \er{fado} has the form
\[\lb{mrSo}
H_\a(k)=A_\a(k)+v, \qqq\textrm{where}\qqq v=\diag(v_x)_{x\in\cV_*}, \qqq v_x=V_x-\vk_x,
\]
and the $\n\ts\n$ matrix $A_\a(k)$ of the fiber magnetic adjacency operator is given by \er{fnml0}.

The next lemma shows that the operator $H_\a(k)$ can be considered as a fiber weighted magnetic adjacency operator on the modified fundamental graph $\wt\cG_*=(\cV_*,\wt\cA_*)$.

\begin{lemma}\lb{LDeH}
The fiber magnetic Schr\"odinger operator $H_\a(k)=\big(H_{\a,xy}(k)\big)_{x,y\in\cV_*}$ given by \er{mrSo}, \er{fnml0} satisfies
\[\lb{mrmL+}
H_{\a,xy}(k)=\sum_{\be=(x,y)\in\wt\cA_*}\o(\be) e^{-i(\a(\be)+\lan\t(\be),\,k\ran)},\qqq \forall\, x,y\in\cV_*, \qqq \forall\,k\in\T^d,
\]
where $\wt\cA_*$ is the set of all edges of the modified fundamental graph $\wt\cG_*$ defined by \er{wtAs}; $\o(\be)$ is given by \er{webe}, and $\t(\be)$  is the index of the edge $\be\in\wt\cA_*$ defined by \er{in}, \er{dco}.
\end{lemma}

\no \textbf{Proof.} Let $x,y\in\cV_*$. If $x\neq y$, then, using
\er{webe}  and \er{mrSo}, \er{fnml0}, we have
$$
\sum_{\be=(x,y)\in\wt\cA_*}\o(\be) e^{-i(\a(\be)+\lan\t(\be),\,k\ran)}
=\sum_{\be=(x,y)\in\cA_*} e^{-i(\a(\be)+\lan\t(\be),\,k\ran)}=H_{\a,xy}(k).
$$
Similarly, if $x=y$, then we obtain
\begin{multline*}
\sum_{\be=(x,x)\in\wt\cA_*}\o(\be) e^{-i(\a(\be)+\lan\t(\be),\,k\ran)}
=\o(\be_x)+\sum_{\be=(x,x)\in\cA_*} e^{-i(\a(\be)+\lan\t(\be),\,k\ran)}
\\=v_x+\sum_{\be=(x,x)\in\cA_*} e^{-i(\a(\be)+\lan\t(\be),\,k\ran)}=H_{\a,xx}(k).
\end{multline*}
Thus, the identity \er{mrmL+} has been proved. \qq \BBox

\medskip

Now we prove Theorem \ref{TPG} about trace formulas for  the fiber magnetic Schr\"odinger operator $H_\a(k)$.

\medskip

\no\textbf{Proof of Theorem \ref{TPG}.} The proof is similar to the proof  of Theorem \ref{TFao}.

\emph{i}) Using \er{mrmL+} and applying the definitions \er{cyin}, \er{cyin2} and \er{Wcy} of the cycle index $\t(\bc)$, the magnetic flux $\a(\bc)$ and the cycle weight $\o(\bc)$, we obtain that for each $n\in\N$
\begin{multline*}
\Tr H_\a^n(k)=
\sum_{x_1,\ldots,x_n\in\cV_*}H_{\a,x_1x_2}(k)H_{\a,x_2x_3}(k)\ldots H_{\a,x_{n-1}x_n}(k)H_{\a,x_nx_1}(k)
\\=\!\!\!\! \sum_{x_1,\ldots,x_n\in\cV_*}
\sum_{\be_1,\ldots,\be_n\in\wt\cA_*}\!\!\o(\be_1)\o(\be_2)\ldots\o(\be_n)
e^{-i(\a(\be_1)+\a(\be_2)+\ldots+\a(\be_n)+
\lan\t(\be_1)+\t(\be_2)+\ldots+\t(\be_n),k\ran)}\\=
\sum_{\bc\in\wt\cC_n}\o(\bc)e^{-i(\a(\bc)+\lan\t(\bc),k\ran)}=\cT_{\a,n}(k),
\end{multline*}
where $\be_s=(x_s,x_{s+1})$, $s\in \N_n$, $x_{n+1}=x_1$, and $\cT_{\a,n}(k)$ is given by the first identity in \er{deTnk}. Since for each cycle $\bc\in\wt\cC_n$ there exists a reverse cycle $\ul\bc\in\wt\cC_n$, and, due to \er{cyin1} and \er{omin}, $\a(\ul\bc\,)=-\a(\bc)$, $\t(\ul\bc\,)=-\t(\bc)$ and $\o(\ul\bc\,)=\o(\bc)$, we have \er{deTnk}.

The identity \er{deTnk} is a finite Fourier series for the $2\pi\Z^d$-periodic function $\cT_{\a,n}(k)$, since $\t(\bc)\in\Z^d$. We rewrite this Fourier series in the standard form
\[\lb{rg11H}
\cT_{\a,n}(k)=\sum_{\bc\in\wt\cC_n}\o(\bc)\,e^{-i(\a(\bc)+\lan\t(\bc),k\ran)}
=\sum_{\mm\in\Z^d}\sum_{\bc\in\wt\cC_n^\mm}\o(\bc)\,e^{-i\a(\bc)}e^{-i\lan\mm,k\ran}
=\sum_{\mm\in\Z^d}\cT_{\a,n,\mm}e^{-i\lan\mm,k\ran},
\]
where the coefficients $\cT_{\a,n,\mm}$ are given in \er{FsTrH}. Due to \er{cNe0}, $\cT_{\a,n,\mm}=0$ for all $\mm\in\Z^d$ such that $\|\mm\|>n\t_+$, i.e., the sum in RHS of \er{rg11H} is finite. Thus, the identity \er{rg11H} can be written in the form \er{FsTrH}.

\emph{ii}) Integrating \er{TrnH}, \er{FsTrH}  over $k\in\T^d$ and using \er{inex}, we obtain \er{ITrH}, where
$$
\cT_{\a,n,0}=\sum_{\bc\in\wt\cC_n^0}\o(\bc)e^{-i\a(\bc)}.
$$
Finally, using that for each cycle $\bc\in\wt\cC_n^0$ there exists a reverse cycle $\ul\bc\in\wt\cC_n^0$ and $\a(\ul\bc\,)=-\a(\bc)$, $\o(\ul\bc\,)=\o(\bc)$, we get \er{IVFu}. \qq \BBox

\medskip

In the next proposition we give explicit expressions for the first and second order traces for the magnetic Schr\"odinger operator.

\begin{proposition}\lb{TrTO}
Let $H_\a(k)=-\D_\a(k)+V$, $k\in\T^d$, be the fiber magnetic Schr\"odinger operator defined by \er{Hvt'} -- \er{fado} on the fundamental graph $\cG_*=(\cV_*,\cA_*)$. Then
\[\lb{TrH123}
\begin{aligned}
&\Tr H_\a(k)=\cT_{\a,1}(k)=\sum_{x\in\cV_*}v_x+\sum_{\bc\in\cC_1}\cos\f(\bc,k),\qqq v_x=V_x-\vk_x,\\
&\Tr H_\a^2(k)=\cT_{\a,2}(k)=\sum_{x\in\cV_*}v_x^2+
2\sum_{\bc\in\cC_1}v_{x_\bc}\cos\f(\bc,k)+
\sum_{\bc\in\cC_2}\cos\f(\bc,k),\\
& \f(\bc,k)=\a(\bc)+\lan\t(\bc),k\ran,
\end{aligned}
\]
\[\lb{TR13IH}
\begin{aligned}
&\frac1{(2\pi)^d}\int_{\T^d}\Tr H_\a(k)dk=\cT_{\a,1,0}=\sum_{x\in\cV_*}v_x+
\sum_{\bc\in\cC_1^0}\cos\a(\bc),
\\
&\frac1{(2\pi)^d}\int_{\T^d}\Tr H_\a^2(k)dk=\cT_{\a,2,0}=\sum_{x\in\cV_*}
v_x^2+2\sum_{\bc\in\cC_1^0}v_{x_\bc}\cos\a(\bc)+
\sum_{\bc\in\cC_2^0}\cos\a(\bc),
\end{aligned}
\]
where $x_\bc$ is the unique vertex of the cycle $\bc\in\cC_1$, and $\cC_n,\cC_n^0$, $n=1,2$, are defined in \er{cNn},\er{cNnm}.

If there are no loops in the fundamental graph $\cG_*$, then
\[\lb{hnl}
\Tr H_\a(k)=\frac1{(2\pi)^d}\int_{\T^d}\Tr H_\a(k)dk=\sum_{x\in\cV_*}v_x,
\]
and if, in addition, there are no multiple edges in $\cG_*$, then
\[\lb{nome}
\Tr H_\a^2(k)=\frac1{(2\pi)^d}\int_{\T^d}\Tr H_\a^2(k)dk=\#\cA_*+\|v\|_*^2,
\]
where $\|v\|_*^2=\sum\limits_{x\in\cV_*}v_x^2$ and
$\#\cA_*=\sum\limits_{x\in\cV_*}\vk_x$.
\end{proposition}

\no \textbf{Proof.} The identities \er{TrnH}, \er{deTnk} when $n=1$ have the form
$$
\Tr H_\a(k)=\cT_{\a,1}(k)=
\sum_{\bc\in\wt\cC_1}\o(\bc)\cos\f(\bc,k)=
\sum_{\bc\in\wt\cC_1\sm\cC_1}\o(\bc)\cos\f(\bc,k)+
\sum_{\bc\in\cC_1}\cos\f(\bc,k),
$$
where we have used \er{oc1C}. Since all cycles in $\wt\cC_1\sm\cC_1$ are loops $\be_x$, $x\in\cV_*$, added to the fundamental graph $\cG_*$, and
$\t(\be_x)=0$, $\a(\be_x)=0$, $\o(\be_x)=v_x$ (see \er{wtAs} and \er{webe}), we obtain the first identity in \er{TrH123}. In particular, if there are no loops in the fundamental graph $\cG_*$, i.e., $\cC_1=\varnothing$, then this identity has the form \er{hnl}. Similarly, the identities \er{TrnH}, \er{deTnk} when $n=2$  have the form
\[\label{THco2}
\Tr H_\a^2(k)=\cT_{\a,2}(k)=\sum_{\bc\in\wt\cC_2}\o(\bc)
\cos\f(\bc,k)=\sum_{\bc\in\wt\cC_2\sm\cC_2}\o(\bc)\cos\f(\bc,k)+
\sum_{\bc\in\cC_2}\cos\f(\bc,k).
\]
Each cycle from $\wt\cC_2\sm\cC_2$ has one of the following forms:

$\bu$ $\bc_1=(\be_x,\be_x)$, $x\in\cV_*$, where $\be_x$ is the added loop at the vertex $x$ with index $\t(\be_x)=0$, magnetic potential $\a(\be_x)=0$ and weight $\o(\be_x)=v_x$;

$\bu$ $\bc_2=(\be_x,\bc)$ or $\bc_3=(\bc,\be_x)$, $x\in\cV_*$, where $\bc\in\cC_1$ is a loop at the vertex $x$ (i.e., a cycle of length one) in $\cG_*$ with weight $\o(\bc)=1$, see \er{oc1C}.
\\
Due to the definitions \er{cyin}, \er{cyin2} and \er{Wcy} of the cycle index, the magnetic flux and the cycle weight, we have
$$
\begin{array}{ll}
\a(\bc_1)=0,\qqq & \a(\bc_2)=\a(\bc_3)=\a(\bc),\\[4pt]
\t(\bc_1)=0,\qqq & \t(\bc_2)=\t(\bc_3)=\t(\bc),\\[4pt]
\o(\bc_1)=v_x^2,\qq & \o(\bc_2)=\o(\bc_3)=v_x=v_{x_\bc}.
\end{array}
$$
Then the first sum in RHS of \er{THco2} has the form
\[\lb{fsu1}
\sum_{\bc\in\wt\cC_2\sm\cC_2}\o(\bc)\cos\f(\bc,k)=
\sum_{x\in\cV_*}v_x^2+2\sum_{\bc\in\cC_1}v_{x_\bc}\cos\f(\bc,k).
\]
Combining \er{THco2} and \er{fsu1}, we obtain the second identity in \er{TrH123}. Integrating \er{TrH123} over $k\in\T^d$ and using \er{inex}, we get \er{TR13IH}.

If there are no loops and multiple edges in the fundamental graph $\cG_*$, then $\cC_1=\varnothing$ and all cycles of length 2 in $\cG_*$ have the form $\bc=(\be,\ul\be\,)$, $\be\in\cA_*$, and $\a(\bc)=0$, $\t(\bc)=0$. Thus,
$$
\sum_{\bc\in\cC_2}\cos\f(\bc,k)=\#\cA_*,
$$
and the second identity in \er{TrH123} yields \er{nome}. Here we have also used the well-known identity $\sum\limits_{x\in\cV_*}\vk_x=\#\cA_*$. \qq \BBox

\medskip

The traces $\Tr H_\a^n(k)$, $(n,k)\in\N\ts\T^d$, given by \er{TrnH}, \er{deTnk}, depend on the magnetic fluxes $\a(\bc)$ through the cycles $\bc$ of the fundamental graph $\cG_*$. Therefore, the effect of the magnetic field on the traces $\Tr H_\a^n(k)$ is fully determined by $\b$ fluxes through a chosen set of basis cycles of $\cG_*$, where $\b$ is the Betti number of the fundamental graph $\cG_*=(\cV_*,\cE_*)$. We show that using a shift of the quasimomentum $k\in\T^d$, we can reduce the number of these independent fluxes to $\b-d$, where $d$ is the dimension of the periodic graph. First we need the following lemma.

\begin{lemma}
\lb{LbocC}
Let $\cG$ be a $\G$-periodic graph with the fundamental graph $\cG_*=(\cV_*,\cA_*)=\cG/\G$. Then

i) There exists a basis $\cB=\{\bc_1,\ldots,\bc_\b\}$ of the cycle space $\cC$ of $\cG_*$ such that
\[\lb{bacC}
\begin{array}{l}
\textrm{$\bu$ $\t(\bc_1),\ldots,\t(\bc_d)\in\Z^d$ form the standard orthonormal basis for the lattice $\Z^d$},\\[2pt]
\textrm{$\bu$  $\{\bc_{d+1},\ldots,\bc_\b\}$ is a basis of the subspace $\cC^0$ of all cycles of $\cG_*$ with zero index.}
\end{array}
\]
Here $d$ is the rank of the lattice $\G$; $\b$ is the Betti number of $\cG_*$, and $\t(\bc)$ is the index of the cycle $\bc$ defined by \er{cyin}.

ii) $\wt\cB=\{\bc_1,\ldots,\bc_\b\}\cup\{\be_x\}_{x\in\cV_*}$ is the basis of the cycle space $\wt\cC$ of the modified fundamental graph
\[\lb{mfg}
\wt\cG_*=(\cV_*,\wt\cA_*), \qqq \wt\cA_*=\cA_*\cup\{\be_x\}_{x\in\cV_*},
\]
where $\be_x$ is the loop at the vertex $x$.
\end{lemma}

\no \textbf{Proof.} \emph{i}) The image of the function $\t:\cC\ra\Z^d$ defined by \er{cyin} is the lattice $\Z^d$ (see, e.g., Lemma 4.3.\emph{i} in \cite{KS22a}). Then there exist cycles $\bc_1,\ldots,\bc_d\in\cC$ such that $\{\t(\bc_s)\}_{s=1}^d$ is the standard orthonormal basis for $\Z^d$. It is known (see \cite{S13}, p.81 and also Theorem 3, p.211) that $\dim\cC^0=\b-d$ and if $\{\bc_{d+1},\ldots,\bc_\b\}$ is a basis of $\cC^0$, then $\{\bc_1,\ldots,\bc_d,\bc_{d+1},\ldots,\bc_\b\}$ is a  basis of $\cC$. By construction, this basis satisfies the condition \er{bacC}.

\emph{ii}) The cycles $\bc_1,\ldots,\bc_\b$, $\be_x$, $x\in\cV_*$, comprise
a linearly independent set, because $\cB=\{\bc_1,\ldots,\bc_\b\}$ is a basis of $\cC$ and each loop edge $\be_x$ belongs to no other cycle from $\wt\cB$. Since
$$
\#\wt\cB=\b+\n=\#\cE_*-\#\cV_*+1+\n=\#\wt\cE_*-\#\cV_*+1,
$$
and this is the dimension of the cycle space $\wt\cC$ of $\wt\cG_*=(\wt\cE_*,\cV_*)$, then $\wt\cB$ is a basis of $\wt\cC$. \qq $\BBox$

\medskip

\medskip

Due to Lemma \ref{LbocC}.\emph{ii}, each cycle $\bc\in\wt\cC$ of the modified fundamental graph has the following unique representation with respect to the basis $\{\bc_1,\ldots,\bc_\b\}\cup\{\be_x\}_{x\in\cV_*}$ of $\wt\cC$:
\[\lb{cyre1}
\textstyle\bc=\sum\limits_{s=1}^\b n_s(\bc)\bc_s+\sum\limits_{x\in\cV_*}n_x(\bc)\be_x \qqq\textrm{for some} \qqq n_1(\bc),\ldots,n_\b(\bc),n_x(\bc)\in\Z.
\]
We define the projection $\cP$ of $\wt\cC$ onto $\cC^0$ by
\[\lb{depr}
\cP:\wt\cC\to\cC^0,\qqq \cP\bc=\sum\limits_{s=d+1}^\b n_s(\bc)\bc_s, \qqq\forall\,\bc\in\wt\cC,
\vspace{-4mm}
\]
where
\begin{itemize}
  \item $\wt\cC$ is the space of all cycles in the modified fundamental graph $\wt\cG_*$;
  \item $\cC^0$ is the space of all cycles with zero index in the fundamental graph $\cG_*$;
  \item $n_s(\bc)\in\Z$ are uniquely determined by \er{cyre1}.
\end{itemize}

\begin{corollary}[\textbf{Minimal magnetic fluxes representation}]
\label{TCo0} Let $H_\a=-\D_\a+V$ be the Schr\"odinger operator defined by \er{Sh} -- \er{ALO}  with a periodic magnetic potential $\a$ and a periodic electric potential $V$ on a periodic graph $\cG$. Then

i) The traces $\Tr H_\a^n(k)$, $n\in\N$, of the fiber operators $H_\a^n(k)$, $k\in\T^d$, given by \er{TrnH}, \er{deTnk}, in the shifted quasimomentum
\[\lb{sfqm}
\wt k=k-k_o\in\T^d, \qqq \textrm{where}\qqq k_o=-(\a(\bc_1),\ldots,\a(\bc_d)),
\]
have the form
\[
\lb{wtl2.13}
\begin{aligned}
&\Tr H_\a^n(\wt k+k_o)=\sum\limits_{\bc\in\wt\cC_n}\o(\bc)\,e^{-i\f_o(\bc,\wt k\,)}=
\sum\limits_{\bc\in\wt\cC_n}\o(\bc)\cos\f_o(\bc,\wt k\,), \\ &\f_o(\bc,\wt k\,)=\f(\bc,\wt k+k_o)=\a(\cP\bc)+\lan\t(\bc),\wt k\,\ran,
\end{aligned}
\]
or, in the form of Fourier series,
\[\lb{TrHM}
\Tr H_\a^n(\wt k+k_o)=\sum\limits_{\mm\in\Z^d\atop \|\mm\|\leq n\t_+}
\cT^o_{\a,n,\mm}e^{-i\lan\mm,\wt k\,\ran},\qqq
\cT^o_{\a,n,\mm}=\sum\limits_{\bc\in\wt\cC_n^\mm}\o(\bc)e^{-i\a(\cP\bc)}.
\]
Here $\cP:\wt\cC\to\cC^0$ is the projection of $\wt\cC$ onto $\cC^0$ defined by \er{depr}, the cycles $\bc_1,\ldots,\bc_d\in\cC$ satisfy the first condition in \er{bacC}, and $\t_+$ is given by \er{tapl}.

In particular, if the Betti number $\b=d$, then
$$
\cP\bc=0,\qqq \forall\,\bc\in\wt\cC,
$$
and
\[\lb{eqTr}
\Tr H_\a^n(\wt k+k_o)=\Tr H_0^n(\wt k\,), \qqq \forall\,n\in\N.
\]

ii) If the magnetic potentials $\a_1$ and $\a_2$ have the same flux through each of $\b-d$ basis cycles of $\cC^0$, then $\s(H_{\a_1})=\s(H_{\a_2})$.
\end{corollary}

\no {\bf Proof.} \emph{i}) Let $n\in\N$ and
$\wt k=k-k_o$, where $k_o$ is defined in \er{sfqm}. Then in the shifted  quasimomentum $\wt k$ the identities \er{TrnH}, \er{deTnk} have the form
$$
\Tr H_\a^n(k)=\Tr H_\a^n(\wt k+k_o)=\sum_{\bc\in\wt\cC_n}\o(\bc)\,e^{-i\f(\bc,\wt k+k_o)}=
\sum_{\bc\in\wt\cC_n}\o(\bc)\cos\f(\bc,\wt k+k_o),
$$
where
\[\lb{fpi0}
\f(\bc,\wt k+k_o)=\a(\bc)+\lan\t(\bc),\wt k+k_o\ran=
\lan\t(\bc),\wt k\,\ran+\a(\bc)+\lan\t(\bc),k_o\ran.
\]

Due to Lemma \ref{LbocC}, there exists a basis $\{\bc_1,\ldots,\bc_\b\}$ of the cycle space $\cC$ of the fundamental graph $\cG_*=(\cV_*,\cA_*)$ satisfying the conditions \er{bacC}, and $\wt\cB=\{\bc_1,\ldots,\bc_\b\}\cup\{\be_x\}_{x\in\cV_*}$ is the basis of the cycle space $\wt\cC$ of the modified fundamental graph $\wt\cG_*$. Then each cycle $\bc\in\wt\cC$ has the unique representation \er{cyre1} with respect to the basis $\wt\cB$. Recall that $\be_x$ is a loop at the vertex $x$ in $\wt\cG_*$ with index $\t(\be_x)=0$ and magnetic potential $\a(\be_x)=0$. Since $\{\bc_{d+1},\ldots,\bc_\b\}$ is a basis of $\cC^0$, then
\[\lb{zein}
\t(\bc_{d+1})=\ldots=\t(\bc_\b)=0.
\]
Since $\{\t(\bc_s)\}_{s=1}^d$ forms the standard orthonormal basis for $\Z^d$ and $k_o=-(\a(\bc_1),\ldots,\a(\bc_d))$, then
\[\lb{wkpm}
\a(\bc_s)+\lan\t(\bc_s),k_o\ran=\a(\bc_s)-\a(\bc_s)=0, \qqq \forall\,s\in\N_d.
\]
Substituting \er{cyre1} into \er{fpi0} and using the conditions \er{zein} and \er{wkpm}, for each cycle $\bc\in\wt\cC$ we obtain
\begin{multline*}
\textstyle\f(\bc,k)=\f(\bc,\wt k+k_o)=\lan \t(\bc),\wt k\,\ran+\sum\limits_{s=1}^\b n_s(\bc)\big(\a(\bc_s)+\lan \t(\bc_s),k_o\ran\big)\\\textstyle=\lan \t(\bc),\wt k\,\ran+\sum\limits_{s=d+1}^\b n_s(\bc)\a(\bc_s)=\lan \t(\bc),\wt k\,\ran+\a\Big(\sum\limits_{s=d+1}^\b n_s(\bc)\bc_s\Big)=\lan\t(\bc),\wt k\,\ran+\a(\cP\bc),
\end{multline*}
where the projection $\cP:\wt\cC\to\cC^0$ is defined by \er{depr}.
Thus, \er{wtl2.13} is proved. Rewriting the Fourier series \er{wtl2.13} in the standard form and using \er{cNe0}, we obtain \er{TrHM}.

If $\b=d$, then $\cP\bc=0$, and consequently $\a(\cP\bc)=0$ for any $\bc\in\wt\cC$. This and \er{TrHM} yield \er{eqTr}.

\emph{ii}) Let the magnetic potentials $\a_1$ and $\a_2$ have the same flux through each of $\b-d$ basis cycles of $\cC^0$. Then $\a_1$ and $\a_2$ have the same flux through each cycle of $\cC^0$. This and \er{wtl2.13} yield that
$$
\Tr H_{\a_1}^n(\wt k+k_{o1})=\Tr H_{\a_2}^n(\wt k+k_{o2}),\qqq \forall\,(n,\wt k\,)\in\N\ts\T^d,
$$
where $k_{os}=-(\a_s(\bc_1),\ldots,\a_s(\bc_d))$, $s=1,2$. Here we have also used that $\cP\bc\in\cC^0$ for all $\bc\in\wt\cC$.
Since the spectrum of a $\n\ts\n$ matrix $A$ is determined by $\Tr A^n$, $n=1,\ldots,\n$, we obtain
$$
\s(H_{\a_1}(\wt k+k_{o1}))=\s(H_{\a_2}(\wt k+k_{o2})),\qqq \forall\,\wt k\in\T^d.
$$
Then, by \er{specH}, $\s(H_{\a_1})=\s(H_{\a_2})$. \qq $\BBox$

\begin{remark}\lb{Re1}
\emph{i}) Corollary \ref{TCo0}.\emph{ii} means that the effect of the magnetic field on the spectrum of the magnetic Schr\"odinger operator $H_\a$ on a periodic graph $\cG$ is fully determined by only $\b-d$ fluxes through some independent cycles of the fundamental graph $\cG_*$ with zero index. Since each cycle with zero index in $\cG_*$ corresponds to some cycle in the periodic graph $\cG$ (see Remark \ref{Re22}), then \textbf{the effect of the magnetic field on the spectrum of $H_\a$ on $\cG$ is fully determined by $\b-d$ fluxes through some non-equivalent independent cycles of the periodic graph $\cG$}. Note that for a finite graph $\cG$, the spectrum of $H_\a$ is defined by $\b$ fluxes through basis cycles of $\cG$ (see \cite{LL93}). For an arbitrary locally finite graph $\cG$, the result "fluxes determine the spectrum of $H_\a$" was proved in \cite{CTT11}. More precisely, it was shown that if $\a_1$ and $\a_2$ have the same flux through each cycle of $\cG$ (which is equivalent to existence of a gauge transformation), then the operators $H_{\a_1}$ and $H_{\a_2}$ are unitarily equivalent.

\emph{ii}) It is known (see Theorem B in \cite{HS99b} or Corollary 2.2 in \cite{KS17}) that if the Betti number $\b$ of the fundamental graph $\cG_*$ coincides with the dimension $d$ of the periodic graph $\cG$, then for any magnetic potential $\a$ the magnetic Schr\"odinger operator $H_\a$ on $\cG$ is unitarily equivalent to the Schr\"odinger operator $H_0$ without magnetic fields, and, consequently, $\s(H_\a)=\s(H_0)$. Corollary \ref{TCo0}.\emph{ii} gives one more proof of this fact, using the trace formulas.
\end{remark}

The following example illustrates Remark \ref{Re1}.\emph{i}.

\begin{example}\lb{EKL}
Let $\D_\a$ be the Laplacian with a periodic magnetic potential $\a$ on the Kagome lattice $\bK$ shown in Fig.~\ref{FKaL}a. We assume that the potential $\a$ has the same periods $\ga_1,\ga_2$ as $\bK$. Then $\b-d=2$ and the spectrum of $\D_\a$ is fully determined by two fluxes $\f_s:=\a(\bc_s)$, $s=1,2$, through the cycles
\[\lb{cc12}
\bc_1=(\be_1,\be_2,\be_3), \qqq  \bc_2=(\be_4,\be_5,\be_6)
\]
on the Kagome lattice $\bK$, see Fig.~\ref{FKaL}a.
\end{example}

\no \textbf{Proof.} The fundamental graph $\bK_*=(\cV_*,\cE_*)$ of the Kagome lattice $\bK$  consists of three vertices $x_1,x_2,x_3$ and six edges $\be_1,\ldots,\be_6$, see Fig.~\ref{FKaL}\emph{b}. Then
$$
\b-d=\#\cE_*-\#\cV_*+1-d=6-3+1-2=2.
$$

The fluxes $\f_s:=\a(\bc_s)$, $s=1,2$, through the cycles $\bc_s$ (see \er{cc12}) on $\bK$ have the form
\[\lb{ff12}
\f_1=\a_1+\a_2+\a_3,\qqq \f_2=\a_4+\a_5+\a_6,
\]
where $\a_s=\a(\be_s)\in\R$, $s\in\N_6$, are the magnetic potentials on edges $\be_1,\ldots,\be_6$ of the fundamental graph $\bK_*$.

The Kagome lattice $\bK$ is a periodic regular graph of degree $\vk_+=4$. Thus, due to \er{suD}, instead of the Laplacian $\D_\a$ we can study the magnetic adjacency operator $A_\a$ given by \er{ALO}. Using \er{fnml0} and the edge indices \er{inKl}, we obtain that the fiber magnetic adjacency operator $A_\a(k)$, $k=(k_1,k_2)\in\T^2$, on $\bK_*$ is given by
\[\lb{effL}
A_\a(k)=\left(
\begin{array}{ccc}
0 & e^{-i\a_2}+e^{i(\a_6+k_2)} & e^{i\a_1}+e^{-i(\a_4-k_1)} \\
e^{i\a_2}+e^{-i(\a_6+k_2)} & 0 & e^{-i\a_3}+e^{i(\a_5+k_1-k_2)}\\
e^{-i\a_1}+e^{i(\a_4-k_1)} & e^{i\a_3}+e^{-i(\a_5+k_1-k_2)}  & 0
\end{array}\right).
\]
We introduce the shifted quasimomentum $\wt k=(\wt k_1,\wt k_2)\in\T^2$ by
$$
\wt k_1=k_1-\a_1-\a_4, \qqq \wt k_2=k_2+\a_2+\a_6.
$$
Then the fiber operator \er{effL} in the quasimomentum $\wt k$ has the form
\begin{multline*}
\wt A_\a(\wt k)=A_\a(\wt k_1+\a_1+\a_4,\wt k_2-\a_2-\a_6)\\=\left(
\begin{array}{ccc}
0 & e^{-i\a_2}+e^{i(\wt k_2-\a_2)} & e^{i\a_1}+e^{i(\wt k_1+\a_1)} \\
e^{i\a_2}+e^{-i(\wt k_2-\a_2)} & 0 & e^{-i\a_3}+e^{i(\wt k_1-\wt k_2+\f_2+\a_1+\a_2)}\\
e^{-i\a_1}+e^{-i(\wt k_1+\a_1)} & e^{i\a_3}+e^{-i(\wt k_1-\wt k_2+\f_2+\a_1+\a_2)}  & 0
\end{array}\right).
\end{multline*}
The characteristic polynomial of $\wt A_\a(\wt k)$ is given by
\begin{multline*}
\det\big(\l I_3-\wt A_\a(\wt k)\big)=\l^3-2\l\big(\cos(k_1-k_2+\f_1+\f_2)+\cos k_1+\cos k_2+3\big)\\-2\big(\cos(k_1-k_2+\f_1)+\cos(k_1-k_2+\f_2)+\cos(k_1+\f_1)+\cos(k_1+\f_2)\\
+\cos(k_2-\f_1)+\cos(k_2-\f_2)+\cos\f_1+\cos\f_2\big),
\end{multline*}
where $I_3$ is the $3\ts3$ identity matrix, and $\f_s$, $s=1,2$, are given by \er{ff12}. This yields that for each $\wt k\in\T^2$ the spectrum of $\wt A_\a(\wt k)$ depends only on the magnetic fluxes $\f_1$ and $\f_2$. Thus, the spectrum of $A_\a$ having the form
$$
\s(A_\a)=\bigcup_{\wt k\in\T^2}\s\big(\wt A_\a(\wt k)\big)
$$
also depends only on $\f_1$ and $\f_2$. \qq $\BBox$

\begin{remark}
The spectra of the magnetic Laplacian $\D_{\f_1,\f_2}:=\D_\a$ with the fluxes $\f_1,\f_2\in\{0,\pi\}$ on the Kagome lattice consist of three spectral bands including one flat band:
$$
\begin{array}{l}
\s(\D_{0,0})=[0,3]\cup[3,6]\cup\{6\},\\[4pt]
\s(\D_{\pi,\pi})=\{2\}\cup[2,5]\cup[5,8],\\[4pt]
\s(\D_{0,\pi})=\s(\D_{\pi,0})=[1,4]\cup\{4\}\cup[4,7].
\end{array}
$$
Note that $\s(\D_{0,0})$ and $\s(\D_{\pi,\pi})$ are symmetric to each other with respect to the point $\vk_+=4$, and $\s(\D_{0,\pi})=\s(\D_{\pi,0})$ is symmetric with respect to $\vk_+=4$.
\end{remark}

In the following statement we compare the traces of the fiber
magnetic Schr\"odinger operators and the traces of the corresponding fiber
Schr\"odinger operators without magnetic fields.

\begin{corollary}\lb{CCrFA}
Let $H_\a(k)$, $k\in\T^d$, be the fiber magnetic Schr\"odinger operator defined by \er{Hvt'} -- \er{fado}. Then for each $n\in\N$
\[\label{TrHn}
\Tr\big(H_0^n(k)-H_\a^n(k+k_o)\big)
=\sum_{\bc\in\wt\cC_n}\o(\bc)\cS(\bc,k), \qq \cS(\bc,k)=2\sin\frac{\a(\cP\bc)}2
\sin\Big(\frac{\a(\cP\bc)}2+\lan\t(\bc),k\ran\Big),
\]
\[\label{TrHnF}
\Tr\big(H_0^n(k)-H_\a^n(k+k_o)\big)
=\sum_{\mm\in\Z^d \atop \|\mm\|\leq n\t_+}\hspace{-3mm}\gt_{\a,n,\mm}e^{-i\lan\mm,k\ran},\qq \gt_{\a,n,\mm}=\sum_{\bc\in\wt\cC_n^\mm}\o(\bc)\big(1-e^{-i\a(\cP\bc)}\big),
\]
\[\label{TrHnI}
\frac1{(2\pi)^d}\int_{\T^d}\Tr\big(H_0^n(k)-H_\a^n(k)\big)dk
=\gt_{\a,n,0}, \qqq \gt_{\a,n,0}=2\sum_{\bc\in\wt\cC_n^0}\o(\bc)\sin^2\frac{\a(\bc)}2\,,
\]
where $k_o\in\T^d$ is defined in Corollary \ref{TCo0}; $\cP:\wt\cC\to \cC^0$ is the projection of $\wt\cC$ onto $\cC^0$ defined by \er{depr}, and $\t_+$ is given in \er{tapl}.  In particular,
\[\label{TrHA1}
\begin{aligned}
&\Tr\big(H_0(k)-H_\a(k+k_o)\big)=\sum_{\bc\in\cC_1}\cS(\bc,k),\\
&\Tr\big(H_0^2(k)-H_\a^2(k+k_o)\big)=2\sum_{\bc\in\cC_1}v_{x_\bc}\cS(\bc,k)+
\sum_{\bc\in\cC_2}\cS(\bc,k),
\end{aligned}
\]
\[\label{TrHA1I}
\begin{aligned}
&\frac{1}{(2\pi)^d}\int_{\T^d}\Tr\big(H_0(k)-H_\a(k)\big)dk=
2\sum_{\bc\in\cC_1^0}\sin^2\frac{\a(\bc)}2\geq0,
\\
&\frac1{(2\pi)^d}\int_{\T^d}\Tr\big(H_0^2(k)-H_\a^2(k)\big)dk=
4\sum_{\bc\in\cC_1^0}v_{x_\bc}\sin^2\frac{\a(\bc)}2+
2\sum_{\bc\in\cC_2^0}\sin^2\frac{\a(\bc)}2\,,
\end{aligned}
\]
where $x_\bc$ is the unique vertex of the cycle $\bc\in\cC_1$, and $v_x=V_x-\vk_x$, $x\in\cV_*$.
\end{corollary}

\no\textbf{Proof.} Let $n\in\N$. Due to \er{TrnH}, \er{deTnk} and \er{wtl2.13}, we have
\begin{multline}\label{difT}
\Tr\big(H_0^n(k)-H_\a^n(k+k_o)\big)=\cT_{0,n}(k)-\Tr H_\a^n(k+k_o)\\=
\sum_{\bc\in\wt\cC_n}\o(\bc)\big(\cos\lan\t(\bc),k\ran
-\cos\big(\a(\cP\bc)+\lan\t(\bc),k\ran\big)\big)\\
=2\sum_{\bc\in\wt\cC_n}\o(\bc)\sin\frac{\a(\cP\bc)}2
\sin\Big(\frac{\a(\cP\bc)}2+\lan\t(\bc),k\ran\Big).
\end{multline}
Thus, \er{TrHn} is proved. Similarly, using \er{TrnH}, \er{FsTrH} and \er{TrHM}, we obtain
\[\label{difT1}
\Tr\big(H_0^n(k)-H_\a^n(k+k_o)\big)=\cT_{0,n}(k)-\Tr H_\a^n(k+k_o)=
\sum_{\mm\in\Z^d\atop \|\mm\|\leq n\t_+}\big(\cT_{0,n,\mm}-
\cT^o_{\a,n,\mm}\big)e^{-i\lan\mm,k\ran}.
\]
Using the definitions of $\cT_{0,n,\mm}$ in \er{FsTrH} and $\cT^o_{\a,n,\mm}$ in \er{TrHM}, we can rewrite \er{difT1} in the form \er{TrHnF}.

Integrating \er{TrHn} over $k\in\T^d$ and using \er{inex} and the fact that $\a(\cP\bc)=\a(\bc)$ for each cycle $\bc\in\wt\cC$ with zero index, we get \er{TrHnI}. Formulas \er{TrHA1} and \er{TrHA1I} follow from \er{TrH123}, \er{TR13IH} and the second identity in \er{wtl2.13}. \qq \BBox

\section{Estimates of bandwidths for magnetic Schr\"odinger operators}\lb{Sec4}
\setcounter{equation}{0}
\subsection{Estimates of bandwidths for magnetic adjacency operators} First we estimate the total bandwidth $\gS(A_\a)=\sum_{j=1}^\n\big|\s_j(A_\a)\big|$ for the magnetic adjacency operator $A_\a$.

\begin{corollary}\lb{CAEs}
Let $A_\a$ be the adjacency operator defined by \er{ALO} with a periodic magnetic potential $\a$ on a periodic graph $\cG$. Then the total bandwidth $\gS(A_\a)$ satisfies
$$
\gS(A_\a)\geq\frac{\max\big\{\wt B_{\a,n}^+,2\wt B_{\a,n}^{odd}\big\}}{n\vk_+^{n-1}}\,, \qqq \forall\,n\in\N,
$$
where
$$
\wt B_{\a,n}^+=\Big|\sum_{\bc\in\cC_n^+}\cos\a(\bc)\Big|\,,\qqq
\wt B_{\a,n}^{odd}=\Big|\sum_{\bc\in\cC_n^{odd}}\cos\a(\bc)\Big|\,.
$$
Here $\cC_n^+$ and $\cC_n^{odd}$ are defined in \er{cNn+} and \er{cNno}, respectively; and $\vk_+$ is given in \er{spAD}.
\end{corollary}

We omit the proof, since it repeats the proof of Theorem \ref{TVECS}.

\subsection{Estimates of bandwidths for magnetic Schr\"odinger operators}
\lb{ssSO} We discuss estimates of the total bandwidth for the magnetic Schr\"odinger operators $H_\a=-\D_\a+V$. We can always shift the periodic potential $V$ by a constant such that
\[
\lb{Qmi0} \min_{x\in\cV_*}V_x=\vk_+,\qqq\textrm{where}\qqq \vk_+=\max_{x\in\cV_*}\vk_x,
\]
and $\vk_x$ is the degree of the vertex $x\in\cV_*$.
From the inclusion $\s(-\D_\a)\subseteq[-2\vk_+,0]$ and under the condition \er{Qmi0} we deduce that
\[\lb{splo}
\s(H_\a)=\s(-\D_\a+V)\subseteq[-2\vk_++\vk_+,\vk_++\diam V]=[-\vk_+,\vk_++\diam V],
\]
where $\diam V=\max_{x\in\cV_*}V_x-\min_{x\in\cV_*}V_x$.

\medskip

In order to prove Theorems \ref{TVECS} and \ref{TEsFK} about lower estimates of the total bandwidth for the magnetic Schr\"odinger operators we need the following lemma.

\begin{lemma}
\lb{Lesbw}
Let $H_\a=-\D_\a+V$ be the Schr\"odinger operator defined by \er{Sh} -- \er{ALO} with a periodic magnetic potential $\a$ and a periodic electric potential $V$ on a periodic graph $\cG$, and let $n\in\N$. Then the total bandwidth $\gS(H_\a^n)$ of $H_\a^n$ satisfies the following estimates
\[\lb{Le11}
\gS(H_\a^n)\leq n (\diam V+\vk_+)^{n-1}\gS(H_\a),
\]
\[\lb{Le12}
\gS(H_\a^n)\geq
\big|\Tr H_\a^n(k_1)-\Tr H_\a^n(k_2)\big|,\qqq \forall\, k_1,k_2\in\T^d,
\]
where $H_\a(k)$, $k\in\T^d$, is the fiber magnetic Schr\"odinger operator given by \er{Hvt'} -- \er{fado}.
\end{lemma}

\no \textbf{Proof.} Let $n\in\N$. By the definition of the spectral bands $\s_j(H_\a^n)$ and $\s_j(H_\a)$, $j=1,\ldots,\n$, see \er{ban.1H}, we get for some $k^\pm\in\T^d$:
\begin{multline}\lb{esba}
|\s_j(H_\a^n)|=\max_{k\in\T^d}\l_{\a,j}^n(k)-\min_{k\in\T^d}\l_{\a,j}^n(k)=
\l^n_{\a,j}(k^+)-\l_{\a,j}^n(k^-)\\\le n(\diam V+\vk_+)^{n-1}\big|\l_{\a,j}(k^+)-\l_{\a,j}(k^-)\big|\le n (\diam V+\vk_+)^{n-1}|\s_j(H_\a)|.
\end{multline}
Here we have also used the inclusion \er{splo}. Summing \er{esba} over $j=1,\ldots,\n$, we obtain \er{Le11}.

The estimate \er{Le12} also follows from the definition of the spectral bands of $H_\a^n$:
\begin{multline*}
\gS(H_\a^n)=\sum_{j=1}^\n\big|\s_j(H_\a^n)\big|=
\sum_{j=1}^\n\Big(\max_{k\in\T^d}\l_{\a,j}^n(k)-
\min_{k\in\T^d}\l_{\a,j}^n(k)\Big)\\\geq
\max_{k\in\T^d}\sum_{j=1}^\n\l_{\a,j}^n(k)-
\min_{k\in\T^d}\sum_{j=1}^\n\l_{\a,j}^n(k)
=\max_{k\in\T^d}\Tr H_\a^n(k)-\min_{k\in\T^d}\Tr H_\a^n(k)\\\geq
\big|\Tr H_\a^n(k_1)-\Tr H_\a^n(k_2)\big|,\qqq \forall\, k_1,k_2\in\T^d. \qqq \BBox
\end{multline*}

Now we prove Theorems \ref{TVECS} and \ref{TEsFK}.

\medskip

\no{\bf Proof of Theorem \ref{TVECS}}. Let $n\in\N$. Due to \er{Le12}, we have
\[\lb{ank}
\gS(H_\a^n)\geq
\big|\Tr H_\a^n(0)-\Tr H_\a^n(k)\big|,\qqq \forall\, k\in\T^d.
\]
By the mean value theorem,
$$
\frac1{(2\pi)^d}\int_{\T^d}\Tr H_\a^n(k)dk=\Tr H_\a^n(\hat k)\qqq \textrm{for some}\qqq \hat k\in\T^d.
$$
Then, using \er{TrnH}, \er{deTnk} and \er{ITrH}, \er{IVFu}, we write the estimate \er{ank} as $k=\hat k$ in the form
\begin{multline}\label{esB1}
\gS(H_\a^n)\geq\big|\Tr H_\a^n(0)-\Tr H_\a^n(\hat k)\big|=
\Big|\cT_{\a,n}(0)-\cT_{\a,n,0}\Big|
\\=\Big|\sum_{\bc\in\wt\cC_n}\o(\bc)\cos\a(\bc)
-\sum_{\bc\in\wt\cC_n^0}\o(\bc)\cos\a(\bc)\Big|=
\Big|\sum_{\bc\in\wt\cC_n^+}\o(\bc)\cos\a(\bc)\Big|.
\end{multline}
Similarly, if $k=\pi\1$, where $\1=(1,\ldots,1)\in\R^d$, then, using \er{TrnH}, \er{deTnk}, we obtain
\begin{multline}\lb{esHan}
\gS(H_\a^n)\geq
\big|\Tr H_\a^n(0)-\Tr H_\a^n(\pi\1)\big|=\Big|\sum_{\bc\in\wt\cC_n}\o(\bc)
\big(\cos\f(\bc,0)-\cos\f(\bc,\pi\1)\big)\Big|\\=
\Big|\sum_{\bc\in\wt\cC_n}\o(\bc)\big(\cos\a(\bc)-\cos\big(\a(\bc)+\lan \t(\bc),\pi\1\ran\big)\big)\Big|\\=
\Big|\sum_{\bc\in\wt\cC_n}\o(\bc)\big(\cos\a(\bc)-(-1)^{\lan \t(\bc),\1\ran}\cos\a(\bc)\big)\Big|=
2\,\Big|\sum_{\bc\in\wt\cC_n^{odd}}\o(\bc)\cos\a(\bc)\Big|.
\end{multline}
Here we have also used that $\t(\bc)\in\Z^d$ for all cycles
$\bc\in\wt\cC$.  Combining \er{esB1} and \er{esHan}, we get
$$
\gS(H_\a^n)\geq\max\big\{B_{\a,n}^+,2B_{\a,n}^{odd}\big\},
$$
where $B_{\a,n}^+,B_{\a,n}^{odd}$ are given by \er{Banpo}. This and \er{Le11} yield \er{leHa}.

Now let $n$ be the length of the shortest cycle with non-zero index in the fundamental graph $\cG_*$. Then each cycle $\bc$ of length $n$ and with non-zero index in the modified fundamental graph $\wt\cG_*$ is also a cycle in $\cG_*$ and, due to \er{oc1C}, $\o(\bc)=1$. Thus, the constants $B_{\a,n}^+$ and $B_{\a,n}^{odd}$ defined by \er{Banpo} have the form \er{Banp1}. \qq \BBox

\medskip

\no \textbf{Proof of Theorem \ref{TEsFK}.} Let $n\in\N$.
For the Fourier coefficients $\cT_{\a,n,\mm}$ of the function $\Tr H_\a^n(k)$, see \er{TrnH}, \er{FsTrH}, we have
\[\lb{FC1}
\cT_{\a,n,\mm}=\frac1{(2\pi)^d}\int_{\T^d}e^{\,i\lan\mm,k\,\ran}\Tr H_\a^n(k)dk,\qqq \mm\in\Z^d.
\]
Let $\mm=(m_s)_{s=1}^d\in\Z^d\sm\{0\}$. Then $m_s\neq0$ for some $s\in\N_d$. Replacing $k$ by $k+\frac\pi{m_s}\mathfrak{e}_s$, where $\{\mathfrak{e}_s\}_{s=1}^d$ is the standard orthonormal basis of $\R^d$, we can write
\begin{multline}\lb{FC2}
\cT_{\a,n,\mm}=\frac1{(2\pi)^d}\int_{\T^d}e^{\,i\lan\mm, k+\frac\pi{m_s}\mathfrak{e}_s\,\ran}\Tr H_\a^n\big(k+{\textstyle\frac\pi{m_s}}\mathfrak{e}_s\big)dk\\=
-\frac1{(2\pi)^d}\int_{\T^d}e^{\,i\lan\mm,k\,\ran}\Tr H_\a^n\big(k+{\textstyle\frac\pi{m_s}}\mathfrak{e}_s\big)\,dk.
\end{multline}
Summing \er{FC1} and \er{FC2}, we obtain
$$
2\cT_{\a,n,\mm}=\frac1{(2\pi)^d}\int_{\T^d}e^{\,i\lan\mm,k\,\ran}\Big(\Tr H_\a^n(k)-\Tr H_\a^n\big(k+{\textstyle\frac\pi{m_s}}\mathfrak{e}_s\big)\Big)\,dk.
$$
This and \er{Le12} yield
$$
2\big|\cT_{\a,n,\mm}\big|\leq\frac1{(2\pi)^d}\int_{\T^d}\Big|\Tr H_\a^n(k)-\Tr H_\a^n\big(k+{\textstyle\frac\pi{m_s}}\mathfrak{e}_s\big)\Big|\,dk\leq\gS(H_\a^n).
$$
Finally, using \er{Le11},  we obtain \er{EsFS}.

Now let $\mm\in\Z^d\sm\{0\}$ and $\frac1q\,\mm\not\in\Z^d$ for any integer $q\geq2$, and let $n(\mm)$ be the length of the shortest cycle with index $\mm$ in the fundamental graph $\cG_*$. Then each cycle $\bc$ from $\wt\cC_{n(\mm)}^\mm$ is a cycle in the fundamental graph $\cG_*$, i.e., $\bc\in\cC_{n(\mm)}^\mm$ and, due to \er{oc1C}, $\o(\bc)=1$. This cycle $\bc\in\cC_{n(\mm)}^\mm$ is a prime cycle, i.e., it is not obtained by repeating $q$ times ($q\geq2$) any cycle in $\cG_*$. All cyclic permutations of edges of a prime cycle $\bc$ give $|\bc|$ distinct prime cycles having the same length $|\bc|$, magnetic flux $\a(\bc)$ and index $\t(\bc)$. Denoting by  $\bc_1,\ldots,\bc_p$ all cycles of length $n(\mm)$ and with index $\mm$ in $\cG_*$ (up to cyclic permutations of edges) and using \er{FsTrH}, we obtain
$$
\big|\cT_{\a,n(\mm),\mm}\big|=
\Big|\sum_{\bc\in\cC_{n(\mm)}^\mm}e^{-i\a(\bc)}\Big|=
n(\mm)\Big|\sum_{j=1}^pe^{-i\a(\bc_j)}\Big|=
n(\mm)\Big|1+\sum_{j=2}^pe^{-i\a(\bc_j-\bc_1)}\Big|\,.
$$
Thus, \er{EsFS1} is proved.

\emph{i}) If $p=1$, then \er{EsFS1} yields
$$
\big|\cT_{\a,n(\mm),\mm}\big|=n(\mm)>0.
$$

\emph{ii}) If $p=2$, then
$$
\big|\cT_{\a,n(\mm),\mm}\big|=
n(\mm)\Big|1+e^{-i\a(\bc_2-\bc_1)}\Big|=2n(\mm)
\Big|\cos\frac{\a(\bc_2-\bc_1)}2\Big|\,.
$$

\emph{iii}) If $p\geq2$ and $\a_+=\max_{2\leq j\leq p}\big|\a(\bc_j-\bc_1)\big|<\pi/2$, then
\begin{multline*}
\big|\cT_{\a,n(\mm),\mm}\big|=n(\mm)\Big|1+\sum_{j=2}^pe^{-i\a(\bc_j-\bc_1)}\Big|\geq n(\mm)\Big|1+\sum_{j=2}^p\cos\a(\bc_j-\bc_1)\Big|\\\geq n(\mm)(1+(p-1)\cos\a_+)>0.
\end{multline*}
Thus, in the cases \emph{i}) and \emph{iii}) and in the case \emph{ii}) when $\a(\bc_2-\bc_1)\neq\pi$ we have $\big|\cT_{\a,n(\mm),\mm}\big|>0$. Then, by \er{EsFS} as $n=n(\mm)$,
$$
\gS(H_\a)\geq \frac{2\big|\cT_{\a,n(\mm),\mm}\big|}{n(\mm)v_*^{n(\mm)-1}}>0.
$$
This yields that $H_\a$ has at least one non-degenerate band. \qq \BBox

\medskip

We prove Example \ref{Exa} about estimates of the total bandwidth of the magnetic Laplacian on the $\Z$-periodic graph $\cG$ shown in Fig. \ref{FEx1}\emph{a}.

\medskip

\no \textbf{Proof of Example \ref{Exa}.} \emph{i--ii})  The fundamental graph $\cG_*$ of the $\Z$-periodic graph $\cG$ shown in Fig.~\ref{FEx1}\emph{a} consists of four vertices $x_1,x_2,x_3,x_4$ with degrees $\vk_{x_1}=\vk_{x_4}=3$ and $\vk_{x_2}=\vk_{x_3}=2$ and five edges
$$
\be_1=(x_4,x_1),\qq \be_2=(x_1,x_3),\qq \be_3=(x_3,x_4),\qq
\be_4=(x_4,x_2), \qq \be_5=(x_2,x_1)
$$
with indices
\[\lb{edinE}
\t(\be_1)=1,\qqq \t(\be_2)=\t(\be_3)=\t(\be_4)=\t(\be_5)=0
\]
and their inverse edges (see Fig. \ref{FEx1}\emph{b}).

It is known that two magnetic Laplacians with the same fluxes through all cycles of a graph are unitarily equivalent (see Proposition 2.1.\emph{ii} in \cite{CTT11}). Thus, without loss of generality, we may assume that the periodic magnetic potential $\alpha$ on the periodic graph $\mathcal{G}$ is supported only on the edge $\mathbf{e}_5$ (and its equivalent edges), i.e.,
\begin{equation}
\label{Emp}
\alpha(\mathbf{e}_1)=\ldots=\alpha(\mathbf{e}_4)=0,\qquad \alpha(\mathbf{e}_5)=\phi,
\end{equation}
for some $\phi\in[0,2\pi)$.

We estimate the total bandwidth for the magnetic Laplacian $\D_\a$ on $\cG$ using Theorems \ref{TVECS} and \ref{TEsFK}. The length of the shortest cycle with non-zero index in $\cG_*$ is $3$. All cycles of length 3 in $\cG_*$ are the cycles
$$
\bc_1=(\be_1,\be_2,\be_3),\qqq \bc_2=(\be_1,\ul\be_5,\ul\be_4),
$$
their cyclic edge permutations (three permutations for each cycle $\bc_s$, $s=1,2$) and their reverse cycles. Using the definitions (\ref{cyin}) and (\ref{cyin2}) of cycle indices and magnetic fluxes and formulas (\ref{edinE}), (\ref{Emp}), we obtain
$$
\tau(\mathbf{c}_1)=\tau(\mathbf{c}_2)=1,\qquad \alpha(\mathbf{c}_1)=0,\qquad \alpha(\mathbf{c}_2)=-\phi.
$$
Then the constants $B_{\a,3}^+,B_{\a,3}^{odd}$ defined by \er{Banp1} when $n=3$ have the form
\[\lb{coBB}
\begin{array}{l}
B_{\a,3}^+=\Big|\sum\limits_{\bc\in\cC_3^+}\cos\a(\bc)\Big|=
6\big|\cos\a(\bc_1)+\cos\a(\bc_2)\big|=
6\big(1+\cos\phi\big)=12\cos^2\dfrac{\phi}2\,,\\[6pt] B_{\a,3}^{odd}=\Big|\sum\limits_{\bc\in\cC_3^{odd}}\cos\a(\bc)\Big|=
B_{\a,3}^+=12\cos^2\dfrac{\phi}2\,.
\end{array}
\]
Here we have also used that $\a(\ul\bc\,)=-\a(\bc)$ for all $\bc\in\cC$. By Theorem \ref{TVECS} when $n=3$, the total bandwidth for $\D_\a$ on $\cG$ satisfies
\[\lb{es1D}
\gS(\D_\a)\geq\frac{2B_{\a,3}^{odd}}{3\vk_+^{3-1}}=\frac{24\cos^2\frac{\phi}2}{3\cdot3^2}
=\frac89\,\cos^2\frac{\phi}2\,.
\]
Since the number of all cycles of length $3$ with index 1 in $\cG_*$ (up to the cyclic permutations) is $p=2$, then, by Theorem \ref{TEsFK}.\emph{ii} when $\mm=1$ and $n(\mm)=3$,  we obtain
$$
\big|\cT_{\a,3,1}\big|=2\cdot3\cdot\Big|\cos\frac{\a(\bc_2-\bc_1)}2\Big|
=6\,\Big|\cos\frac{\a(\ul\bc_0)}2\Big|=6\,\Big|\cos\frac{\phi}2\Big|,
$$
where $\bc_0=(\be_2,\be_3,\be_4,\be_5)$ is the cycle of the fundamental graph $\cG_*$ with flux $\a(\bc_0)=\phi$.

Then, due to \er{EsFS}, the total bandwidth for $\D_\a$ on $\cG$ satisfies
\[\lb{es1D'}
\gS(\D_\a)\geq\frac{2\big|\cT_{\a,3,1}\big|}{3\vk_+^{3-1}}=
\frac{4}{9}\,\Big|\cos\frac{\phi}2\Big|.
\]
Combining \er{es1D} and \er{es1D'}, we obtain the lower bound in \er{esti}. For the graph $\cG$, the invariant $\cI$ given by \er{dIm} is equal to one (see \er{edinE}). This and the estimate \er{ues1} yield the upper bound in \er{esti}.

The fiber magnetic Laplacian $\D_\a(k)$, $k\in[0,2\pi)$, defined by \er{Hvt'}, \er{fado} on the fundamental graph $\cG_*=(\cV_*,\cA_*)$ in the standard orthonormal basis of $\ell^2(\cV_*)$ is given by the $4\ts4$ matrix
\[\lb{efHSL1}
\D_\a(k)=\diag(3,2,2,3)-A_\a(k),\qqq A_\a(k)=\left(
\begin{array}{cccc}
0 & e^{i\phi} & 1 & e^{ik} \\
e^{-i\phi} & 0 & 0 & 1 \\
1 & 0 & 0 & 1 \\
e^{-ik} & 1 & 1 & 0 \\
\end{array}\right).
\]
The characteristic polynomial of $\D_\a(k)$ has the form
\begin{multline*}
\det(\l I_4-\D_\a(k))=\l^4-10\l^3+32\l^2+2\l(\cos k+\cos(k-\phi)-18)\\-4\cos k-4\cos(k-\phi)-2\cos\phi+10,
\end{multline*}
where $I_4$ is the $4\ts4$ identity matrix.

If the magnetic flux $\phi=\pi$, then the eigenvalues of $\D_\a(k)$ do not depend on $k$ and are given by $2\pm\sqrt{2}\,,3\pm\sqrt{3}$. Thus, the spectrum of $\D_\a$ degenerates to four eigenvalues with infinity multiplicity, and $\gS(\D_\a)=0$. The lower bound in \er{esti} as $\phi=\pi$ is also zero. Thus, the lower estimate in \er{esti} is sharp.

Let $\phi=0$. Then $\a=0$. It is known \cite{EKW10} that the spectral bands of the Laplacian $\D_0$ without magnetic fields on a $\Z$-periodic graph with the invariant $\cI=1$ do not overlap and their endpoints are the eigenvalues of $\D_0(0)$ and $\D_0(\pi)$. By direct calculations we obtain
$$
\s\big(\D_0(0)\big)=\{0,2,4,4\},\qqq \s\big(\D_0(\pi)\big)=\{3-\sqrt{5}\;,2,2,3+\sqrt{5}\;\}.
$$
Thus, the spectral bands of $\D_0$ are given by
$$
\s_1(\D_0)=[0,3-\sqrt{5}\;],\qq \s_2(\D_0)=\{2\},\qq
\s_3(\D_0)=[2,4],\qq \s_4(\D_0)=[4,3+\sqrt{5}\;],
$$
and $\gS(\D_0)=\sum_{j=1}^{4}|\s_j(\D_0)|=4$. Therefore, the upper estimate in \er{esti} is also sharp. \qq \BBox

\medskip

In order to prove Proposition \ref{Tfmp} about the flat band spectrum of the magnetic Schr\"odinger operators we need the following simple lemma.

\begin{lemma}
\lb{Lcond}
Let $H_\a=-\D_\a+V$ be the Schr\"odinger operator defined by \er{Sh} -- \er{ALO} with  a periodic magnetic potential $\a$ and a periodic electric potential $V$ on a periodic graph $\cG$. Then all spectral bands of $H_\a$ are degenerate iff
\[\lb{cTe0}
\cT_{\a,n,\mm}=0\qq\textrm{ for all } \qq (n,\mm)\in\N_\n\ts\big(\Z^d\sm\{0\}\big),
\]
where $\cT_{\a,n,\mm}$ are the Fourier coefficients defined in \er{FsTrH}.
\end{lemma}

\no \textbf{Proof.} The characteristic polynomial of the fiber operator $H_\a(k)$ defined by \er{Hvt'} -- \er{fado} has the form
\[\lb{Det}
\det\big(\l I_\n-H_\a(k)\big)=\prod_{j=1}^\n(\l-\l_{\a,j}(k))=
\l^\n+\x_1\l^{\n-1}+\x_2\l^{\n-2}+\ldots+\x_{\n-1}\l+\x_\n,
\]
where the coefficients $\x_j$ are given by (see, e.g., p. 87--88 in \cite{Ga60})
\[\lb{coDe}
\x_n=-{1\/n}\rt(T_n+\sum_{j=1}^{n-1}T_{n-j}\,\x_j \rt), \qqq T_n=\Tr
H^n_\a(k),\qqq n\in \N_\n,
\]
and, in particular, $\x_1=-T_1$, $\x_2=-{1\/2}(T_2-T_1^2),\,\ldots\,$.

Let all spectral bands of $H_\a$ be degenerate. Then all band functions $\l_{\a,j}(\cdot)$, $j\in\N_\n$, are constant, and, consequently, $\Tr H_\a^n(k)=\sum_{j=1}^\n\l_{\a,j}^n(k)$ does not depend on $k$ for each $n\in\N_\n$. Then, using the Fourier series
\er{TrnH}, \er{FsTrH}, we obtain \er{cTe0}.

Conversely, let the condition \er{cTe0} hold true.  Then, using \er{TrnH}, \er{FsTrH}, we obtain
$$
\Tr H_\a^n(k)=\cT_{\a,n,0}=\sum_{\bc\in\wt\cC_n^0}\o(\bc)e^{-i\a(\bc)},
\qqq \forall\,n\in\N_\n,
$$
i.e., the traces $\Tr H_\a^n(k)$, $n\in\N_\n$, do not depend on $k$.
From this  and \er{Det}, \er{coDe} it follows that the
determinant $\det\big(\l I_\n-H_\a(k)\big)$ does not depend on $k$.
Then all band functions $\l_{\a,j}(\cdot)$, $j\in\N_\n$, are
constant and all spectral bands $\s_j(H_\a)$ are degenerate. \qq \BBox

\medskip

\no \textbf{Proof of Proposition \ref{Tfmp}.} The first spectral band
$\s_1(H_0)$  of the Schr\"odinger operator $H_0$ without magnetic
field is non-degenerate (see Theorem 2.1.\emph{ii} in \cite{KS19}).
Then, by Lemma \ref{Lcond}, there exists $(n,\mm)\in\N_\n\ts\big(\Z^d\sm\{0\}\big)$ such that the Fourier coefficient
$\cT_{0,n,\mm}\neq0$, where $\cT_{0,n,\mm}$ is defined by \er{FsTrH}. For this $(n,\mm)$ and for the magnetic potential $t\a$ we define the function
$$
f(t):=\cT_{t\a,n,\mm}=\sum_{\bc\in\wt\cC_n^\mm}\o(\bc)e^{-it\a(\bc)}, \qq t\in\R,
$$
and note that the sum is finite. This function has an analytic extension to the whole complex plane. If $f=\const$, then we obtain
$f=f(0)=\cT_{0,n,\mm}\ne0$ for such specific $(n,\mm)$. Then Lemma
\ref{Lcond} yields that there exists at least one non-degenerate band, i.e., the absolutely continuous spectrum of the operator $H_{t\a}$ is not empty for all $t\in \R$.

If $f\ne \const$, then $f$ has only a finite number of zeros on any bounded interval. Then Lemma \ref{Lcond} yields that the absolutely continuous spectrum of  $H_{t\a}$ is not empty for all except finitely many $t\in[0,1]$.
\qq \BBox

\medskip

\footnotesize
 \textbf{Acknowledgments.}  Evgeny Korotyaev was
supported by the RSF grant  No. 23-21-00023. We would like to thank a referee for thoughtful comments that helped us to improve the manuscript.

\end{document}